\documentstyle[a4,mypic,amssymb,epsf,twoside,amstex,epic,eepic]{article}
\input psfig.tex

\def\proof#1{\noindent{\bf Proof #1:}\enspace}

\makeatletter

\@ifundefined{eqref}{
\def\rom#1{\leavevmode\skip@\lastskip\unskip\/%
  \ifdim\skip@=\z@\else\hskip\skip@\fi
  {\normalshape#1}}
\def\eqref#1{\rom{\tagform@{\ref{#1}}}}
}{}

\def\myitem{\@ifnextchar[{\@myitem}{\@myitem[3cm]}}
\def\@myitem[#1]#2{\noindent\hangafter1\hangindent#1\hbox to #1{#2\hfil}\kern0pt}

\def\myitemm#1{{\setbox0=\hbox{#1}\noindent\global
\hangafter1\global\hangindent\wd0\box0\kern0pt}\ignorespaces}

\def\ta{\raisebox{-0.8ex}{\~{}}}

\def\prtag{\tag*{\llap{$\Box$\hskip-\displaywidth}}}

\let\ot\otimes
\def\standardunitlength{0.08em}

\def\pentagon@pic{\begin{picture}(20,35)(0,-10)
\path(5,0)(15,0)(20,10)
	(10,18)(0,10)(5,0)
\end{picture}}

\def\hexagon@pic{\begin{picture}(18,35)(0,-10)
\path(9,20)(0,15)(0,5)
	(9,0)(18,5)(18,15)(9,20)
\end{picture}}

\def\silenteepic#1#2{%
	 \relax\unitlength\standardunitlength\relax
	 \unitlength#2\unitlength\relax
	\begin{array}{c}  \hspace{-0.5em}%
        	\raisebox{-0.3em}{\csname #1@pic\endcsname }%
        	\hspace{-0.4em}%
	\end{array}%
}

\def\pentagon{\silenteepic{pentagon}{0.4}}      % A pentagon
\def\hexagon{\silenteepic{hexagon}{0.4}}        % A hexagon
\def\hexagons{\hexagon\null_\pm}                % Both hexagons

\@addtoreset {equation}{section}\expandafter \xdef
\csname theequation\endcsname
{\expandafter \noexpand \csname thesection\endcsname
\@thmcountersep \@thmcounter {equation}}

\newlength{\weite}
\def\vcbox#1{{\setbox1=\vbox{#1}\parbox{\wd1}{\box1}}}
\def\vctext#1#2#3{\vcbox{\hbox{#3\begin{tabular}{#1}#2\end{tabular}}}}

\newcommand{\brac}[6]{%picture-x picture-y unitlength picture-text
{\unitlength#3\setlength{\weite}{#2\unitlength}
\addtolength{\weite}{1ex}
\left#5\fbox{\hbox to #1\unitlength{\vbox to \weite{\vss
\mbox{\rx{0.2ex}%
\begin{picture}(#1,#2)
#4
\end{picture}
\rx{0.2ex}}%
\vss}}}
\right#6
}}

\newdimen\mydim

\def\placetext#1#2#3#4#5{%
\mbox{{\setbox0=\hbox{#5}%
\mydim=\ht0\advance\mydim by #3\advance\mydim by #4%
\rx{#1}\zker\parbox{\wd0}{\rule[-#3]{0pt}{\mydim}\zker\box0}\zker\rx{#2}}}}
\def\diag#1#2#3#4{%
\placetext{-0.3ex}{-0.3ex}{0.2ex}{0.2ex}{
\unitlength#1
\begin{picture}(#2,#3)
\picinit
#4
\pictail
\end{picture}
}}

\def\snakest{\pic@pss{pt}}

\def\snakeend{\diag{1cm}{0}{0}{
\pic@grend
\pss{pt 2 c 2{5 4 r x - 10 : x}R
[/pic@unity
/pic@unitx
/pic@originy
/pic@originx]{ed}fa
}
\wdtoPSunit{0.15ex}
\pic@pss{pic@unitx :
d 5 * /pic@veclength ed
d 2 * /pic@vecwidth ed
d /pic@linew ed
3 *}
\picPSgraphics{cm /pic@matrix ed
}
\picstroke{
\piclinedash{3 2 r d}{0}
\picarcangle{10 10}{10.9 10}{10.9 8.6}{0.3}
\picarcangleto{10.9 8.1}{8.5 8.1}{0.3}
\picarcangleto{7.6 8.1}{7.6 7}{0.3}
\picarcangleto{7.6 5}{3.8 5}{0.3}
\picarcangleto{3.4 5}{3.4 2.2}{0.3}
\picarcangleto{3.4 1.7}{-0.5 1.7}{0.3}
\picarcangleto{-0.9 1.7}{-0.9 0.3}{0.3}
\picarcangleto{-0.9 -0.1}{-0.1 -0.1}{0.3}
}
\piclinedash{}{0}
\picvecrline{-0.1 -0.1}{0.1 0}
\pic@grst
}
}

\def\sumint{\setbox0=\hbox{$\displaystyle\sum$}\mathop{\rlap{\copy0}
\kern0pt \hbox to \wd 0{\hss$\displaystyle\int$\hss}}}

\def\faktor#1#2{%
\raisebox{0.6ex}{\mbox{$#1\mbox{\raisebox{-0.5ex}{\bigg /}}%
\mbox{\raisebox{-1.2ex}{$#2$}}$}}}

\def\rightarrowto#1{\hbox to #1{\rightarrowfill}}

\def\text#1{{\ifmmode{\mathchoice{\mbox{\normalsize #1}}{%
\mbox{\normalsize #1}}{\mbox{\scriptsize #1}}{\mbox{\tiny #1}}}\else\hbox{#1}\fi}}

\def\restrsize{\small}
\def\restrto#1{\kern0.2ex\raisebox{-0.6em}{\parbox[b]{0.1ex}{\raisebox{-0.2em}%
{\rule{0.1ex}{1.6em}}}\kern0.2ex\restrsize\mbox{$#1$}}}
\begin{document}

\newenvironment{DrorRemark}{\begingroup\small\sf\leftskip2\parindent\relax{\bf
Dror's Remark: }}{%
\vskip\parskip\endgroup}

\newcommand{\Span}{\text{\rm span}}

\newcommand{\psdraw}[2]
        {\begin{array}{c} \hspace{-1.3mm}
        \raisebox{-4pt}{\psfig{figure=draws/#1.ps,width=#2}}
        \hspace{-1.9mm}\end{array}}

\newcommand{\mathmode}[1]{$#1$}
\newlength{\eepiclength}
\setlength{\eepiclength}{0.00083333in}
\newcommand{\eepic}[2]{
        \setlength{\unitlength}{#2\eepiclength}
        \begin{array}{c}  \hspace{-1.7mm}
                \raisebox{-8pt}{\input draws/#1.tex }
                \hspace{-1.9mm}
        \end{array}
}

\def\qed{{\hfill\text{$\Box$}}}

\newenvironment{myitemize}{
        \begin{list}{$\bullet$}{\setlength{\leftmargin}{16pt}
        \setlength{\labelwidth}{12pt}
        \setlength{\labelsep}{4pt}}
}{
        \end{list}
}

\catcode`\@=11
\long\def\@makecaption#1#2{%
    \vskip 10pt
    \setbox\@tempboxa\hbox{%\ifvoid\tinybox\else\box\tinybox\fi
      \small\sf{\bfcaptionfont #1. }\ignorespaces #2}%
    \ifdim \wd\@tempboxa >\captionwidth {%
        \rightskip=\@captionmargin\leftskip=\@captionmargin
        \unhbox\@tempboxa\par}%
      \else
        \hbox to\hsize{\hfil\box\@tempboxa\hfil}%
    \fi}
\font\bfcaptionfont=cmssbx10 scaled \magstephalf
\newdimen\@captionmargin\@captionmargin=2\parindent
\newdimen\captionwidth\captionwidth=\hsize
\catcode`\@=12

\newtheorem{theo}{Theorem}[section]
\newtheorem{conj}[theo]{Conjecture}
\newtheorem{lemma}[theo]{Lemma}
\newtheorem{prop}[theo]{Proposition}
\newtheorem{propdef}[theo]{Proposition-Definition}
\newtheorem{propproof}[theo]{Proposition \& Proof}

\newtheorem{defi}[theo]{Definition}
\newtheorem{problem}[theo]{Problem}
\newtheorem{ques}[theo]{Question}

\newtheorem{example}[theo]{Example}
\newtheorem{exercise}[theo]{Exercise}
\newtheorem{hint}[theo]{Hint}
\newtheorem{remark}[theo]{Remark}

\def\nup{n\!\!\uparrow}
\def\ndown{n\!\!\downarrow}

\title{\underline{\Large\bf The Fundamental Theorem of Vassiliev
Invariants}}

\author{%
\small lecture notes by\\
\normalsize\sc Dror Bar-Natan,\\
\small The Hebrew University, Jerusalem\\
\null\\
\normalsize Odense, July 1995\\
\null\\
\small revised and prepared for publication by\\
\normalsize\sc Alexander Stoimenow,\\
\small Humboldt University, Berlin
}

\def\lb{\linebreak[0]}
\date{
  \parbox{4.8in}{\small
    These notes appeared in {\em Geometry and Physics} (J.~.E.~Andersen,
    J.~Dupont, H.~Pedersen, and A.~Swann, eds.), lecture notes in pure
    and applied mathematics {\bf 184}, Marcel Dekker, New-York 1997. They
    are also available electronically at
    {\tt http:\lb //\lb www.\lb ma.\lb huji.\lb ac.\lb il/\lb \ta drorbn},
    {\tt http:\lb //\lb www.\lb informatik.\lb hu-berlin.\lb de/\lb
    \ta stoimeno}, at
    {\tt ftp:\lb //\lb ftp.\lb ma.\lb huji.\lb ac.\lb il/\lb drorbn},
    and at
    {\tt http:\lb //\lb xxx.\lb lanl.\lb gov/\lb abs/\lb q-alg/\lb 9702009}.
  \\ \ } \\
  This edition: Feb.~06,~1997; \ \ First edition: Jan.~1,~1996.
}

\maketitle
\makeatletter
\sloppy

\def\lectitle#1{{\section{#1}
\markboth{The Fundamental Theorem of {\/\sc Vassiliev} Invariants}%
{{\rm Lecture \thesection: #1}}
}}

\tableofcontents

\section*{Introduction \normalsize (by the first author)}

These notes grew out of four lectures I gave in a summer school
titled ``Geometry and Physics'' in Odense, Denmark, in July 1995.
I had two purposes in giving these lectures. The first was to expose
the students to the theory of Vassiliev invariants and to some of its
numerous connections with other parts of mathematics and mathematical
physics. I chose to concentrate on only one theorem, the basic existence
theorem for invariants with a given ``$m$th derivative'' (which I call
``The Fundamental Theorem'' both for its fundamental nature and for its
similarity with the fundamental theorem of calculus). Each lecture was
a brief exposition of one of the four approaches I know for proving the
theorem, with each approach related to a different branch of mathematics.

My second purpose in giving these lectures was to draw attention to
the fact that even though the Fundamental Theorem is fundamental
and is proven, we still don't know the ``right'' proof. The naive
and most natural topological approach discussed in the first lecture
is not yet complete, and the slightly stronger theorem it requires
(conjecture~\ref{conjthree}) may well be false. Each of the other three
approaches does succeed, but always at some cost. Always the method is
indirect and very complicated, and/or some a-priori unnatural choices
have to be made, and/or the ground ring has to be limited. It seems like a
conspiracy, and I hope that it really is a conspiracy.  Maybe some small
perturbation(s) of the theorem is(are) false? Light travels on straight
lines, but not near very heavy objects. Maybe there's some heavy object
around here too, that prevents us from finding a direct proof? I hope
that that object will be found one day. It may be fertile. Is it near
conjecture~\ref{H1Vanishes}?

As it's often the case with lecture notes, these notes are not quite
perfectly organized, and many of the details are insufficiently explained.
I do hope, though, that they are clear enough at least to whet the reader's
appetite to read some of the references scattered within. The only new
mathematics in these notes is the repackaging of Hutchings' argument in
terms of the snake lemma in section~\ref{Hutchings}.

\noindent {\bf Acknowledgement.} We wish to thank J.~Andersen,
H.~Munkholm, H.~Pedersen, and A.~Swann, the organizers of the Odense
summer school, for caring for all our special needs (especially the first
author's), for feeding us good food, for the T-shirts, and for bringing us
(and all the others) together for a very enjoyable and productive period
of time in Denmark.

\lectitle{Topology (and Combinatorics)}

\subsection{Vassiliev invariants and the Fundamental Theorem}

Any invariant $V$ of oriented knots in oriented space can be extended
to an invariant of singular knots (allowing finitely many transverse
double points as singularities) by inductive use of the formula\footnote{
Here and throughout these notes we use the standard convention in knot
theory, that if several almost equal knots (or singular knots) appear in
an equation, only the parts in which they differ are drawn.
}:

\[
V\left(\diag{6mm}{1}{1}{
\picvecline{0 0}{1 1}
\picvecline{1 0}{0 1}
\picfillgraycol{0}
\picfilledcircle{0.5 0.5}{0.06}{}
% \piccircle{0.5 0.5}{0.76}{}
}\right)\,:=\,
V\left(\diag{6mm}{1}{1}{
\picmultivecline{-4 1 -1 0}{1 0}{0 1}
\picmultivecline{-4 1 -1 0}{0 0}{1 1}
}\right)\,-\,
V\left(\diag{6mm}{1}{1}{
\picmultivecline{-4 1 -1 0}{0 0}{1 1}
\picmultivecline{-4 1 -1 0}{1 0}{0 1}
}\right)\qquad\qquad\text{(verify consistency)}.
\]

Differences are cousins of derivatives, and it is tempting to think of $V$
evaluated on an $m$-singular knot (a knot with exactly $m$ double points)
as ``the $m$-th derivative of the original $V$''. In analogy with polynomials
of degree $m$ we define:

\begin{defi} (Goussarov~\cite{Goussarov:New, Goussarov:nEquivalence},
Vassiliev~\cite{Vassiliev:CohKnot, Vassiliev:Book})
$V$ is called ``a Vassiliev invariant of type m'', if
\[ V
\left(\raisebox{2.5mm}{$\displaystyle
  \underbrace{
  \diag{6mm}{4.2}{1}{
  \picmultigraphics{2}{3.2 0}{
  \picvecline{0 0}{1 1}
  \picvecline{1 0}{0 1}
  \picfillgraycol{0}
  \picfilledcircle{0.5 0.5}{0.06}{}
  }
  \picmultigraphics{4}{0.5 0}{
  \picfilledcircle{1.45 0.5}{0.03}{}
  }
  }
  }_{\mbox{\raisebox{0mm}{$\textstyle m+1$}}}
$} \right)\,=\,0
\]

\noindent (that is, if $V$ vanishes when evaluated on a knot with more
then m double points)

\end{defi}

It is easy to show that {\em many} known knot invariants are
Vassiliev, including, for example, all coefficients (in proper
parametrizations) of the Conway, Jones, and HOMFLY polynomials. (See e.g.\
\cite{Bar-Natan:Weights, Birman:Bulletin, BirmanLin:Vassiliev,
Goussarov:New}.)

With polynomials in mind, the following conjecture is just a
variation of Taylor's theorem:

\begin{conj}\label{conjone}
Vassiliev invariants separate knots.
\end{conj}

Little is known about conjecture \ref{conjone}. If ``knots'' are replaced
by ``braids'' \cite{Bar-Natan:Homotopy, Bar-Natan:glN, Kohno:deRham} or
``string links up to homotopy\footnote{allowing change of
self-crossings of the strands}'' \cite{Bar-Natan:Homotopy, Lin:Milnor,
Lin:Expansions}, it is verified. As it stands it sounds very appealing,
but unfortunately,
% it does not live together well with another very
% appealing conjecture, which we will discuss briefly later.
we cannot even yet affirm the following weaker

\begin{ques}\label{quest1}
(see \cite[sect.~7.2]{Bar-Natan:Vassiliev})
Do Vassiliev invariants distinguish knot orientation?
\end{ques}

We will come back to this question in the next lecture.

Whatever you think of conjecture \ref{conjone}, it would clearly be nice to
know what is the set  of all Vassiliev invariants. Let us start:

\begin{defi}
${\cal K}^0_m\,=\,\faktor{\Span\{\mbox{\rm $m$-singular knots}\}}{%
\mbox{\small $\begin{array}{c}
\mbox{differentiability}\\
\mbox{relation}
\end{array}$}}
$,\newline
where the differentiability relation is

\[
\diag{6mm}{2.5}{1}{
\picvecline{0 0}{1 1}
\picvecline{1 0}{0 1}
\picfillgraycol{0}
\picfilledcircle{0.5 0.5}{0.06}{}
\pictranslate{1.5 0}{
\picmultivecline{-4   1 -1 0}{1 0}{0 1}
\picmultivecline{-4   1 -1 0}{0 0}{1 1}
}}\,-\,
\diag{6mm}{2.5}{1}{
\picvecline{0 0}{1 1}
\picvecline{1 0}{0 1}
\picfillgraycol{0}
\picfilledcircle{0.5 0.5}{0.06}{}
\pictranslate{1.5 0}{
\picmultivecline{-4   1 -1 0}{0 0}{1 1}
\picmultivecline{-4   1 -1 0}{1 0}{0 1}
}}\enspace=\enspace
\diag{6mm}{2.5}{1}{
\picmultivecline{-4   1 -1 0}{1 0}{0 1}
\picmultivecline{-4   1 -1 0}{0 0}{1 1}
\pictranslate{1.5 0}{
\picvecline{0 0}{1 1}
\picvecline{1 0}{0 1}
\picfillgraycol{0}
\picfilledcircle{0.5 0.5}{0.06}{}
}}\,-\,
\diag{6mm}{2.5}{1}{
\picmultivecline{-4   1 -1 0}{0 0}{1 1}
\picmultivecline{-4   1 -1 0}{1 0}{0 1}
\pictranslate{1.5 0}{
\picvecline{0 0}{1 1}
\picvecline{1 0}{0 1}
\picfillgraycol{0}
\picfilledcircle{0.5 0.5}{0.06}{}
}}
\]

\end{defi}

\begin{defi}
Let $\delta: {\cal K}^0_{m+1}\,\longrightarrow\,{\cal K}^0_m$ be defined by

\[
\diag{6mm}{1}{1}{
\picvecline{0 0}{1 1}
\picvecline{1 0}{0 1}
\picfillgraycol{0}
\picfilledcircle{0.5 0.5}{0.06}{}
% \piccircle{0.5 0.5}{0.76}{}
}\,\rightarrowto{7mm}\,
\diag{6mm}{1}{1}{
\picmultivecline{-4   1 -1 0}{1 0}{0 1}
\picmultivecline{-4   1 -1 0}{0 0}{1 1}
}\,-\,
\diag{6mm}{1}{1}{
\picmultivecline{-4 1 -1 0}{0 0}{1 1}
\picmultivecline{-4 1 -1 0}{1 0}{0 1}
}\,.
\]
(The differentiability relation ensures that this is well defined.)
\end{defi}

Again use the analogy between Vassiliev invariants and polynomials. If
the $m+1$-st derivative vanishes, the $m$-th derivative should be
a constant. Classifying these ``constants'' ($sym^mE^*$ in the case
of polynomials in Euclidean spaces) is the same as classifying all
polynomials.  That is, if $V$ is Vassiliev of type $m$, it is enough to
understand $V$ on ${\cal K}^0_m$. Thus restricted, it is a linear functional
on ${\cal K}^0_m$, which vanishes on $\delta {\cal K}^0_{m+1}$. That is,

\begin{prop}
To every type $m$ invariant $V$ corresponds an element $W_V$ of
$({\cal K}^0_m/\delta {\cal K}^0_{m+1})^*$.
\end{prop}

\def\entspr{\stackrel{\mbox{\scriptsize def}}{=}}

\begin{propdef} \label{DefCD}
\[
{\cal K}^0_m/\delta {\cal K}^0_{m+1}={\cal D}^0_m\entspr \Span\left\{
\parbox{12mm}{
\vtop{\hbox{\diag{6mm}{2}{2}{
\piccirclevecarc{1 1}{1}{30 390}
%\piclinedash{0.2 0.1}{0.25}
\pictranslate{1 1}{\picline{1 120 polar}{1 240 polar}
\picline{1 100 polar}{1 190 polar}
\picline{1 70 polar}{1 290 polar}
\picline{1 10 polar}{1 200 polar}
\picline{1 260 polar}{1 340 polar}
}}}
\vskip3pt
\hbox to 13.5mm{\hfil\small \,$m$ chords\hfil}
}}
\right\}\,=\,\Span\left\{
  \begin{array}{c}
    \text{degree }m \\
    \text{chord diagrams}
  \end{array}
\right\}.
\]
If a chord diagram $D\in{\cal D}^0_m$ is the image of an $m$-singular knot
$K\in{\cal K}^0_m$ via the projection $F:{\cal K}^0_m\to{\cal K}^0_m/\delta
{\cal K}^0_{m+1}={\cal D}^0_m$ we say that $D$ is the chord diagram
\underline{underlying} $K$ and that the knot $K$ \underline{represents} the
diagram $D$. 
\end{propdef}

We leave the (easy) proof of the assertion in~\ref{DefCD} to the reader.

\begin{ques}\label{quesone}
When does $W\in ({\cal D}^0_m)^*$ integrate to a type $m$ Vassiliev invariant $V$?
\end{ques}

Today's approach is: use induction. Set $V=W$ on ${\cal K}^0_m$, try your luck
integrating it to ${\cal K}^0_{m-1}$, try it again to go on to ${\cal K}^0_{m-2}$,
and keep your fingers crossed hoping to meet no obstruction until you
reach the goal ---~${\cal K}^0_0$.

The obvious question that comes in mind is:

\begin{ques}
When does an invariant in $({\cal K}^0_m)^\star$ integrate one step to
an invariant in $({\cal K}^0_{m-1})^\star$?
\end{ques}

The complete answer to this question is given by the following
theorem, which appears implicitly in Vassiliev~\cite{Vassiliev:CohKnot,
Vassiliev:Book} and explicitly in Stanford~\cite{Stanford:FiniteType}, and
is written as Mike Hutchings~\cite{Hutchings:SingularBraids} writes it.

\begin{theo} \label{theoone}
The sequence
\[
{\cal K}^1_m\,\stackrel\partial\longrightarrow\,{\cal K}^0_m\,\stackrel\delta\longrightarrow\,
{\cal K}^0_{m-1}
\]
is exact. Here ${\cal K}^1_m$ is the space spanned by all singular knots that
have $m-2$ double points and one triple point in which one of the strands is
marked by a $\star$ (called ``Topological 4-Term'' or ``$T4T$'' knots), and
by all singular knots that have $m-1$ double points and one marked point
somewhere on them but not on a double point (called ``Topological Framing
Independence'' or ``$TFI$'' knots), with everything moded out by the
differentiability relation. In pictures,
\[
{\cal K}^1_m\,=\,\faktor{\Span\left\{
\enspace
\diag{7mm}{2}{2.7}{
\pictranslate{0 0.7}{
\picvecline{0.3 0}{1.7 2}
\picvecline{1 0}{1 2}
\picvecline{0 1}{2 1}
\picfillgraycol{0}
\picfilledcircle{1 1}{0.06}{}
}
\picmultigraphics{3}{0.35 0}{
\picfilledcircle{0.6 0.2}{0.03}{}
}
\picmultigraphics{2}{1.5 0}{
\picscale{0.4 0.4}{
\picvecline{0 0}{1 1}
\picvecline{1 0}{0 1}
\picfillgraycol{0}
\picfilledcircle{0.5 0.5}{0.06}{}
}
}\pictext{$\underbrace{\kern1.9\unitlength}_{m-2}$}{0 0}{0 0}
\pictext{$\star$}{1.1 2.5}{0 0}
}\enspace,\enspace
\diag{7mm}{3}{1.1}{
\picline{0 1}{3 1}
\picfillgraycol{0}
\picfilledcircle{1.5 1}{0.06}{}
\picmultigraphics{2}{2.3 0}{
\picscale{0.7 0.7}{
\picvecline{0 0}{1 1}
\picvecline{1 0}{0 1}
\picfilledcircle{0.5 0.5}{0.06}{}
}}
%\tracingmacros1
\picmultigraphics{3}{0.4 0}{
\picfilledcircle{1.13 0.4}{0.04}{}
}
\pictext{$\underbrace{\kern3\unitlength}_{m-1}$}{0 0}{0 0}
}\enspace
\right\}}{\mbox{\small differentiability.}}
\]

\noindent The map $\partial$ is given by
\begin{eqnarray*}
\partial\left(
\diag{7mm}{2}{2}{
\picvecline{0.3 0}{1.7 2}
\picvecline{1 0}{1 2}
\picvecline{0 1}{2 1}
\picfillgraycol{0}
\picfilledcircle{1 1}{0.06}{}
\pictext{$\star$}{1.1 1.8}{0 0}
}
\right)\,
& = &
% \,\sum\,\pm 
% \diag{7mm}{2.8}{2}{
% \picvecline{0.8 0}{0.8 2}
% \picmultivecline{-4 1 -1 0}{0.2 0.4}{2.1 2.0}
% \picvecline{0 1.3}{2.4 1.3}
% \picfillgraycol{0}
% \picfilledcircle{0.8 1.3}{0.06}{}
% \picfilledcircle{1.3 1.3}{0.06}{}
% \picellipsevecarc{1.3 1.3}{0.5 0.3}{149 150}
% \piclinedash{0.08 0.04}{0.025}
% \picellipsearc{1.3 1.3}{0.5 0.3}{30 390}
% \pictext{\mbox{\small loop around}}{1.3 0.8}{0 0}
% }

        \setlength{\unitlength}{0.65\eepiclength}
        \begin{array}{c}  \hspace{-1.7mm}
                \raisebox{-8pt}{\begingroup\makeatletter\ifx\SetFigFont\undefined%
\gdef\SetFigFont#1#2#3#4#5{%
  \reset@font\fontsize{#1}{#2pt}%
  \fontfamily{#3}\fontseries{#4}\fontshape{#5}%
  \selectfont}%
\fi\endgroup%
\begin{picture}(6174,932)(0,-10)
\put(3011,384){\makebox(0,0)[lb]{\smash{{\mathmode{-}}}}}
\put(4661,384){\makebox(0,0)[lb]{\smash{{\mathmode{-}}}}}
\thicklines
\put(612,458){\blacken\ellipse{80}{80}}
\put(612,458){\ellipse{80}{80}}
\put(836,608){\blacken\ellipse{80}{80}}
\put(836,608){\ellipse{80}{80}}
\put(2262,458){\blacken\ellipse{80}{80}}
\put(2262,458){\ellipse{80}{80}}
\put(2486,458){\blacken\ellipse{80}{80}}
\put(2486,458){\ellipse{80}{80}}
\put(3911,459){\blacken\ellipse{80}{80}}
\put(3911,459){\ellipse{80}{80}}
\put(3686,309){\blacken\ellipse{80}{80}}
\put(3686,309){\ellipse{80}{80}}
\put(5336,459){\blacken\ellipse{80}{80}}
\put(5336,459){\ellipse{80}{80}}
\put(5562,459){\blacken\ellipse{80}{80}}
\put(5562,459){\ellipse{80}{80}}
\path(12,458)(612,458)
\path(612,458)(1212,458)
\path(1092.000,428.000)(1212.000,458.000)(1092.000,488.000)
\path(612,458)(1062,758)
\path(978.795,666.474)(1062.000,758.000)(945.513,716.397)
\path(612,458)(162,158)
\path(1662,458)(2262,458)
\path(2262,458)(2862,458)
\path(2742.000,428.000)(2862.000,458.000)(2742.000,488.000)
\path(2262,458)(1812,158)
\path(3911,459)(4511,459)
\path(4391.000,429.000)(4511.000,459.000)(4391.000,489.000)
\path(3911,459)(4361,759)
\path(4277.795,667.474)(4361.000,759.000)(4244.513,717.397)
\path(3911,459)(3461,159)
\path(4962,459)(5562,459)
\path(5562,459)(6162,459)
\path(6042.000,429.000)(6162.000,459.000)(6042.000,489.000)
\path(5562,459)(6012,759)
\path(5928.795,667.474)(6012.000,759.000)(5895.513,717.397)
\path(5562,459)(5112,159)
\path(2265,459)(2420,564)
\path(2490,609)(2710,759)
\path(2627.753,666.613)(2710.000,759.000)(2593.953,716.186)
\path(3310,459)(3610,459)
\path(3685,459)(3910,459)
\path(2261,8)	(2261.000,56.986)
	(2261.000,83.000)

\path(2261,83)	(2256.418,119.495)
	(2261.000,158.000)

\path(2261,158)	(2335.137,199.269)
	(2378.906,213.084)
	(2411.000,233.000)

\path(2411,233)	(2440.400,280.002)
	(2464.216,340.850)
	(2480.175,404.022)
	(2486.000,458.000)

\path(2486,458)	(2480.175,511.977)
	(2464.216,575.150)
	(2440.400,635.997)
	(2411.000,683.000)

\path(2411,683)	(2378.906,702.916)
	(2335.137,716.731)
	(2261.000,758.000)

\path(2261,758)	(2256.418,796.505)
	(2261.000,833.000)

\path(2261,833)	(2261.000,859.014)
	(2261.000,908.000)

\path(2291.000,788.000)(2261.000,908.000)(2231.000,788.000)
\path(3911,9)	(3911.000,57.986)
	(3911.000,84.000)

\path(3911,84)	(3914.889,120.645)
	(3911.000,159.000)

\path(3911,159)	(3865.333,202.077)
	(3798.500,234.000)
	(3731.667,265.923)
	(3686.000,309.000)

\path(3686,309)	(3662.679,377.663)
	(3656.848,417.540)
	(3654.905,459.000)
	(3656.848,500.460)
	(3662.679,540.337)
	(3686.000,609.000)

\path(3686,609)	(3731.667,652.077)
	(3798.500,684.000)
	(3865.333,715.923)
	(3911.000,759.000)

\path(3911,759)	(3914.889,797.355)
	(3911.000,834.000)

\path(3911,834)	(3911.000,860.014)
	(3911.000,909.000)

\path(3941.000,789.000)(3911.000,909.000)(3881.000,789.000)
\path(835,499)	(844.219,568.219)
	(836.000,608.000)

\path(836,608)	(790.405,651.271)
	(723.564,683.173)
	(656.690,714.988)
	(611.000,758.000)

\path(611,758)	(607.111,796.355)
	(611.000,833.000)

\path(611,833)	(611.000,859.014)
	(611.000,908.000)

\path(641.000,788.000)(611.000,908.000)(581.000,788.000)
\path(835,424)	(814.924,350.980)
	(796.904,297.942)
	(760.000,234.000)

\path(760,234)	(727.783,214.717)
	(684.029,200.831)
	(610.000,159.000)

\path(610,159)	(605.418,120.495)
	(610.000,84.000)

\path(610,84)	(610.000,57.986)
	(610.000,9.000)

\path(5365,369)	(5342.276,425.737)
	(5336.000,459.000)

\path(5336,459)	(5341.952,512.893)
	(5357.896,576.075)
	(5381.643,636.969)
	(5411.000,684.000)

\path(5411,684)	(5443.094,703.916)
	(5486.863,717.731)
	(5561.000,759.000)

\path(5561,759)	(5565.583,797.505)
	(5561.000,834.000)

\path(5561,834)	(5561.000,860.014)
	(5561.000,909.000)

\path(5591.000,789.000)(5561.000,909.000)(5531.000,789.000)
\path(5560,9)	(5560.000,57.986)
	(5560.000,84.000)

\path(5560,84)	(5563.052,120.859)
	(5560.000,159.000)

\path(5560,159)	(5511.605,222.728)
	(5465.579,261.288)
	(5400.000,309.000)

\put(1360,384){\makebox(0,0)[lb]{\smash{{\mathmode{+}}}}}
\end{picture} }
                \hspace{-1.9mm}
        \end{array}

\\
\partial\left(
\diag{7mm}{1.5}{1.5}{
\picvecline{0 0}{1.5 1.5}
\picfillgraycol{0}
\picfilledcircle{0.75 0.75}{0.06}{}
}
\right)\,
&=
&\,
\diag{7mm}{1.5}{1.5}{
\piccurve{0 0}{0.7 0.8}{1.3 0.9}{1.4 0.7}
\piccurve{1.5 1.5}{0.5 0.3}{1.6 0.3}{1.4 0.7}
\picfillgraycol{0}
\picfilledcircle{1.1 0.75}{0.06}{}
}.
\end{eqnarray*}

\end{theo}

The fact that $\delta\circ\partial=0$ is easy, and it already implies a
partial answer to question \ref{quesone}:

\def\CD#1{%
  \diag{6mm}{2}{2}{\piccircle{1 1}{1}{}\pictranslate{1 1}{\picscale{-1 1}{#1}}}
}

\begin{prop}\label{propthree}
A necessary condition for $W\in ({\cal D}^0_m)^*$ to integrate to a
Vassiliev invariant is that it vanishes on $\partial {\cal D}^1_m$, where
\[ {\cal D}^1_m\,=\, \Span\left\{ \diag{6mm}{2}{2}{
  \pictext{$\star$}{0.5 1.2}{0 0} \piccircle{1 1}{1}{} \pictranslate{1 1}{
  \picline{1 80 polar}{1 -80 polar}
  \piccurve{1 60 polar}{0.2 0.6}{-0.3 0.6}{1 120 polar}
  \piccurve{1 30 polar}{0.7 0.2}{0.7 -0.2}{1 -30 polar}
  \piccurve{1 -10 polar}{0.6 -0.2}{-0.1 -0.6}{1 260 polar}
  \picvecline{1 190 polar}{1 80 polar p 0.2} } }\,,\,
  \diag{6mm}{2}{2}{ \piccircle{1 1}{1}{} \pictranslate{1 1}{
  \piccurve{1 10 polar}{0.3 0.2}{-0.2 -0.1}{1 230 polar}
  \piccurve{1 130 polar}{-0.3 0.5}{-0.2 0}{1 260 polar}
  \picfillgraycol{0} \picfilledcircle{1 80 polar}{0.09}{} } } \right\}
\]
\[ \partial\left( \diag{6mm}{2}{2}{ \pictext{$\star$}{0.5 1.2}{0 0}
  \piccircle{1 1}{1}{} \pictranslate{1 1}{ \picline{1 80 polar}{1 -80 polar}
  \picvecline{1 190 polar}{1 80 polar p 0} } } \right) \enspace=\enspace
  \CD{\picline{1 110 polar}{1 250 polar} \picline{1 230 polar}{1 0 polar} }
  \enspace-\enspace \CD{\picline{1 130 polar}{1 0 polar}
  \picline{1 110 polar}{1 250 polar} } \enspace+\enspace
  \CD{\picline{1 125 polar}{1 235 polar} \picline{1 0 polar}{1 110 polar} }
  \enspace-\enspace \CD{\picline{1 125 polar}{1 235 polar}
  \picline{1 0 polar}{1 250 polar}
  }\,\kern-2cm\mbox{\raisebox{-9mm}{\small\rm (the $4T$ relation)}},
\]
and
\[ \partial\left( \diag{6mm}{2}{2}{ \piccircle{1 1}{1}{} \pictranslate{1 1}{
  \picfillgraycol{0} \picfilledcircle{1 0 polar}{0.09}{} } } \right)
  \enspace=\enspace \diag{6mm}{2}{2}{ \piccircle{1 1}{1}{} \pictranslate{1 1}{
  \piccurve{1 30 polar}{0.7 0.2}{0.7 -0.2}{1 -30 polar}}
  }\,\mbox{\raisebox{-4mm}{\small\rm (the $FI$ relation)}}
\]
\end{prop}

\proof{}
Just consider the chord diagrams underlying the knots in $\partial
{\cal K}^1_m$.
\qed

\begin{remark} \label{SmallCommutativeDiagram} Notice that ${\cal
D}^1_m={\cal K}^1_m/\delta{\cal K}^1_{m+1}$, where the map $\delta:{\cal
K}^1_{m+1}\to{\cal K}^1_m$ is defined in the same way as the map
$\delta:{\cal K}^0_{m+1}\to{\cal K}^0_m$, and that the map $\partial:{\cal
D}^1_m\to{\cal D}^0_m$ is the only map that makes the following diagram
commutative (with exact rows):
\[ \begin{array}{*7c}
  {\cal K}^1_{m+1}&\stackrel\delta\longrightarrow&
    {\cal K}^1_{m}&\stackrel F\longrightarrow&{\cal D}^1_m&
    \longrightarrow&0\\
  \partial\Big\downarrow&&\partial\Big\downarrow&&\partial\Big\downarrow\\[2mm]
  {\cal K}^0_{m+1}& \stackrel\delta\longrightarrow&{\cal K}^0_m&\stackrel
    F\longrightarrow&{\cal D}^0_m&\longrightarrow&0
\end{array} \]
\end{remark}

\begin{figure}[htpb]
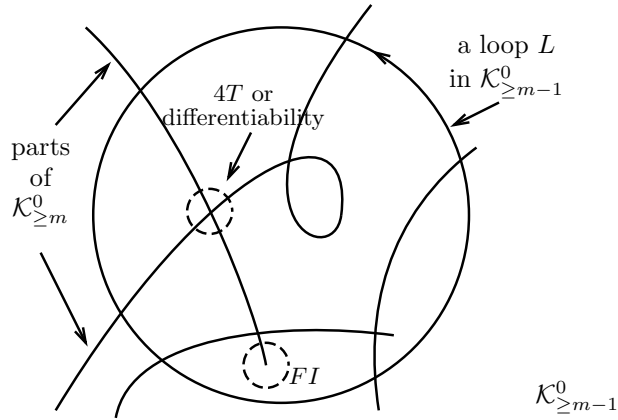

\[
  \diag{1cm}{9}{6}
  {
    \piccirclevecarc{4 3}{2.5}{60 420}
    \pictext{\shortstack{parts\\of\\${\cal K}^0_{\ge m}$}}{0.8 3.5}{-0.5 -0.5}
    \pictext{\small\shortstack{$4T$ or\\differentiability}}{3.5 4.2}{-0.5 0}
    \pictext{\small$FI$}{4.1 0.75}{0 0}
    \pictext{\shortstack{a loop $L$ \\ in ${\cal K}^0_{\ge m-1}$}}{7.0 4.7}{-0.5 0}
    \pictext{${\cal K}^0_{\ge m-1}$}{8.5 0.5}{-1 0}
    \pictranslate{1 0}
    {
      \piccurve{2.8 1.0}{2.7 1.6}{1.8 4.2}{0.4 5.5}
      \piccurve{4.3 0.4}{4.1 1.7}{4.4 3.0}{5.6 3.9}
      \piccurve{0.8 0.3}{1.0 1.4}{3.0 1.6}{4.5 1.4}
      \piccurve{0.0 0.4}{2.0 3.7}{4.0 4.6}{3.8 3.0}
      \piccurve{3.8 3.0}{3.7 2.2}{2.0 3.0}{4.2 5.8}
      \picvecline{0.0 4.2}{0.7 5.0}
      \picvecline{-0.2 2.5}{0.4 1.3}
      \picvecline{2.6 4.1}{2.3 3.5}
      \picvecline{5.9 4.5}{5.3 4.2}
      \piclinedash{0.15 0.07}{0.1}
      \piccircle{2.8 1.0}{0.3}{}
      \piccircle{2.05 3.05}{0.3}{}
    }
  }
\]
\caption{The proof of theorem~\protect\ref{theoone}.} \label{DiskOfKnots}
\end{figure}

\proof{of theorem \ref{theoone} (sketch)} We only need to show that
$\ker\delta\subset\text{im }\partial$.  Take a generic loop $L$ in the
set ${\cal K}^0_{\ge m-1}$ of all parametrized knots with at least $m-1$
double points, and possibly some worse singularities. Such a loop meets
${\cal K}^0_{\ge m}$ in finitely many points, that are $m$-singular
knots. Let $S_L$ be the (properly signed) sum of these $m$-singular
knots. It is not hard to show that $\ker\delta$ is spanned by these
$S_L$'s, so it is enough to show that $S_L$ is in $\text{im }\partial$
for any $L$. Now notice that ${\cal K}^0_{\ge m-1}$ is simply connected,
so $L$ bounds some generic disk $D$ in ${\cal K}^0_{\ge m-1}$. The
intersection of $D$ with the codimension $1$ set of knots of a higher
singularity is some graph $G$ on $D$ (see figure~\ref{DiskOfKnots}),
and the vertices of $G$ correspond to points in the codimension $2$ set
of generic knots of an even higher singularities. One can check that
this set is exactly the set of generators of ${\cal K}^1_m$, and that
$S_L=\delta S_D$ where $S_D$ is the (properly signed) sum in ${\cal
K}^1_m$ corresponding to the vertices of $G$.  \qed

\noindent\underline{\bf The Fundamental Theorem of Vassiliev invariants}
{\it The condition in proposition~\ref{propthree} is also sufficient.}

Let
\[
{\cal A}^r_m={\cal D}^0_m/\delta {\cal D}^1_m= \left(\vcbox{\hbox{chord diagrams}\hbox{$\bmod$
 $4T$ \& $FI$}}\right).
\]

Then every weight system $W$ (an element in $({\cal A}^r_m)^*$) integrates to a
Vassiliev invariant. It follows that the associated graded vector space
of the filtered space of all  Vassiliev invariants is
\[
  ({\cal A}^r)^*\entspr \left(\bigoplus_{m=0}^{\infty} {\cal A}^r_m \right)^*
  \qquad
  \text{(duals are taken in the graded sense)}.
\]

There are two problems with this lovely theorem
\begin{enumerate}
\item Although much is known about ${\cal A}^r$ (and its equivalent but
  friendlier version ${\cal A}$ in which the $FI$ relation is not imposed),
  we are far from understanding it.
\item As indicated in the introduction, we know at least four approaches
  to the proof. The topological approach of this lecture, which fails,
  but comes close. And three other approaches, geometrical, physical,
  and algebraic, that all work, but have other defects.
\end{enumerate}

\subsection{Hutchings' combinatorial-topological approach}
\label{Hutchings}

In view of theorem~\ref{theoone}, the Fundamental Theorem follows from
the following:

\begin{conj}\label{conjthree} Any invariant satisfying the $T4T$ and $TFI$
can be integrated one step to an invariant \underline{that does the same}.
\end{conj}

In~\cite{Hutchings:SingularBraids}, M.~Hutchings was able to reduce
this conjecture to a statement that appears to be easier to verify
(``Hutchings' condition'', below), and to show that this statement follows
from a completely combinatorial statement (conjecture~\ref{H1Vanishes}).

\subsubsection{Hutchings' condition.} Consider the following commutative
diagram:

\[
\begin{array}{*9c}
&&0&&0&&0\\
&&\Big\downarrow&&\Big\downarrow&&\Big\downarrow&\\[2mm]
&&{\cal K}^1_{m+1}&\stackrel\delta\longrightarrow&
\ker \partial|_{{\cal K}^1_{m}}&\stackrel F\longrightarrow&
\ker \partial|_{{\cal D}^1_m}\raisebox{0.6ex}{\snakest}\\
&&\Big\downarrow&&\Big\downarrow&&\Big\downarrow&\\[2mm]
&&{\cal K}^1_{m+1}&\stackrel\delta\longrightarrow&
{\cal K}^1_{m}&\stackrel F\longrightarrow&{\cal D}^1_m&\longrightarrow&0\\
&&0\Big\downarrow&&\partial\Big\downarrow&&\partial\Big\downarrow\\[2mm]
0&\longrightarrow&{\cal K}^0_{m+1}/\partial {\cal K}^1_{m+1}&
\stackrel\delta\longrightarrow&{\cal K}^0_m&\stackrel
F\longrightarrow&{\cal D}^0_m&(\,\longrightarrow&0\,)\\
&&\Big\downarrow&&\Big\downarrow&&\Big\downarrow&\\[2mm]
&&
%\mbox{\hskip-1ex\hbox to 0pt{
\raisebox{0.3ex}{\snakeend}
\,\,\,
%\hss}\hskip1ex}
{\cal K}^0_{m+1}/\partial
{\cal K}^1_{m+1}&\stackrel\delta\longrightarrow&
{\cal K}^0_{m}/\partial {\cal K}^1_{m}&\longrightarrow&{\cal A}_m^r\\
&&\Big\downarrow&&\Big\downarrow&&\Big\downarrow&\\[2mm]
&&0&&0&&0
\end{array}
\]

The columns of this diagrams are exact by definition. The second row is
exact as in remark~\ref{SmallCommutativeDiagram}. The third row is exact
(though we will not use its exactness at the right end) because it is a
folding (on the left) of the sequence
\[ {\cal K}^1_{m+1}\stackrel\partial\longrightarrow
  {\cal K}^0_{m+1}\stackrel\delta\longrightarrow
  {\cal K}^0_m\stackrel F\longrightarrow
  {\cal D}^0_m\longrightarrow 0,
\]
whose left half is exact by theorem~\ref{theoone} and whose right half is
exact as in remark~\ref{SmallCommutativeDiagram}.

Rephrased in an algebraic language, conjecture \ref{conjthree}
says that $\delta^*:({\cal K}^0_m/\partial {\cal K}^1_m)^*\to ({\cal
K}^0_{m+1}/\partial {\cal K}^1_{m+1})^*$ is surjective.  This is equivalent
to the injectivity of $\delta: {\cal K}^0_{m+1}/\partial {\cal
K}^1_{m+1}\to {\cal K}^0_{m}/\partial {\cal K}^1_{m}$. By the snake
lemma applied to the above diagram, this is equivalent to the surjectivity
of
\[
F: \ker \partial|_{{\cal K}^1_{m}}\longrightarrow\ker \partial|_{{\cal D}^1_m}
\]
In other words, it is enough to prove ``Hutchings' condition'', saying
that
\begin{itemize}
\item[*]
Every relation between $4T$ \& $FI$ relations (on the level of diagrams)
lifts to a relation between $T4T$ \& $TFI$ relations (on the level of knots,
and mod the differentiability relation).
\end{itemize}

\subsubsection{A possible strategy.}

\begin{enumerate}
\item Find many elements of $\ker \partial|_{{\cal D}^1_m}$. Namely, find a big
${\cal D}^2_m$ and a map $\partial:{\cal D}^2_m\to{\cal D}^1_m$ so that
\[
  {\cal D}^2_m\,\stackrel\partial\longrightarrow\,
  {\cal D}^1_m\,\stackrel\partial\longrightarrow\,
  {\cal D}^0_m\,\longrightarrow\, {\cal A}^r_m
\]
is exact.
\item Prove that $F$ is onto $\partial{\cal D}^2_m$.
\item Compute $H_\partial^1({\cal D}^*_m)$. If it is $0$, you win.
\end{enumerate}

\begin{remark} In~\cite{Hutchings:SingularBraids}, M.~Hutchings proved
conjecture \ref{conjthree} for braids (and hence the Fundamental Theorem
for braids) by following this strategy.
\end{remark}

Anyway, ignoring $FI$ for simplicity, here's a candidate for ${\cal D}^2_m$
(which worked well for braids):
\[
{\cal D}^2_m=\Span\left\{
\vtop{\hbox{\diag{6mm}{2}{2}{
\piccircle{1 1}{1}{}
\pictranslate{1 1}{
\picline{0 0}{1 90 polar}
\picline{0 0}{1 210 polar}
\picline{0 0}{1 330 polar}
}}}
\vbox to 0pt{
\vskip2pt\hbox to 12mm{\hss $3T$\hss}\vss}
}
\,,\,
\vtop{\hbox{\diag{6mm}{2}{2}{
\piccircle{1 1}{1}{}
\pictranslate{1 1}{
\piccurve{1 -30 polar}{0.5 -0.1}{0.4 0.3}{1 70 polar}
\picveccurve{1 -80 polar}{0.2 -0.7}{0.4 0.1}{0.45 0.25}
\picscale{-1 1}{
\piccurve{1 -30 polar}{0.5 -0.1}{0.4 0.3}{1 70 polar}
\picveccurve{1 -80 polar}{0.2 -0.7}{0.4 0.1}{0.45 0.25}
}
}}}
\vbox to 0pt{
\vskip2pt\hbox to 12mm{\hss $8T$\hss}\vss}
}
\,,\,
\vtop{\hbox{\diag{6mm}{2}{2}{
\piccircle{1 1}{1}{}
\pictranslate{1 1}{
\picline{1 60 polar}{1 -60 polar}
\picveccurve{1 120 polar}{-0.3 0.5}{0.1 0.2}{0.5 0}
\picveccurve{1 240 polar}{-0.5 -0.6}{-0.3 -0.1}{0.0 0.3}
}
}}
\vbox to 0pt{
\vskip2pt\hbox to 12mm{\hss $14T$\hss}\vss}
}
\right\},
\]
where
\begin{eqnarray*}
\partial\left(
\diag{6mm}{2}{2}{
\piccircle{1 1}{1}{}
\pictranslate{1 1}{
\picline{0 0}{1 90 polar}
\picline{0 0}{1 210 polar}
\picline{0 0}{1 330 polar}
}}
\right)&\enspace=\enspace&
\diag{6mm}{2}{2}{
\piccircle{1 1}{1}{}
\pictranslate{1 1}{
\picline{1 90 polar}{0 0}
\picline{1 210 polar}{0 0}
\picvecline{1 330 polar}{0 0}
}}\,+\,
\diag{6mm}{2}{2}{
\piccircle{1 1}{1}{}
\pictranslate{1 1}{
\picvecline{1 90 polar}{0 0}
\picline{1 210 polar}{0 0}
\picline{1 330 polar}{0 0}
}}\,+\,
\diag{6mm}{2}{2}{
\piccircle{1 1}{1}{}
\pictranslate{1 1}{
\picline{1 90 polar}{0 0}
\picvecline{1 210 polar}{0 0}
\picline{1 330 polar}{0 0}
}}
\\
\partial\left(
\diag{6mm}{2}{2}{
\piccircle{1 1}{1}{}
\pictranslate{1 1}{
\piccurve{1 -30 polar}{0.5 -0.1}{0.4 0.3}{1 70 polar}
\picveccurve{1 -80 polar}{0.2 -0.7}{0.4 0.1}{0.45 0.25}
\picscale{-1 1}{
\piccurve{1 -30 polar}{0.5 -0.1}{0.4 0.3}{1 70 polar}
\picveccurve{1 -80 polar}{0.2 -0.7}{0.4 0.1}{0.45 0.25}
}
}
}
\right)&\enspace=\enspace&
\diag{6mm}{2}{2}{
\piccircle{1 1}{1}{}
\pictranslate{1 1}{
\picline{1 -60 polar}{1 60 polar}
\picline{1 130 polar}{1 230 polar}
\picveccurve{1 -80 polar}{0.1 -0.7}{0.3 -0.3}{0.5 0.2}
\piccurve{1 -95 polar}{-0.1 -0.7}{-0.3 -0.6}{1 240 polar}
}
}
\,-\,
\diag{6mm}{2}{2}{
\piccircle{1 1}{1}{}
\pictranslate{1 1}{
\picline{1 -60 polar}{1 60 polar}
\picline{1 130 polar}{1 230 polar}
\picveccurve{1 -80 polar}{0.1 -0.7}{0.3 -0.3}{0.5 0.2}
\piccurve{1 -95 polar}{-0.3 -0.5}{-0.5 -0.4}{1 215 polar}
}
}
\,+\,
\diag{6mm}{2}{2}{
\piccircle{1 1}{1}{}
\pictranslate{1 1}{
\picline{1 -60 polar}{1 60 polar}
\picline{1 130 polar}{1 230 polar}
\picveccurve{1 -80 polar}{0.1 -0.7}{0.3 -0.3}{0.5 0.2}
\piccurve{1 -95 polar}{-0.3 -0.1}{-0.6 0.3}{1 145 polar}
}
}
\,-\,
\diag{6mm}{2}{2}{
\piccircle{1 1}{1}{}
\pictranslate{1 1}{
\picline{1 -60 polar}{1 60 polar}
\picline{1 130 polar}{1 230 polar}
\picveccurve{1 -80 polar}{0.1 -0.7}{0.3 -0.3}{0.5 0.2}
\piccurve{1 -95 polar}{-0.2 -0.1}{-0.3 0.4}{1 125 polar}
}
}\\
&&
\llap{$-\enspace$}
\diag{6mm}{2}{2}{
\piccircle{1 1}{1}{}
\pictranslate{1 1}{
\picscale{-1 1}{
\picline{1 -60 polar}{1 60 polar}
\picline{1 130 polar}{1 230 polar}
\picveccurve{1 -80 polar}{0.1 -0.7}{0.3 -0.3}{0.5 0.2}
\piccurve{1 -95 polar}{-0.1 -0.7}{-0.3 -0.6}{1 240 polar}
}}
}
\,+\,
\diag{6mm}{2}{2}{
\piccircle{1 1}{1}{}
\pictranslate{1 1}{
\picscale{-1 1}{
\picline{1 -60 polar}{1 60 polar}
\picline{1 130 polar}{1 230 polar}
\picveccurve{1 -80 polar}{0.1 -0.7}{0.3 -0.3}{0.5 0.2}
\piccurve{1 -95 polar}{-0.3 -0.5}{-0.5 -0.4}{1 215 polar}
}}
}
\,-\,
\diag{6mm}{2}{2}{
\piccircle{1 1}{1}{}
\pictranslate{1 1}{
\picscale{-1 1}{
\picline{1 -60 polar}{1 60 polar}
\picline{1 130 polar}{1 230 polar}
\picveccurve{1 -80 polar}{0.1 -0.7}{0.3 -0.3}{0.5 0.2}
\piccurve{1 -95 polar}{-0.3 -0.1}{-0.6 0.3}{1 145 polar}
}}
}
\,+\,
\diag{6mm}{2}{2}{
\piccircle{1 1}{1}{}
\pictranslate{1 1}{
\picscale{-1 1}{
\picline{1 -60 polar}{1 60 polar}
\picline{1 130 polar}{1 230 polar}
\picveccurve{1 -80 polar}{0.1 -0.7}{0.3 -0.3}{0.5 0.2}
\piccurve{1 -95 polar}{-0.2 -0.1}{-0.3 0.4}{1 125 polar}
}}
}
\\
\partial\left(
\diag{6mm}{2}{2}{
\piccircle{1 1}{1}{}
\pictranslate{1 1}{
\picline{1 60 polar}{1 -60 polar}
\picveccurve{1 120 polar}{-0.3 0.5}{0.1 0.2}{0.5 0}
\picveccurve{1 240 polar}{-0.5 -0.6}{-0.3 -0.1}{0.0 0.3}
}
}
\right)&\enspace=\enspace& 
        \setlength{\unitlength}{0.5\eepiclength}
        \begin{array}{c}  \hspace{-1.7mm}
                \raisebox{-8pt}{\begingroup\makeatletter\ifx\SetFigFont\undefined%
\gdef\SetFigFont#1#2#3#4#5{%
  \reset@font\fontsize{#1}{#2pt}%
  \fontfamily{#3}\fontseries{#4}\fontshape{#5}%
  \selectfont}%
\fi\endgroup%
\begin{picture}(5183,1229)(0,-10)
\put(1800,531){\makebox(0,0)[lb]{\smash{{\mathmode{-}}}}}
\put(3599,532){\makebox(0,0)[lb]{\smash{{\mathmode{+}}}}}
\thicklines
\put(974,607){\ellipse{1200}{1200}}
\put(2775,607){\ellipse{1200}{1200}}
\put(4575,607){\ellipse{1200}{1200}}
\path(1424,982)(1424,232)
\path(524,982)(1424,607)
\path(1301.692,625.462)(1424.000,607.000)(1324.769,680.846)
\path(3225,982)(3225,232)
\path(2325,982)(3225,607)
\path(3102.692,625.462)(3225.000,607.000)(3125.769,680.846)
\path(5025,982)(5025,232)
\path(4125,982)(5025,607)
\path(4902.692,625.462)(5025.000,607.000)(4925.769,680.846)
\path(524,232)(449,907)
\path(2324,232)(2399,1057)
\path(4124,232)(4949,1057)
\put(0,531){\makebox(0,0)[lb]{\smash{{\mathmode{\ }}}}}
\end{picture} }
                \hspace{-1.9mm}
        \end{array}
 \\
& & 
        \setlength{\unitlength}{0.5\eepiclength}
        \begin{array}{c}  \hspace{-1.7mm}
                \raisebox{-8pt}{\begingroup\makeatletter\ifx\SetFigFont\undefined%
\gdef\SetFigFont#1#2#3#4#5{%
  \reset@font\fontsize{#1}{#2pt}%
  \fontfamily{#3}\fontseries{#4}\fontshape{#5}%
  \selectfont}%
\fi\endgroup%
\begin{picture}(5183,1229)(0,-10)
\put(0,532){\makebox(0,0)[lb]{\smash{{\mathmode{-}}}}}
\put(1800,532){\makebox(0,0)[lb]{\smash{{\mathmode{+}}}}}
\put(3600,532){\makebox(0,0)[lb]{\smash{{\mathmode{-}}}}}
\thicklines
\put(4575,607){\ellipse{1200}{1200}}
\put(2775,607){\ellipse{1200}{1200}}
\put(974,607){\ellipse{1200}{1200}}
\path(5025,982)(5025,232)
\path(4125,982)(5025,607)
\path(4902.692,625.462)(5025.000,607.000)(4925.769,680.846)
\path(3225,982)(3225,232)
\path(2325,982)(3225,607)
\path(3102.692,625.462)(3225.000,607.000)(3125.769,680.846)
\path(1424,982)(1424,232)
\path(524,982)(1424,607)
\path(1301.692,625.462)(1424.000,607.000)(1324.769,680.846)
\path(524,233)(1499,908)
\path(2324,233)(3299,308)
\path(4124,233)(4949,158)
\end{picture} }
                \hspace{-1.9mm}
        \end{array}
 \\
& & 
        \setlength{\unitlength}{0.5\eepiclength}
        \begin{array}{c}  \hspace{-1.7mm}
                \raisebox{-8pt}{\begingroup\makeatletter\ifx\SetFigFont\undefined%
\gdef\SetFigFont#1#2#3#4#5{%
  \reset@font\fontsize{#1}{#2pt}%
  \fontfamily{#3}\fontseries{#4}\fontshape{#5}%
  \selectfont}%
\fi\endgroup%
\begin{picture}(6983,1229)(0,-10)
\put(0,533){\makebox(0,0)[lb]{\smash{{\mathmode{-}}}}}
\put(1800,533){\makebox(0,0)[lb]{\smash{{\mathmode{+}}}}}
\put(3600,533){\makebox(0,0)[lb]{\smash{{\mathmode{-}}}}}
\put(5400,533){\makebox(0,0)[lb]{\smash{{\mathmode{+}}}}}
\thicklines
\put(6375,607){\ellipse{1200}{1200}}
\put(4575,607){\ellipse{1200}{1200}}
\put(2775,607){\ellipse{1200}{1200}}
\put(974,607){\ellipse{1200}{1200}}
\path(6825,982)(6825,232)
\path(5025,982)(5025,232)
\path(3225,982)(3225,232)
\path(1424,982)(1424,232)
\path(524,976)(1349,1051)
\path(2324,1000)(3299,925)
\path(4124,987)(5099,312)
\path(5924,983)(6749,158)
\path(524,233)(936,1006)
\path(906.032,885.992)(936.000,1006.000)(853.083,914.213)
\path(2324,233)(2781,961)
\path(2742.608,843.416)(2781.000,961.000)(2691.791,875.316)
\path(4124,233)(4649,608)
\path(4568.789,513.839)(4649.000,608.000)(4533.915,562.663)
\path(5924,233)(6374,533)
\path(6290.795,441.474)(6374.000,533.000)(6257.513,491.397)
\end{picture} }
                \hspace{-1.9mm}
        \end{array}
 \\
& & 
        \setlength{\unitlength}{0.5\eepiclength}
        \begin{array}{c}  \hspace{-1.7mm}
                \raisebox{-8pt}{\begingroup\makeatletter\ifx\SetFigFont\undefined%
\gdef\SetFigFont#1#2#3#4#5{%
  \reset@font\fontsize{#1}{#2pt}%
  \fontfamily{#3}\fontseries{#4}\fontshape{#5}%
  \selectfont}%
\fi\endgroup%
\begin{picture}(6984,1229)(0,-10)
\put(0,457){\makebox(0,0)[lb]{\smash{{\mathmode{-}}}}}
\put(1800,532){\makebox(0,0)[lb]{\smash{{\mathmode{+}}}}}
\put(3600,532){\makebox(0,0)[lb]{\smash{{\mathmode{-}}}}}
\put(5400,532){\makebox(0,0)[lb]{\smash{{\mathmode{+}}}}}
\thicklines
\put(975,607){\ellipse{1200}{1200}}
\put(2776,607){\ellipse{1200}{1200}}
\put(4576,607){\ellipse{1200}{1200}}
\put(6376,607){\ellipse{1200}{1200}}
\path(1425,982)(1425,232)
\path(3226,982)(3226,232)
\path(5026,982)(5026,232)
\path(6826,982)(6826,232)
\path(525,982)(1350,1057)
\path(2325,982)(3300,907)
\path(4125,982)(5100,307)
\path(5925,982)(6750,157)
\path(525,232)(1425,607)
\path(1325.769,533.154)(1425.000,607.000)(1302.692,588.538)
\path(2321,232)(3221,607)
\path(3121.769,533.154)(3221.000,607.000)(3098.692,588.538)
\path(4117,232)(5017,607)
\path(4917.769,533.154)(5017.000,607.000)(4894.692,588.538)
\path(5914,232)(6814,607)
\path(6714.769,533.154)(6814.000,607.000)(6691.692,588.538)
\end{picture} }
                \hspace{-1.9mm}
        \end{array}
 \\
\end{eqnarray*}

It is not hard to lift $3T$, $8T$ and $14T$ to a $T3T$, $T8T$ and $T14T$ in
$\ker \partial|_{{\cal K}^1_m}$.\vspace{3mm}

\begin{conj} \label{H1Vanishes}
$H_\partial^1({\cal D}^*_m)=0$ \hspace{3cm}$\left(\parbox{2.2in}{To
be honest, we hope it's false. This will make life more
interesting!}\right)$\vspace{3mm}
\end{conj}

Notice that this is a diagram level statement, which implies the
Fundamental Theorem!

\begin{conj}\label{conjfour}
$H_\partial^1({\cal D}^*_m)$ is isomorphic to (a certain twist of)
Kontsevich's graph homology.
\end{conj}

Proving conjecture \ref{conjfour} appears to be only a matter of labor.
\vskip5mm

\begin{remark} See Domergue-Donato~\cite{DomergueDonato:Integrating}
and Willerton~\cite{Willerton:HalfIntegration} for some other partial
results on the combinatorial-topological approach.
Some enumerative results on chord diagrams appear
in~\cite{Stoimenow:Number}.
\end{remark}

\subsection{Why are we not happy?}

\begin{enumerate}
\item The construction of the diagram on which the snake lemma was applied was
  somewhat artificial. Is there something more basic going on?
\item We don't know that $H_\partial^1({\cal D}^*_m)=0$.  We believe, our
  ${\cal D}^2_m$ is the right one, but it may well be that $H^1_\partial$
  \underline{does not} vanish, and that its non-triviality means
  something. What does it mean?
\end{enumerate}

\lectitle{Geometry}
\def\missingref#1{{\typeout{*** Missing ref #1 ***}\mbox{\tiny\tt #1}}}

\def\cA{{\cal A}}
\def\tZ{{\tilde Z}}

\subsection{A short review of lecture 1.}

Generalize a knot invariant $V$ (a map $\{\text{knots up to
isotopy}\}\to\bold C$) to singular knots by
\begin{equation}\label{diff}
V\left(\diag{6mm}{1}{1}{
\picvecline{0 0}{1 1}
\picvecline{1 0}{0 1}
\picfillgraycol{0}
\picfilledcircle{0.5 0.5}{0.06}{}
% \piccircle{0.5 0.5}{0.76}{}
}\right)\,:=\,
V\left(\diag{6mm}{1}{1}{
\picmultivecline{-4 1 -1 0}{1 0}{0 1}
\picmultivecline{-4 1 -1 0}{0 0}{1 1}
}\right)\,-\,
V\left(\diag{6mm}{1}{1}{
\picmultivecline{-4 1 -1 0}{0 0}{1 1}
\picmultivecline{-4 1 -1 0}{1 0}{0 1}
}\right)\,,
\end{equation}
and then define a Vassiliev invariant
\[
\left(\vcbox{\hbox{\shortstack{$V$ is of\\type $m$}}}\right)\,\iff\,
V
\kern1.6ex\underbrace{
\kern-1.6ex
\left(
\,\diag{6mm}{4.2}{1}{
\picmultigraphics{2}{3.2 0}{
\picvecline{0 0}{1 1}
\picvecline{1 0}{0 1}
\picfillgraycol{0}
\picfilledcircle{0.5 0.5}{0.06}{}
}
\picmultigraphics{4}{0.5 0}{
\picfilledcircle{1.45 0.5}{0.03}{}
}
}
\,\right)\kern-1.6ex}%
_{\mbox{\raisebox{-0mm}{$\textstyle m+1$}}}
\kern1.6ex\,=\,0\,.
\]

On can think of (\ref{diff}) as of ``differentiating'' an invariant
and of a Vassiliev invariant as of a ``polynomial''. So, to understand
them we would like to know their ``coefficients''. Here is a nice
candidate.
\[
(V \text{ of type $m$})\,\Longrightarrow
V
%\underbrace{
\left(
\diag{6mm}{4.5}{1}{
\picmultigraphics{2}{2.2 0}{
\picvecline{0 0}{1 1}
\picvecline{1 0}{0 1}
\picfillgraycol{0}
\picfilledcircle{0.5 0.5}{0.06}{}
}
\picmultigraphics{3}{0.4 0}{
\picfilledcircle{1.25 0.5}{0.03}{}
}
\pictranslate{3.5 0.0}{
\picmultivecline{-4 1 -1 0}{1 0}{0 1}
\picmultivecline{-4 1 -1 0}{0 0}{1 1}
}
}
\right)%}_{\mbox{\raisebox{-3mm}{$\textstyle m+1$}}}\,=\,0\,.
\,=\,
V
%\underbrace{
\left(
\diag{6mm}{4.5}{1}{
\picmultigraphics{2}{2.2 0}{
\picvecline{0 0}{1 1}
\picvecline{1 0}{0 1}
\picfillgraycol{0}
\picfilledcircle{0.5 0.5}{0.06}{}
}
\picmultigraphics{3}{0.4 0}{
\picfilledcircle{1.25 0.5}{0.03}{}
}
\pictranslate{3.5 0.0}{
\picmultivecline{-4 1 -1 0}{0 0}{1 1}
\picmultivecline{-4 1 -1 0}{1 0}{0 1}
}
}
\right)%}_{\mbox{\raisebox{-3mm}{$\textstyle m+1$}}}\,=\,0\,.
\,,
\]

and that's why $V$ defines
\[
W_V\,:\,
\Span\left\{
\parbox{8mm}{
\vtop{\hbox{\diag{4mm}{2}{2}{
\piccirclevecarc{1 1}{1}{30 390}
%\piclinedash{0.2 0.1}{0.25}
\pictranslate{1 1}{\picline{1 120 polar}{1 240 polar}
\picline{1 100 polar}{1 190 polar}
\picline{1 70 polar}{1 290 polar}
\picline{1 10 polar}{1 200 polar}
\picline{1 260 polar}{1 340 polar}
}}}
\vskip3pt
\hbox to 9mm{\hfil\tiny $\!m$ chords\hfil}
}}
\right\}\,\longrightarrow\,\bold C\,.
\]
This $W_V$ satisfies two relations ($4T$ and $FI$) because
of topological reasons and hence it becomes a {\em
weight system}\/ $W_V\in(\cA^r_m)^*$, where

\def\CD#1{%
\diag{3mm}{2}{2}{
\piccircle{1 1}{1}{}
%\piccirclevecarc{1 1}{1}{100 145}
%\piccirclevecarc{1 1}{1}{220 265}
%\piclinedash{0.1}{0.05}
%\piccirclearc{1 1}{1}{25 100}
%\piccirclearc{1 1}{1}{145 220}
%\piccirclearc{1 1}{1}{265 340}
%\piclinedash{0.2 0.1}{0.25}
\pictranslate{1 1}{\picscale{-1 1}{#1}}}}

\[
\cA^r\enspace=\enspace
\Span\left\{%
\diag{6mm}{2}{2}{
\piccirclevecarc{1 1}{1}{30 390}
%\piclinedash{0.2 0.1}{0.25}
\pictranslate{1 1}{\picline{1 120 polar}{1 240 polar}
\picline{1 100 polar}{1 190 polar}
\picline{1 70 polar}{1 290 polar}
\picline{1 10 polar}{1 200 polar}
\picline{1 260 polar}{1 340 polar}
}}
\right\}\,\Bigg/\,
\begin{array}{l}
4T\,:\enspace
\CD{\picline{1 110 polar}{1 250 polar}
\picline{1 230 polar}{1 0 polar}
}
\,-\,
\CD{\picline{1 130 polar}{1 0 polar}
\picline{1 110 polar}{1 250 polar}
}
\,+\,
\CD{\picline{1 125 polar}{1 235 polar}
\picline{1 0 polar}{1 110 polar}
}
\,-\,
\CD{\picline{1 125 polar}{1 235 polar}
\picline{1 0 polar}{1 250 polar}
}\,=0\\[2mm]
FI\,:\enspace\diag{3mm}{2}{2}{
\piccircle{1 1}{1}{}
\pictranslate{1 1}{
\piccurve{1 30 polar}{0.7 0.2}{0.7 -0.2}{1 -30 polar}}
}\,=0
\end{array}
\]

Now the following theorem tells us that this is exactly what
we were looking for.

\noindent\underline{\bf The Fundamental Theorem}
{\it Every $W\in(\cA^r_m)^*$ is $W_V$ for some type $m$ invariant $V$.}

It turns out to be interesting to explore these combinatorial
objects. So, before we start proving the Fundamental Theorem,
let's say something more about them.

%here comes sth.

\subsection{A word about Lie algebras}  (drop $FI$ for convenience,
i.~e., consider {\em framed}\/ knots)

There is a way to construct a weight system out of a Lie algebra
representation. First we need the following

\begin{theo}\label{TrivalentVertices} (\cite{Bar-Natan:Vassiliev})
There is an equivalent representation of our diagram space $\cA$ in terms
of diagrams in which some number of oriented internal trivalent vertices
are also allowed. Namely

\def\CD#1{%
\diag{2.5mm}{2}{2}{
\piccircle{1 1}{1}{}
\pictranslate{1 1}{\picscale{-1 1}{#1}}}}

\begin{eqnarray*}
\cA\enspace&=&\enspace
\Span\left\{%
\diag{6mm}{2}{2}{
\piccirclevecarc{1 1}{1}{30 390}
%\piclinedash{0.2 0.1}{0.25}
\pictranslate{1 1}{\picline{1 120 polar}{1 240 polar}
\picline{1 100 polar}{1 190 polar}
\picline{1 70 polar}{1 290 polar}
\picline{1 10 polar}{1 200 polar}
\picline{1 260 polar}{1 340 polar}
}}
\right\}\,\Bigg/\,
\begin{array}{l}
\quad 4T\,:\enspace\\[1mm]
\CD{\picline{1 110 polar}{1 250 polar}
\picline{1 230 polar}{1 0 polar}
}
\,-\,
\CD{\picline{1 130 polar}{1 0 polar}
\picline{1 110 polar}{1 250 polar}
}
\,+\,
\CD{\picline{1 125 polar}{1 235 polar}
\picline{1 0 polar}{1 110 polar}
}
\,-\,
\CD{\picline{1 125 polar}{1 235 polar}
\picline{1 0 polar}{1 250 polar}
}\,=0\\[2mm]
\end{array}\\
&\cong&\enspace
\Span\left\{
\diag{6mm}{2}{2}{
\piccirclevecarc{1 1}{1}{30 390}
\pictranslate{1 1}{
\picveclength{0.1}
\picvecwidth{0.05}
\piccirclevecarc{-0.2 0}{0.2}{0 250}
%\piclinedash{0.2 0.1}{0.25}
\picline{1 120 polar}{-0.2 0}
\picline{1 240 polar}{-0.2 0}
\picline{1 0 polar}{-0.2 0}
\picline{1 70 polar}{1 -70 polar}
}
}
\right\}\,\Bigg/\,
\begin{array}{ll}
AS\,:&
\diag{2mm}{2}{2}
{%\piclinedash{0.2 0.1}{0.25}
\picline{1 1}{2 2}
\picline{1 1}{0 2}
\picline{1 1}{1 0}
%\piclinedash{}{0}
% \picveclength{0.1}
% \picvecwidth{0.05}
\piccirclevecarc{1 1}{0.4}{45 270}
}
\,=\,
-\,
\diag{2mm}{2}{2}
{%\piclinedash{0.2 0.1}{0.25}
\picline{1 1}{2 2}
\picline{1 1}{0 2}
\picline{1 1}{1 0}
%\piclinedash{}{0}
% \picveclength{0.1}
% \picvecwidth{0.05}
\piccirclevecarc{1 1}{0.4}{-90 -225}
}\,,\\
IHX\,:&
\diag{2mm}{2}{2}
{%\piclinedash{0.2 0.1}{0.25}
\picline{0 2}{2 2}
\picline{1 0}{1 2}
\picline{0 0}{2 0}
}\,=\,
\diag{2mm}{2}{2}
{%\piclinedash{0.2 0.1}{0.25}
\picline{0 2}{0 0}
\picline{0 1}{2 1}
\picline{2 2}{2 0}
}\,-\,
\diag{2mm}{2}{2}
{%\piclinedash{0.2 0.1}{0.25}
\picline{0 0}{2 2}
\picline{0.6 0.6}{1.4 0.6}
\picline{0 2}{2 0}
}\,
\\
STU\,:&
\diag{2mm}{2}{2}
{\piccirclevecarc{1 2}{2}{240 300}
%\piclinedash{0.2 0.1}{0.25}
\picline{1.9 2}{1 1.2}
\picline{0.1 2}{1 1.2}
\picline{1 1.2}{1 0}
%\piclinedash{}{0}
%\picveclength{0.15}
%\picvecwidth{0.08}
%\piccirclevecarc{1 1.2}{0.3}{135 45}
}
\,=\,
\diag{2mm}{2}{2}
{\piccirclevecarc{1 2}{2}{240 300}
%\piclinedash{0.2 0.1}{0.25}
\picline{1.9 2}{1.2 0}
\picline{0.1 2}{0.8 0}
}
\,-\,
\diag{2mm}{2}{2}
{\piccirclevecarc{1 2}{2}{240 300}
%\piclinedash{0.2 0.1}{0.25}
\picline{1.9 2}{0.8 0}
\picline{0.1 2}{1.2 0}
}\,,\\
\end{array}
\end{eqnarray*}
\end{theo}

\begin{remark} In fact, $AS$ and $IHX$ are consequences of $STU$, so
they need not to be imposed explicitly here (they are more important
in connection with another $3^\text{rd}$ representation of $\cA$, as
in~\cite[section~5]{Bar-Natan:Vassiliev}).  However, we will use $AS$
to turn every trivalent vertex to be oriented counterclockwise and then
drop all orientation arrows.
\end{remark}

\proof{} This is basically a consequence of the T-shirt identity

\def\CD#1{%
\diag{4mm}{2}{2}{
\piccircle{1 1}{1}{}
%\piclinedash{0.2 0.1}{0.25}
\pictranslate{2 0}{\picscale{-1 1}{\pictranslate{1 1}{#1}}}}}

\[
\CD{\picline{1 125 polar}{1 235 polar}
\picline{1 0 polar}{1 250 polar}
}
\enspace-\enspace
\CD{\picline{1 110 polar}{1 250 polar}
\picline{1 230 polar}{1 0 polar}
}
\enspace=\enspace
\CD{\picline{1 240 polar}{0 0}
\picline{1 0 polar}{0 0}
\picline{1 120 polar}{0 0}}
\enspace=\enspace
%-\,
\CD{\picline{1 125 polar}{1 235 polar}
\picline{1 0 polar}{1 110 polar}
}
\enspace-\enspace
\CD{\picline{1 130 polar}{1 0 polar}
\picline{1 110 polar}{1 250 polar}
}
\,,
\]
(with some more technical details). \qed

\def\fg{{\frak g}}
\def\tr{\text{tr}}

Now given a finite-dimensional Lie-algebra $\fg$ with a metric
and an orthonormal basis $\{\fg_a\}_{a=1}^{\dim\fg}$ and a
finite dimensional representation $R$, set
\[
W_{\fg,R}\left(
\diag{6mm}{5}{2}%
{\picvecline{0 0.5}{5 0.5}
%\piclinedash{0.2 0.1}{0.25}
\picellipsearc{3 0.5}{1.6 1.0}{0 180}
\piccirclearc{1.3 0.5}{0.8}{0 180}
\picline{3.3 0.5}{3.3 1.5}
\pictext{$a$}{1.3 1.4}{-0.5 0}
\pictext{$b$}{2.4 1.5}{-0.5 0}
\pictext{$c$}{3.4 0.7}{0 0}
\pictext{$d$}{4.5 1.0}{0 0}
}
\right)
\quad=\quad\sum_{a,b,c,d=1}^{\dim\fg}\,f_{bcd}\tr_R(\fg_a\fg_b
\fg_a\fg_c\fg_d)\,,
\]
where $f_{bcd}$ are the structure constants of $\fg$ relative to the basis
$\{\fg_a\}$. It should be clear how to extend this example and define
$W_{\fg,R}(D)$ for any diagram $D$ of the kind appearing in
theorem~\ref{TrivalentVertices}.

\begin{propproof} $W_{\fg,R}$ is well defined (i.~e., independent of the
choice of the basis $\{\fg_a\}$) and satisfies:
\begin{itemize}
\item[*] The $AS$ relation by the anti-symmetry of the bracket.
\item[*] The $IHX$ relation because of the Jacobi identity.
\item[*] The $STU$ relation because representations represent.
\end{itemize}
\end{propproof}

\begin{conj}\label{conj_two}
All weight systems ($\entspr$ elements of $\cA^*$)
come from this construction.
\end{conj}

\def\tempp#1#2{%
\vtop{\hbox{#1}\vbox to 0pt{\hbox{small #2}\vss}}}
\def\tempp#1#2{\hbox{\shortstack{#1\\\small\strut #2%
}}}
\def\tempq{%
{\setbox9=\tempp{2}{sec}\vbox to \ht9{\vss\hbox{---}\vss}}}

{\samepage
  \noindent\underline{\bf A word about numbers}
  \begin{center}
  \begin{tabular}{c|*{10}c}
  &\\[-2mm]
  $m$&0&1&2&3&4&5&6&7&8&9\\[2mm]
  \hline
  \\[-2mm]
  $\dim \cA^r_m$&1&0&1&1&3&4&9&14&27&44\\[2mm]
  \hline
  \\[-2mm]
  $\dim \cA_m$&1&1&2&3&6&10&19&33&60&104\\[2mm]
  \hline
  \\[-2mm]
  $\dim\left(\vctext{@{}c@{}}{span of\\all $W_{\fg,R}$}{\tiny}\right)%
  $&1&1&2&3&6&10&19&33&60&104\\[2mm]
  \hline
  \\[-2mm]
  \tempp{CPU time}{190MHz Digital alpha Workstation}&
  \tempq&\tempq&\tempq&\tempq&\tempq&\tempq&
  \tempp{0.64}{sec}&
  \tempp{27}{sec}&
  \tempp{19}{min}&
  \tempp{2.7}{days}\\[1mm]
  \end{tabular}
  \end{center}
}

Looking at this table, the case for conjecture~\ref{conj_two} appears to be
convincing. However,

\noindent{\bf Warning:}
% Conjecture \ref{conjone} (that Vassiliev invariants separate
% knots) and conjecture \ref{conj_two} (that they all come from
% Lie algebras, at least in the stronger form, where only
% semi-simple \& Abelian algebras are allowed) {\bf contradict.} Furthermore,
% all the Lie algebraic weight systems appearing in the table were generated
% using the Lie algebras $so(N)$ and $gl(N)$, which are known not to span the
% space of all Lie algebraic weight systems~\cite{Vogel:Structures}. So funny
% things are {\em known} to happen beyond degree~9.
From \cite{Bar-Natan:Vassiliev} is was known that Conjecture \ref{conj_two},
at least in the somewhat stronger form, where only semi-simple \& Abelian
algebras are allowed, would answer negatively question \ref{quest1}
and therefore contradicts Conjecture \ref{conjone}. Finally,
recently Vogel ~\cite{Vogel:Structures} {\em disproved} this
stronger version of Conjecture \ref{conj_two}. However,
all the Lie algebraic weight systems appearing in the table were generated
using only the Lie algebras $so(N)$ and $gl(N)$. Beyond degree~9
we will have to deal with nilpotent (and ev.~exceptional) Lie algebras
too. But Vogel even announced to the second author that
Conjecture \ref{conj_two} is wrong {\em in full generality}.

Anyway, the answer to question \ref{quest1} and the fate of
Conjecture \ref{conjone} remain unclear.

In a way, this is good news. It means that we don't understand something,
which means that we still have something left to do!

Now let's come back to our Fundamental Theorem. We will use the following

\noindent\underline{Equivalent Reformulation}\\[2mm]
There exists a ``universal Vassiliev invariant''
\[
\tZ\,:\,\bigl\{\text{ knots }\bigr\}\,\longrightarrow
\bar{\cA^r}\qquad
\left(\vcbox{\hbox{\shortstack{the graded completion of $\cA^r$}}}\right)
\]
such that if $D$ is the degree $m$ chord diagram underlying
an $m$-singular knot $K$, then
\[
\tZ(K)\,=\,D+
\left(\vcbox{\hbox{\shortstack{\small higher degree\\
\small diagrams}}}\right)
\]

\proof{of equivalence}
\[
\diag{3mm}{15}{7}
{\pictext{$\bold C$}{7.5 0}{-0.5 0}
\pictext{$\bigl\{\text{ knots }\bigr\}$}{2.0 5.5}{-0.5 0}
\pictext{$\left\{\vcbox{\hbox{\footnotesize \shortstack{chord\\
diagrams}}}\right\}\,=\,\cA^r$}{13.0 5.5}{-0.5 0}
\pictext{$\stackrel{\tZ}{\rightarrowto{3.2\unitlength}}$%
}{7.0 5.5}{-0.7 0}
\pictranslate{2.2 4.8}{\picrotate{-40}{\pictext
{\rightarrowto{5.5\unitlength}}{0 0}{0 0}
}}
\pictranslate{12.0 4.9}{\picrotate{-140}{\pictext
{\rightarrowto{5.0\unitlength}}{0 0}{0 0}
}}
\pictext{\footnotesize $V$}{5.4 2.6}{0 0}
\pictext{\footnotesize $W_V$}{11.3 2.6}{0 0}
}
\]
If you have $\tZ$ and you're given a $W$, define $V$ to be the obvious
composition. If you know how to associate a $V$ to any $W$ in a basis
of $\cA^r$, there's a unique $\tZ$ making the diagram commutative. \qed

Here we will present Kontsevich's geometric approach for constructing
such a $\tZ$. %First of all,

\subsection{Connections, curvature, and holonomy}

Up to some (important, but not here) subtlety, a connection is a $1$-form
whose values are in the algebra of endomorphisms of the fiber. One would
like to know how much of the theory of connections can be generalized to
the case of $1$-forms with values in an arbitrary associative algebra. As
was shown by K-T.~Chen \cite{Chen:Iterated}, much of the theory persists
in the more general case.  Let us briefly review some aspects of Chen's
theory.

Let $X$ be a smooth manifold and let $\frak A$ be a topological algebra
over the real numbers ${\bold R}$ (or the complex numbers ${\bold C}$),
with a unit $1$. An {\em $\frak A$-valued connection} $\Omega$ on $X$
is an $\frak A$-valued $1$-form $\Omega$ on $X$. Its curvature
$F_\Omega$ is the $\frak A$-valued $2$-form
$F_\Omega=d\Omega+\Omega\wedge\Omega$, where the definitions of the
exterior differentiation operator $d$ and of the wedge product $\wedge$
are precisely the same as the corresponding definitions in the case of
matrix valued forms. The notion of ``parallel transport'' also has a
generalization in the new context: Let $B:I\rightarrow X$ be a
smooth map from some interval $I=[a,b]$ to $X$. Define the {\em holonomy}
$hol_B(\Omega)$ of $\Omega$ along $B$ to be the function
$hol_{B}(\Omega):I\rightarrow\frak A$ which satisfies
\[ hol_{B}(\Omega)(a)=1;\qquad 
	\frac{\partial}{\partial t}hol_{B}(\Omega)(t)
	=\Omega\left(\dot{B}(t)\right)hol_{B}(\Omega)(t),\quad (t\in I) \]
if such a function exists and is unique.
In many interesting cases, $hol_{B}(\Omega)$ exists and is given (see {\it e.g.}
\cite{Chen:Iterated}) by the following ``iterated integral'' formula:
\begin{equation} \label{iterated_integral}
	hol_{B}(\Omega)(t)=1+\sum_{m=1}^\infty
	\mathop{\int}_{a\leq t_1\leq\ldots\leq t_m\leq t}
	(B^\ast\Omega)(t_m)\cdot\ldots\cdot(B^\ast\Omega)(t_1).
\end{equation}
(In this formula $B^\ast\Omega$ denotes the pullback of $\Omega$ to $I$ via
$B$). Furthermore, just like in the standard theory of connections,
if $F_\Omega\equiv 0$ (`$\Omega$ is flat'), then $h_{B,\Omega}$ is
invariant under homotopies of $B$ that preserve its endpoints.

In the case of interest for us, $\frak A$ will be the completion of
a graded algebra of finite type over the complex numbers --- the direct
product of the finite dimensional (over ${\bold C}$) homogeneous
components of a graded algebra. The connection $\Omega$ will be
homogeneous of degree $1$. In this case the $m$th term $hol_{B}(\Omega)_m$
in \eqref{iterated_integral} is homogeneous of degree $m$, and there is
no problem with the convergence of the sum there. Also, as each term
lives in a different degree, Chen's theory implies that each term is
invariant under homotopies of $B$ that preserve its endpoints. These
assertions are not very hard to verify directly from the definition of
$hol_{B}(\Omega)_m$ as a multiple integral.

\subsection{The formal Knizhnik-Zamolodchikov connection}
Let ${\cal D}(\nup)$ be the collection of all diagrams made of $n$ ordered
upward pointing arrows, and chords and oriented vertices as in the
definition of ${\cal A}$, with the standard conventions about higher than
trivalent vertices and about the orientation of vertices:
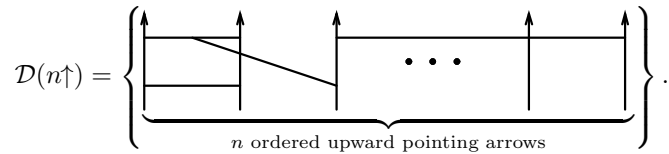
\[ {\cal D}(\nup)=\left\{\raisebox{2mm}{$
	\underbrace{
        \setlength{\unitlength}{0.5\eepiclength}
        \begin{array}{c}  \hspace{-1.7mm}
                \raisebox{-8pt}{\begingroup\makeatletter\ifx\SetFigFont\undefined%
\gdef\SetFigFont#1#2#3#4#5{%
  \reset@font\fontsize{#1}{#2pt}%
  \fontfamily{#3}\fontseries{#4}\fontshape{#5}%
  \selectfont}%
\fi\endgroup%
\begin{picture}(6024,1239)(0,-10)
\thicklines
\put(3912,612){\blacken\ellipse{60}{60}}
\put(3912,612){\ellipse{60}{60}}
\put(3612,612){\blacken\ellipse{60}{60}}
\put(3612,612){\ellipse{60}{60}}
\put(3312,612){\blacken\ellipse{60}{60}}
\put(3312,612){\ellipse{60}{60}}
\path(6012,12)(6012,1212)
\path(6049.500,1077.000)(6012.000,1212.000)(5974.500,1077.000)
\path(4812,12)(4812,1212)
\path(4849.500,1077.000)(4812.000,1212.000)(4774.500,1077.000)
\path(2412,12)(2412,1212)
\path(2449.500,1077.000)(2412.000,1212.000)(2374.500,1077.000)
\path(1212,12)(1212,1212)
\path(1249.500,1077.000)(1212.000,1212.000)(1174.500,1077.000)
\path(12,12)(12,1212)
\path(49.500,1077.000)(12.000,1212.000)(-25.500,1077.000)
\path(12,912)(1212,912)
\path(612,912)(2412,312)
\path(12,312)(1212,312)
\path(2412,912)(6012,912)
\end{picture} }
                \hspace{-1.9mm}
        \end{array}
}_{\text{$n$ ordered
	upward pointing arrows}}$} \right\}. \]
Let the ground field be ${\bold C}$ and
let ${\cal A}(\nup)$ be the quotient
\[ {\cal A}(\nup) = \Span({\cal D}(\nup))\left/
	\{STU \text{ relations}\}\right. . \]
${\cal A}(\nup)$ is an algebra with `composition' as its product:
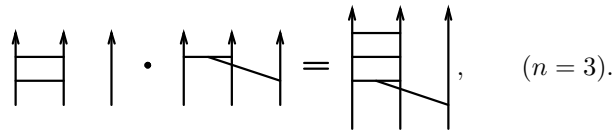
\[ 
        \setlength{\unitlength}{0.5\eepiclength}
        \begin{array}{c}  \hspace{-1.7mm}
                \raisebox{-8pt}{\begingroup\makeatletter\ifx\SetFigFont\undefined%
\gdef\SetFigFont#1#2#3#4#5{%
  \reset@font\fontsize{#1}{#2pt}%
  \fontfamily{#3}\fontseries{#4}\fontshape{#5}%
  \selectfont}%
\fi\endgroup%
\begin{picture}(5424,1539)(0,-10)
\thicklines
\put(1662,807){\blacken\ellipse{60}{60}}
\put(1662,807){\ellipse{60}{60}}
\path(3612,837)(3912,837)
\path(3612,762)(3912,762)
\path(5412,12)(5412,1512)
\path(5449.500,1377.000)(5412.000,1512.000)(5374.500,1377.000)
\path(4812,12)(4812,1512)
\path(4849.500,1377.000)(4812.000,1512.000)(4774.500,1377.000)
\path(4212,12)(4212,1512)
\path(4249.500,1377.000)(4212.000,1512.000)(4174.500,1377.000)
\path(3312,312)(3312,1212)
\path(3349.500,1077.000)(3312.000,1212.000)(3274.500,1077.000)
\path(2712,312)(2712,1212)
\path(2749.500,1077.000)(2712.000,1212.000)(2674.500,1077.000)
\path(2112,312)(2112,1212)
\path(2149.500,1077.000)(2112.000,1212.000)(2074.500,1077.000)
\path(1212,312)(1212,1212)
\path(1249.500,1077.000)(1212.000,1212.000)(1174.500,1077.000)
\path(612,312)(612,1212)
\path(649.500,1077.000)(612.000,1212.000)(574.500,1077.000)
\path(12,312)(12,1212)
\path(49.500,1077.000)(12.000,1212.000)(-25.500,1077.000)
\path(12,912)(612,912)
\path(12,612)(612,612)
\path(2112,912)(2712,912)
\path(2412,912)(3312,612)
\path(4212,1212)(4812,1212)
\path(4212,912)(4812,912)
\path(4212,612)(4812,612)
\path(4512,612)(5412,312)
\end{picture} }
                \hspace{-1.9mm}
        \end{array}
,\qquad (n=3). \]
${\cal A}(\nup)$ is graded by half the number of vertices in a diagram,
excluding the $2n$ endpoints of the $n$ arrows; the degree of the above
product is $4$.

For $1\leq i<j\leq n$ define $\Omega_{ij}\in{\cal A}(\nup)$ by
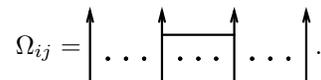
\[ \Omega_{ij}=
        \setlength{\unitlength}{0.5\eepiclength}
        \begin{array}{c}  \hspace{-1.7mm}
                \raisebox{-8pt}{\begingroup\makeatletter\ifx\SetFigFont\undefined%
\gdef\SetFigFont#1#2#3#4#5{%
  \reset@font\fontsize{#1}{#2pt}%
  \fontfamily{#3}\fontseries{#4}\fontshape{#5}%
  \selectfont}%
\fi\endgroup%
\begin{picture}(2724,939)(0,-10)
\thicklines
\put(2487,312){\blacken\ellipse{30}{30}}
\put(2487,312){\ellipse{30}{30}}
\put(2037,312){\blacken\ellipse{30}{30}}
\put(2037,312){\ellipse{30}{30}}
\put(2262,312){\blacken\ellipse{30}{30}}
\put(2262,312){\ellipse{30}{30}}
\put(1587,312){\blacken\ellipse{30}{30}}
\put(1587,312){\ellipse{30}{30}}
\put(1137,312){\blacken\ellipse{30}{30}}
\put(1137,312){\ellipse{30}{30}}
\put(1362,312){\blacken\ellipse{30}{30}}
\put(1362,312){\ellipse{30}{30}}
\put(687,312){\blacken\ellipse{30}{30}}
\put(687,312){\ellipse{30}{30}}
\put(462,312){\blacken\ellipse{30}{30}}
\put(462,312){\ellipse{30}{30}}
\put(237,312){\blacken\ellipse{30}{30}}
\put(237,312){\ellipse{30}{30}}
\path(2712,12)(2712,912)
\path(2749.500,777.000)(2712.000,912.000)(2674.500,777.000)
\path(1812,12)(1812,912)
\path(1849.500,777.000)(1812.000,912.000)(1774.500,777.000)
\path(912,12)(912,912)
\path(949.500,777.000)(912.000,912.000)(874.500,777.000)
\path(12,12)(12,912)
\path(49.500,777.000)(12.000,912.000)(-25.500,777.000)
\path(912,612)(1812,612)
\end{picture} }
                \hspace{-1.9mm}
        \end{array}
. \]
Let $X_n$ be the configuration space of $n$ distinct points in ${\bold C}$;
$X_n=\{(z_1,\ldots,z_n)\in{\bold C}^n: z_i=z_j\Rightarrow i=j\}$, and let
$\omega_{ij}$ be the complex $1$-form on $X_n$ defined by
\[ \omega_{ij}=d(\log z_i-z_j)=\frac{dz_i-dz_j}{z_i-z_j}. \]
The formal Knizhnik-Zamolodchikov connection is the ${\cal
A}(\nup)$-valued connection $\displaystyle\Omega_n=\sum_{1\leq i<j\leq
n}\Omega_{ij}\omega_{ij}$ on $X_n$.

\begin{prop} The formal Knizhnik-Zamolodchikov connection $\Omega_n$
is flat.
\end{prop}

\proof{} Clearly $d\Omega_n=0$. Let us check that
\begin{equation} \Omega_n\wedge\Omega_n=\sum_{i<j;i'<j'}
	\Omega_{ij}\Omega_{i'j'} \omega_{ij}\wedge\omega_{i'j'} = 0.
\label{omega_sqr} \end{equation}
The above sum can be separated into three parts, according to the
cardinality of the set $\{i,j,i',j'\}$. If this cardinality is $2$ or $4$
then $\Omega_{ij}$ and $\Omega_{i'j'}$ commute, while $\omega_{ij}$ and
$\omega_{i'j'}$ anti-commute. It is easy to check that this implies that
the corresponding parts of the sum \eqref{omega_sqr} vanish. The only
interesting case is when $|\{i,j,i',j'\}|=3$, say
$\{i,j,i',j'\}=\{1,2,3\}$. In this case,
\[ \sum_{\{i,j,i',j'\}=\{1,2,3\}} \hspace{-12pt} \Omega_{ij}\Omega_{i'j'}
	\omega_{ij}\wedge\omega_{i'j'}
	= (\Omega_{12}\Omega_{23}-\Omega_{23}\Omega_{12})
		\omega_{12}\wedge\omega_{23}
	+ (\text{cyclic permutations}). \]
By the $STU$ relation this is 
\begin{equation} \label{KZSTU}
	= \Omega_{123}(\omega_{12}\wedge\omega_{23} + 
	(\text{cyclic permutations})) = 0, \end{equation}
where $\Omega_{123}$ is given by
\[ \Omega_{123}= 
        \setlength{\unitlength}{0.5\eepiclength}
        \begin{array}{c}  \hspace{-1.7mm}
                \raisebox{-8pt}{\begingroup\makeatletter\ifx\SetFigFont\undefined%
\gdef\SetFigFont#1#2#3#4#5{%
  \reset@font\fontsize{#1}{#2pt}%
  \fontfamily{#3}\fontseries{#4}\fontshape{#5}%
  \selectfont}%
\fi\endgroup%
\begin{picture}(2424,939)(0,-10)
\thicklines
\put(1812,462){\blacken\ellipse{30}{30}}
\put(1812,462){\ellipse{30}{30}}
\put(2112,462){\blacken\ellipse{30}{30}}
\put(2112,462){\ellipse{30}{30}}
\put(1542,462){\blacken\ellipse{30}{30}}
\put(1542,462){\ellipse{30}{30}}
\path(12,12)(12,912)
\path(49.500,777.000)(12.000,912.000)(-25.500,777.000)
\path(612,12)(612,912)
\path(649.500,777.000)(612.000,912.000)(574.500,777.000)
\path(1212,12)(1212,912)
\path(1249.500,777.000)(1212.000,912.000)(1174.500,777.000)
\path(2412,12)(2412,912)
\path(2449.500,777.000)(2412.000,912.000)(2374.500,777.000)
\path(612,312)(312,612)(1212,312)
\path(312,612)(12,312)
\end{picture} }
                \hspace{-1.9mm}
        \end{array}
 
	\in{\cal A}(\nup). \]
The vanishing of
$\omega_{12}\wedge\omega_{23}+(\text{cyclic permutations})$ is called
`Arnold's identity' \cite{Arnold:DyedBraids} and can be easily verified
by a direct computation. \qed

\begin{remark} The connection $\Omega_n$ has a simple generalization to
the case when the underlying algebra is ${\cal A}(\nup\!\!\amalg\ndown)$,
the algebra generated by diagrams having $2n$ arrows, whose first $n$
arrows point upward and whose next $n$ arrows point downward.  The only
difference is a sign difference in the application of the $STU$ relation
in \eqref{KZSTU}. Therefore if one defines
\[ \Omega_{n,n}=\sum_{1\leq i\leq j\leq 2n}s_is_j\Omega_{ij}\omega_{ij}, \]
where
$s_i=\left\{\begin{array}{cl}
  +1 & \,i\leq n \\
  -1 & \,i>n
\end{array}\right.$,
then the connection $\Omega_{n,n}$ is flat.
\end{remark}

\subsection{Kontsevich's integral invariants}
Choose a decomposition ${\bold R}^3={\bold C}_z\times{\bold R}_t$ of
${\bold R}^3$ to a product of a complex plane ${\bold C}_z$ parametrized
by $z$ and a real line ${\bold R}_t$ parametrized by $t$ and let
$K:S^1\rightarrow{\bold R}^3$ be a parametrized knot on which the function
$t$ is a Morse function. Consider the following series, whose precise
definition will be discussed below:
\begin{equation} \label{ZDef}
  Z(K)=\sum_{m=0}^\infty (2\pi i)^{-m}
	\mathop{\int}_{t_1<\ldots<t_m}
	\sum_{\begin{array}{c}
		\text{\scriptsize applicable pairings} \\
		\scriptstyle P=\{(z_i,z'_i)\}
	\end{array}}
	(-1)^{\#P_{\downarrow}}D_P
	\bigwedge_{i=1}^{m}\frac{dz_i-dz'_i}{z_i-z'_i}
	\in{\bar{\cal A}}^r. \label{Kints}
\end{equation}

\par
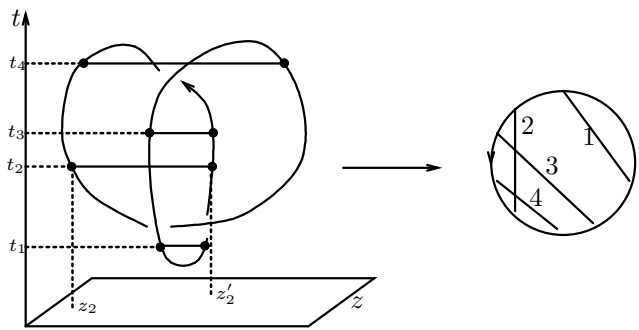
\begin{figure}[htpb]
\[ 
        \setlength{\unitlength}{0.5\eepiclength}
        \begin{array}{c}  \hspace{-1.7mm}
                \raisebox{-8pt}{\begingroup\makeatletter\ifx\SetFigFont\undefined%
\gdef\SetFigFont#1#2#3#4#5{%
  \reset@font\fontsize{#1}{#2pt}%
  \fontfamily{#3}\fontseries{#4}\fontshape{#5}%
  \selectfont}%
\fi\endgroup%
\begin{picture}(7854,3939)(0,-10)
\thicklines
\put(6945.000,2052.000){\arc{1800.000}{6.2832}{9.4248}}
\path(6007.500,2187.000)(6045.000,2052.000)(6082.500,2187.000)
\put(6945.000,2052.000){\arc{1800.000}{3.1416}{6.2832}}
\put(960,3297){\blacken\ellipse{90}{90}}
\put(960,3297){\ellipse{90}{90}}
\put(2580,2427){\blacken\ellipse{90}{90}}
\put(2580,2427){\ellipse{90}{90}}
\put(1785,2427){\blacken\ellipse{90}{90}}
\put(1785,2427){\ellipse{90}{90}}
\put(2565,2007){\blacken\ellipse{90}{90}}
\put(2565,2007){\ellipse{90}{90}}
\put(810,2007){\blacken\ellipse{90}{90}}
\put(810,2007){\ellipse{90}{90}}
\put(2475,1017){\blacken\ellipse{90}{90}}
\put(2475,1017){\ellipse{90}{90}}
\put(1920,1002){\blacken\ellipse{90}{90}}
\put(1920,1002){\ellipse{90}{90}}
\put(3465,3297){\blacken\ellipse{90}{90}}
\put(3465,3297){\ellipse{90}{90}}
\path(240,12)(240,3912)
\path(277.500,3777.000)(240.000,3912.000)(202.500,3777.000)
\path(240,12)(3765,12)
\path(240,12)(1065,612)(4590,612)(3765,12)
\path(4200,1992)(5385,1992)
\path(5250.000,1954.500)(5385.000,1992.000)(5250.000,2029.500)
\dottedline{90}(1920,1002)(240,1002)
\dottedline{90}(825,1992)(825,252)
\dottedline{90}(2550,1992)(2550,417)
\dottedline{90}(810,2007)(240,2007)
\dottedline{90}(1800,2427)(240,2427)
\dottedline{90}(930,3297)(240,3297)
\path(1005,3297)(3465,3297)
\path(1815,2427)(2565,2427)
\path(825,2007)(2535,2007)
\path(1965,1017)(2505,1017)
\path(6345,2727)(6345,1452)
\path(6120,2427)(7320,1302)
\path(6945,2952)(7770,1827)
\path(6120,1827)(6870,1227)
\path(2055,1257)	(2102.932,1260.104)
	(2149.173,1263.242)
	(2193.767,1266.421)
	(2236.752,1269.646)
	(2278.172,1272.926)
	(2318.066,1276.266)
	(2356.477,1279.674)
	(2393.445,1283.157)
	(2463.218,1290.374)
	(2527.714,1297.973)
	(2587.264,1306.007)
	(2642.198,1314.533)
	(2692.843,1323.604)
	(2739.531,1333.276)
	(2782.590,1343.603)
	(2822.351,1354.642)
	(2893.295,1379.070)
	(2955.000,1407.000)

\path(2955,1407)	(3002.078,1432.227)
	(3052.920,1461.504)
	(3106.714,1494.581)
	(3162.647,1531.210)
	(3219.908,1571.141)
	(3277.683,1614.125)
	(3335.160,1659.912)
	(3391.526,1708.253)
	(3445.970,1758.898)
	(3497.679,1811.599)
	(3545.841,1866.105)
	(3589.642,1922.168)
	(3628.271,1979.539)
	(3660.915,2037.967)
	(3686.763,2097.204)
	(3705.000,2157.000)

\path(3705,2157)	(3713.712,2197.858)
	(3720.853,2239.871)
	(3726.432,2282.934)
	(3730.462,2326.938)
	(3732.951,2371.778)
	(3733.909,2417.345)
	(3733.347,2463.533)
	(3731.276,2510.235)
	(3727.704,2557.343)
	(3722.643,2604.752)
	(3716.103,2652.353)
	(3708.093,2700.040)
	(3698.624,2747.706)
	(3687.706,2795.244)
	(3675.349,2842.546)
	(3661.564,2889.506)
	(3646.360,2936.017)
	(3629.748,2981.972)
	(3611.737,3027.263)
	(3592.339,3071.784)
	(3571.563,3115.428)
	(3549.419,3158.087)
	(3525.918,3199.655)
	(3501.070,3240.025)
	(3474.884,3279.089)
	(3447.372,3316.741)
	(3388.407,3387.380)
	(3324.256,3451.085)
	(3255.000,3507.000)

\path(3255,3507)	(3188.884,3542.511)
	(3152.244,3554.606)
	(3113.757,3563.309)
	(3073.826,3568.905)
	(3032.856,3571.679)
	(2991.251,3571.916)
	(2949.416,3569.903)
	(2907.756,3565.923)
	(2866.675,3560.263)
	(2826.577,3553.209)
	(2787.868,3545.045)
	(2716.231,3526.529)
	(2655.000,3507.000)

\path(2655,3507)	(2609.211,3489.384)
	(2560.240,3467.609)
	(2508.738,3441.998)
	(2455.357,3412.875)
	(2400.746,3380.562)
	(2345.557,3345.382)
	(2290.440,3307.659)
	(2236.046,3267.716)
	(2183.026,3225.876)
	(2132.029,3182.461)
	(2083.708,3137.794)
	(2038.712,3092.200)
	(1997.693,3046.000)
	(1961.301,2999.518)
	(1930.186,2953.077)
	(1905.000,2907.000)

\path(1905,2907)	(1885.363,2862.777)
	(1867.650,2816.034)
	(1851.795,2766.946)
	(1837.734,2715.687)
	(1825.399,2662.430)
	(1814.725,2607.352)
	(1805.648,2550.625)
	(1798.100,2492.425)
	(1792.016,2432.926)
	(1787.330,2372.301)
	(1783.978,2310.727)
	(1781.892,2248.377)
	(1781.007,2185.425)
	(1781.258,2122.046)
	(1782.578,2058.414)
	(1784.902,1994.704)
	(1788.165,1931.090)
	(1792.300,1867.746)
	(1797.242,1804.847)
	(1802.926,1742.568)
	(1809.284,1681.082)
	(1816.252,1620.564)
	(1823.764,1561.189)
	(1831.754,1503.131)
	(1840.157,1446.564)
	(1848.906,1391.662)
	(1857.936,1338.601)
	(1867.181,1287.554)
	(1876.576,1238.696)
	(1886.055,1192.202)
	(1895.551,1148.245)
	(1905.000,1107.000)

\path(1905,1107)	(1923.885,1034.412)
	(1937.355,991.699)
	(1953.844,947.809)
	(1973.602,905.067)
	(1996.880,865.801)
	(2055.000,807.000)

\path(2055,807)	(2124.119,779.499)
	(2163.824,772.624)
	(2205.000,770.333)
	(2246.176,772.624)
	(2285.881,779.499)
	(2355.000,807.000)

\path(2355,807)	(2399.928,845.297)
	(2437.560,903.375)
	(2454.739,942.029)
	(2471.412,988.266)
	(2488.019,1042.963)
	(2505.000,1107.000)

\path(2505,1407)	(2514.490,1477.044)
	(2523.367,1544.673)
	(2531.633,1609.949)
	(2539.286,1672.933)
	(2546.327,1733.688)
	(2552.755,1792.274)
	(2558.572,1848.754)
	(2563.776,1903.190)
	(2568.367,1955.643)
	(2572.347,2006.175)
	(2575.714,2054.848)
	(2578.470,2101.724)
	(2580.612,2146.865)
	(2582.143,2190.332)
	(2583.061,2232.187)
	(2583.367,2272.492)
	(2583.061,2311.309)
	(2582.143,2348.700)
	(2578.470,2419.449)
	(2572.347,2485.234)
	(2563.776,2546.549)
	(2552.755,2603.889)
	(2539.286,2657.748)
	(2523.367,2708.620)
	(2505.000,2757.000)

\path(2505,2757)	(2471.299,2816.829)
	(2414.753,2882.085)
	(2375.715,2918.945)
	(2328.330,2959.798)
	(2271.718,3005.524)
	(2205.000,3057.000)

\path(2335.229,3005.311)(2205.000,3057.000)(2289.912,2945.550)
\path(1905,3207)	(1853.524,3273.718)
	(1807.798,3330.330)
	(1766.945,3377.715)
	(1730.085,3416.753)
	(1664.829,3473.299)
	(1605.000,3507.000)

\path(1605,3507)	(1536.744,3523.134)
	(1496.265,3526.367)
	(1454.040,3526.804)
	(1411.935,3524.765)
	(1371.816,3520.572)
	(1305.000,3507.000)

\path(1305,3507)	(1251.589,3487.885)
	(1191.063,3460.948)
	(1126.573,3427.472)
	(1061.273,3388.744)
	(998.313,3346.047)
	(940.846,3300.668)
	(892.024,3253.891)
	(855.000,3207.000)

\path(855,3207)	(830.635,3164.930)
	(808.285,3117.671)
	(787.945,3065.929)
	(769.612,3010.412)
	(753.283,2951.826)
	(738.954,2890.879)
	(726.624,2828.276)
	(716.287,2764.725)
	(707.942,2700.933)
	(701.585,2637.606)
	(697.211,2575.451)
	(694.820,2515.176)
	(694.406,2457.486)
	(695.967,2403.089)
	(699.499,2352.691)
	(705.000,2307.000)

\path(705,2307)	(720.517,2237.739)
	(732.367,2199.407)
	(746.587,2159.240)
	(762.912,2117.713)
	(781.075,2075.302)
	(800.807,2032.480)
	(821.842,1989.724)
	(843.914,1947.507)
	(866.755,1906.305)
	(890.098,1866.593)
	(913.676,1828.845)
	(960.470,1761.144)
	(1005.000,1707.000)

\path(1005,1707)	(1050.679,1662.066)
	(1105.698,1616.536)
	(1172.254,1569.091)
	(1210.545,1544.238)
	(1252.545,1518.413)
	(1298.528,1491.449)
	(1348.768,1463.183)
	(1403.540,1433.449)
	(1463.120,1402.083)
	(1527.781,1368.920)
	(1597.798,1333.795)
	(1634.901,1315.445)
	(1673.446,1296.543)
	(1713.468,1277.068)
	(1755.000,1257.000)

\put(4290,237){\makebox(0,0)[lb]{\smash{{\mathmode{z}}}}}
\put(60,3762){\makebox(0,0)[lb]{\smash{{\mathmode{t}}}}}
\put(7185,2277){\makebox(0,0)[lb]{\smash{{\mathmode{1}}}}}
\put(6420,2412){\makebox(0,0)[lb]{\smash{{\mathmode{2}}}}}
\put(6720,1917){\makebox(0,0)[lb]{\smash{{\mathmode{3}}}}}
\put(6525,1512){\makebox(0,0)[lb]{\smash{{\mathmode{4}}}}}
\put(30,3237){\makebox(0,0)[lb]{\smash{{\mathmode{\scriptstyle t_4}}}}}
\put(15,2382){\makebox(0,0)[lb]{\smash{{\mathmode{\scriptstyle t_3}}}}}
\put(0,1962){\makebox(0,0)[lb]{\smash{{\mathmode{\scriptstyle t_2}}}}}
\put(30,957){\makebox(0,0)[lb]{\smash{{\mathmode{\scriptstyle t_1}}}}}
\put(885,207){\makebox(0,0)[lb]{\smash{{\mathmode{\scriptstyle z_2}}}}}
\put(2610,357){\makebox(0,0)[lb]{\smash{{\mathmode{\scriptstyle z'_2}}}}}
\end{picture} }
                \hspace{-1.9mm}
        \end{array}
 \]
\caption{$m=4$: a knot $K$ with a pairing $P$ and the corresponding
	chord diagram $D_P$. Notice that $D_P=0$ in ${\bar{\cal A}}^r$
	due to the isolated chord marked by $1$.}
\label{DP} \end{figure}

\noindent In the above equation,
\begin{itemize}
\item an `applicable pairing' is a choice of an unordered
	pair $(z_i,z'_i)$ for every $1\leq i\leq m$, for which
	$(z_i,t_i)$ and $(z'_i,t_i)$ are {\em distinct} points on $K$.
\item $\#P_{\downarrow}$ is the number of points of the form
	$(z_i,t_i)$ or $(z'_i,t_i)$ at which $K$ is decreasing. Remember
	that in this article we are only considering {\em oriented} knots.
\item $D_P$ is the chord diagram naturally associated with $K$ and $P$ as
	in figure~\ref{DP}. It is to be regarded as an element of
	${\bar{\cal A}}^r$.
\item every pairing defines a map $\{t_i\}\mapsto\{(z_i,z'_i)\}$ locally
	around the current values of the $t_i$'s. Use this map to pull the
	$dz_i$'s and $dz'_i$'s to the $m$-simplex $t_{\text{\scriptsize
	min}}< t_1<\ldots<t_m<t_{\text{\scriptsize max}}$ (where
	$t_{\text{\scriptsize min}}$ ($t_{\text{\scriptsize max}}$) is
	the minimal (maximal) value of $t$ on $K$) and then integrate
	the indicated wedge product over that simplex.
\end{itemize}

\subsubsection{Finiteness} \label{KZfiniteness}
Properly interpreted, the integrals in \eqref{Kints} are finite. There appears
to be a problem in the denominator when $z_i-z'_i$ is small for some $i$.
This can happen in either of two ways:
\begin{enumerate}
\item $\quad
        \setlength{\unitlength}{0.6\eepiclength}
        \begin{array}{c}  \hspace{-1.7mm}
                \raisebox{-8pt}{\begingroup\makeatletter\ifx\SetFigFont\undefined%
\gdef\SetFigFont#1#2#3#4#5{%
  \reset@font\fontsize{#1}{#2pt}%
  \fontfamily{#3}\fontseries{#4}\fontshape{#5}%
  \selectfont}%
\fi\endgroup%
\begin{picture}(995,931)(0,-10)
\thicklines
\path(278,818)(983,818)
\path(68,563)(848,563)
\path(908,8)	(905.052,78.315)
	(901.745,143.726)
	(898.025,204.481)
	(893.840,260.825)
	(889.137,313.008)
	(883.862,361.275)
	(877.962,405.875)
	(871.385,447.054)
	(855.985,520.138)
	(837.237,582.505)
	(814.717,636.133)
	(788.000,683.000)

\path(788,683)	(729.612,757.403)
	(691.192,795.074)
	(648.327,830.386)
	(602.297,861.338)
	(554.380,885.926)
	(505.854,902.148)
	(458.000,908.000)

\path(458,908)	(407.931,902.148)
	(357.640,885.926)
	(308.229,861.338)
	(260.803,830.386)
	(216.463,795.074)
	(176.314,757.403)
	(113.000,683.000)

\path(113,683)	(89.048,636.133)
	(68.687,582.505)
	(51.667,520.138)
	(37.741,447.054)
	(31.861,405.875)
	(26.660,361.275)
	(22.109,313.008)
	(18.175,260.825)
	(14.829,204.481)
	(12.038,143.726)
	(9.772,78.315)
	(8.000,8.000)

\path(-26.647,143.760)(8.000,8.000)(48.337,142.180)
\put(233,713){\makebox(0,0)[lb]{\smash{{\mathmode{\scriptstyle z_{i+1}}}}}}
\put(83,383){\makebox(0,0)[lb]{\smash{{\mathmode{\scriptstyle z_i}}}}}
\put(608,383){\makebox(0,0)[lb]{\smash{{\mathmode{\scriptstyle z'_i}}}}}
\end{picture} }
                \hspace{-1.9mm}
        \end{array}
 \ \ \ 
	\parbox{4in}{in this case the integration
	domain for $z_{i+1}$ is as small as $z_i-z'_i$, and its `smallness'
	cancels the singularity coming from the denominator.}$
\item $\quad
        \setlength{\unitlength}{0.6\eepiclength}
        \begin{array}{c}  \hspace{-1.7mm}
                \raisebox{-8pt}{\begingroup\makeatletter\ifx\SetFigFont\undefined%
\gdef\SetFigFont#1#2#3#4#5{%
  \reset@font\fontsize{#1}{#2pt}%
  \fontfamily{#3}\fontseries{#4}\fontshape{#5}%
  \selectfont}%
\fi\endgroup%
\begin{picture}(931,931)(0,-10)
\thicklines
\path(68,563)(848,563)
\path(908,8)	(906.215,78.315)
	(903.934,143.726)
	(901.127,204.481)
	(897.763,260.825)
	(893.812,313.008)
	(889.245,361.275)
	(884.031,405.875)
	(878.139,447.054)
	(864.203,520.138)
	(847.195,582.505)
	(826.875,636.133)
	(803.000,683.000)

\path(803,683)	(738.817,757.403)
	(698.354,795.074)
	(653.832,830.386)
	(606.395,861.338)
	(557.182,885.926)
	(507.336,902.148)
	(458.000,908.000)

\path(458,908)	(408.693,902.148)
	(359.053,885.926)
	(310.396,861.338)
	(264.042,830.386)
	(221.308,795.074)
	(183.511,757.403)
	(128.000,683.000)

\path(128,683)	(99.174,636.133)
	(75.133,582.505)
	(55.459,520.138)
	(39.736,447.054)
	(33.226,405.875)
	(27.546,361.275)
	(22.645,313.008)
	(18.471,260.825)
	(14.972,204.481)
	(12.095,143.726)
	(9.789,78.315)
	(8.000,8.000)

\path(-26.626,143.766)(8.000,8.000)(48.357,142.173)
\put(83,383){\makebox(0,0)[lb]{\smash{{\mathmode{\scriptstyle z_i}}}}}
\put(608,383){\makebox(0,0)[lb]{\smash{{\mathmode{\scriptstyle z'_i}}}}}
\end{picture} }
                \hspace{-1.9mm}
        \end{array}
 \ \ \ 
	\parbox{4in}{in this case the
	corresponding diagram $D_P$ has an isolated chord, and
	so it is $0$ in ${\bar{\cal A}}^r$.}$
\end{enumerate}

\subsubsection{Invariance under horizontal deformations}
For times $t_{\text{\scriptsize min}}\leq a<b\leq t_{\text{\scriptsize max}}$
define $Z(K,[a,b])$ in exactly the same way as \eqref{Kints}, only
restricting the domain of integration to be $a<t_1<\ldots<t_m<b$. Of course,
$Z(K,[a,b])$ will not be in ${\bar{\cal A}}^r$, but rather in the
completed vector space
\[ {\bar{\cal A}}^{K,[a,b]}=
	\Span\left\{\parbox{1.55in}{diagrams whose solid lines are
	as in the part of $K$ on which $a\leq t\leq b$} \right\} \left/
	\left\{\parbox{1.5in}{$STU$ relations and diagrams with
	subdiagrams like $
        \setlength{\unitlength}{0.5\eepiclength}
        \begin{array}{c}  \hspace{-1.7mm}
                \raisebox{-8pt}{%
\begingroup\makeatletter\ifx\SetFigFont\undefined%
\gdef\SetFigFont#1#2#3#4#5{%
  \reset@font\fontsize{#1}{#2pt}%
  \fontfamily{#3}\fontseries{#4}\fontshape{#5}%
  \selectfont}%
\fi\endgroup%
\begin{picture}(2124,336)(0,-10)
\thicklines
\put(1062.000,-175.500){\arc{975.000}{3.5364}{5.8884}}
\path(1512,12)(2112,12)
\path(1977.000,-25.500)(2112.000,12.000)(1977.000,49.500)
\path(12,12)(612,12)
\path(612,12)(1512,12)
\end{picture}
 }
                \hspace{-1.9mm}
        \end{array}
$}\right\}\right..
\]
For example, if $t_1$, $t_4$, and $K$ are as in figure~\ref{DP}, then
the following is a diagram in ${\bar{\cal A}}^{K,[t_1,t_4]}$:
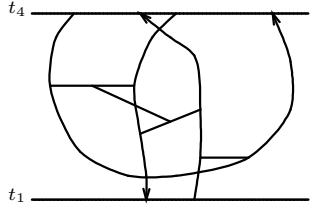
\[ 
        \setlength{\unitlength}{0.5\eepiclength}
        \begin{array}{c}  \hspace{-1.7mm}
                \raisebox{-8pt}{\begingroup\makeatletter\ifx\SetFigFont\undefined%
\gdef\SetFigFont#1#2#3#4#5{%
  \reset@font\fontsize{#1}{#2pt}%
  \fontfamily{#3}\fontseries{#4}\fontshape{#5}%
  \selectfont}%
\fi\endgroup%
\begin{picture}(3762,2520)(0,-10)
\thicklines
\dottedline{48}(3750,36)(300,36)
\dottedline{48}(300,2361)(3750,2361)
\path(525,1461)(1575,1461)
\path(2025,1011)(1050,1461)
\path(2400,1161)(1650,861)
\path(2400,561)(3000,561)
\path(2325,36)	(2333.975,90.800)
	(2342.250,141.772)
	(2349.852,189.106)
	(2356.810,232.996)
	(2363.150,273.635)
	(2368.900,311.213)
	(2378.741,377.959)
	(2386.552,434.772)
	(2392.552,483.192)
	(2396.962,524.755)
	(2400.000,561.000)

\path(2400,561)	(2403.111,620.549)
	(2405.181,693.360)
	(2405.849,733.357)
	(2406.286,775.010)
	(2406.501,817.767)
	(2406.503,861.075)
	(2406.301,904.380)
	(2405.906,947.130)
	(2405.327,988.772)
	(2404.573,1028.753)
	(2402.579,1101.517)
	(2400.000,1161.000)

\path(2400,1161)	(2398.083,1227.919)
	(2397.769,1267.075)
	(2397.589,1309.243)
	(2397.321,1353.835)
	(2396.742,1400.262)
	(2395.630,1447.936)
	(2393.764,1496.269)
	(2390.920,1544.671)
	(2386.878,1592.554)
	(2381.414,1639.330)
	(2374.306,1684.410)
	(2365.334,1727.206)
	(2354.273,1767.128)
	(2325.000,1836.000)

\path(2325,1836)	(2266.207,1901.306)
	(2225.581,1932.465)
	(2181.296,1962.217)
	(2136.225,1990.255)
	(2093.237,2016.270)
	(2025.000,2061.000)

\path(2025,2061)	(1969.456,2104.511)
	(1934.880,2131.942)
	(1894.301,2164.328)
	(1846.623,2202.547)
	(1790.745,2247.480)
	(1725.571,2300.005)
	(1689.154,2329.389)
	(1650.000,2361.000)

\path(1778.588,2305.355)(1650.000,2361.000)(1731.467,2247.007)
\path(2100,2361)	(2036.091,2304.602)
	(1981.485,2255.361)
	(1935.304,2212.398)
	(1896.668,2174.835)
	(1838.516,2112.392)
	(1800.000,2061.000)

\path(1800,2061)	(1767.148,2004.416)
	(1732.245,1934.456)
	(1714.557,1895.728)
	(1696.993,1855.178)
	(1679.767,1813.312)
	(1663.091,1770.638)
	(1647.178,1727.662)
	(1632.241,1684.893)
	(1618.491,1642.836)
	(1606.141,1602.000)
	(1595.405,1562.891)
	(1586.494,1526.017)
	(1575.000,1461.000)

\path(1575,1461)	(1572.526,1398.786)
	(1576.878,1324.688)
	(1581.076,1284.469)
	(1586.336,1242.822)
	(1592.442,1200.262)
	(1599.180,1157.303)
	(1606.334,1114.459)
	(1613.689,1072.245)
	(1621.030,1031.176)
	(1628.142,991.765)
	(1640.819,919.978)
	(1650.000,861.000)

\path(1650,861)	(1657.474,809.030)
	(1667.273,745.494)
	(1678.506,674.225)
	(1684.382,636.888)
	(1690.283,599.055)
	(1696.096,561.204)
	(1701.711,523.815)
	(1711.900,452.338)
	(1719.960,388.456)
	(1725.000,336.000)

\path(1725,336)	(1727.059,291.831)
	(1727.745,231.795)
	(1727.573,193.630)
	(1727.059,148.862)
	(1726.201,96.611)
	(1725.000,36.000)

\path(1690.355,171.761)(1725.000,36.000)(1765.338,170.179)
\path(825,2361)	(775.391,2304.579)
	(733.238,2255.321)
	(697.883,2212.349)
	(668.666,2174.783)
	(626.011,2112.352)
	(600.000,2061.000)

\path(600,2061)	(579.837,2001.674)
	(561.570,1929.899)
	(553.332,1890.588)
	(545.792,1849.659)
	(539.022,1807.611)
	(533.096,1764.941)
	(528.090,1722.148)
	(524.076,1679.731)
	(521.130,1638.186)
	(519.325,1598.013)
	(519.435,1523.774)
	(525.000,1461.000)

\path(525,1461)	(533.010,1416.898)
	(544.163,1368.558)
	(558.228,1316.677)
	(574.971,1261.949)
	(594.159,1205.068)
	(615.559,1146.729)
	(638.939,1087.627)
	(664.065,1028.456)
	(690.705,969.912)
	(718.625,912.689)
	(747.594,857.482)
	(777.377,804.986)
	(807.742,755.895)
	(838.456,710.903)
	(900.000,636.000)

\path(900,636)	(953.947,588.679)
	(1022.412,541.537)
	(1060.761,518.444)
	(1101.145,495.888)
	(1143.032,474.034)
	(1185.892,453.045)
	(1229.193,433.086)
	(1272.404,414.321)
	(1314.992,396.914)
	(1356.427,381.030)
	(1396.176,366.832)
	(1433.709,354.485)
	(1500.000,336.000)

\path(1500,336)	(1560.774,325.701)
	(1634.069,319.500)
	(1674.090,317.756)
	(1715.656,316.820)
	(1758.239,316.619)
	(1801.309,317.081)
	(1844.338,318.135)
	(1886.798,319.708)
	(1928.159,321.727)
	(1967.895,324.121)
	(2040.371,329.745)
	(2100.000,336.000)

\path(2100,336)	(2143.100,341.517)
	(2191.772,348.201)
	(2245.197,356.047)
	(2302.559,365.048)
	(2363.038,375.198)
	(2425.816,386.491)
	(2490.074,398.922)
	(2554.995,412.485)
	(2619.760,427.173)
	(2683.550,442.980)
	(2745.547,459.901)
	(2804.933,477.930)
	(2860.890,497.060)
	(2912.599,517.286)
	(2959.242,538.601)
	(3000.000,561.000)

\path(3000,561)	(3068.940,608.922)
	(3106.031,638.601)
	(3144.258,671.514)
	(3183.164,707.228)
	(3222.291,745.309)
	(3261.182,785.325)
	(3299.378,826.841)
	(3336.421,869.425)
	(3371.854,912.643)
	(3405.219,956.062)
	(3436.058,999.248)
	(3463.914,1041.768)
	(3488.328,1083.189)
	(3508.843,1123.078)
	(3525.000,1161.000)

\path(3525,1161)	(3541.973,1223.454)
	(3551.976,1297.410)
	(3554.695,1337.464)
	(3556.068,1378.909)
	(3556.230,1421.250)
	(3555.311,1463.993)
	(3553.446,1506.641)
	(3550.765,1548.702)
	(3547.403,1589.679)
	(3543.491,1629.079)
	(3534.549,1701.164)
	(3525.000,1761.000)

\path(3525,1761)	(3515.596,1805.041)
	(3502.161,1854.412)
	(3484.037,1910.871)
	(3460.564,1976.175)
	(3446.615,2012.693)
	(3431.083,2052.082)
	(3413.883,2094.561)
	(3394.935,2140.350)
	(3374.154,2189.669)
	(3351.460,2242.737)
	(3326.769,2299.774)
	(3300.000,2361.000)

\path(3388.612,2252.468)(3300.000,2361.000)(3319.941,2222.315)
\put(0,2361){\makebox(0,0)[lb]{\smash{{\mathmode{\scriptstyle t_4}}}}}
\put(0,36){\makebox(0,0)[lb]{\smash{{\mathmode{\scriptstyle t_1}}}}}
\end{picture} }
                \hspace{-1.9mm}
        \end{array}
. \]
The same reasoning as in section \ref{KZfiniteness} shows that $Z(K,[a,b])$
is finite. For
$t_{\text{\scriptsize min}}\leq a<b<c\leq t_{\text{\scriptsize max}}$,
there is an obvious product
${\bar{\cal A}}^{K,[a,b]}\otimes{\bar{\cal A}}^{K,[b,c]}
	\rightarrow{\bar{\cal A}}^{K,[a,c]}$,
and it is easy to show that with this product
$Z(K,[a,b])Z(K,[b,c])=Z(K,[a,c])$.

Let $t_{\text{\scriptsize min}}< a<b< t_{\text{\scriptsize max}}$ be times
for which $K$ has no critical points in the time slice $a\leq t\leq b$, and
let $n$ be the number of upward (or downward) pointing strands of $K$ in
that slice. Then
${\bar{\cal A}}^{K,[a,b]}\equiv{\bar{\cal A}}(\nup\!\!\amalg\ndown)$, and comparing
with \eqref{iterated_integral} and the definition of $\Omega_{n,n}$ we see
that $Z(K,[a,b])$ is the holonomy of $\Omega_{n,n}$ along the braid defined
by the intersection of $K$ with the slice $a\leq t\leq b$. The flatness of
$\Omega_{n,n}$ implies that this holonomy is invariant under horizontal
deformations of that piece of $K$, and together with
\begin{equation} \label{multiplicative}
  Z(K)=Z(K,[t_{\text{\scriptsize min}},t_{\text{\scriptsize max}}])
	=Z(K,[t_{\text{\scriptsize min}},a]) Z(K,[a,b])
	 Z(K,[b,t_{\text{\scriptsize max}}])
\end{equation}
we see that $Z(K)$ is invariant under horizontal deformations of $K$
which `freeze' the time slices in which $K$ has a critical point.

\subsubsection{Moving critical points} \label{mov_crit_sec}
In this section we will show that (subject to some restrictions) $Z(K)$ is
also invariant under deformations of $K$ that do move critical points. The
idea is to narrow the parts near critical points to sharp needles
using horizontal deformations, and then show that very sharp needles
contribute almost nothing to $Z(K)$ and therefore can be moved around
freely. For example, here's how this trick allows us to move two critical
points across each other:
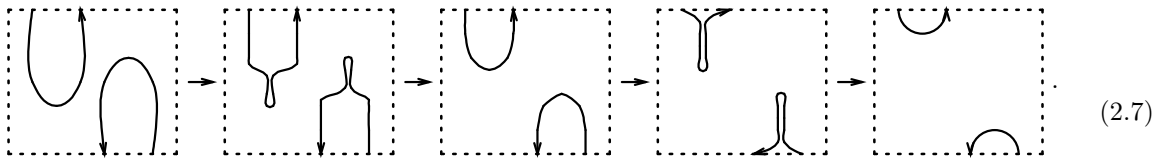
\begin{equation} \label{mov_crit}
  
        \setlength{\unitlength}{0.5\eepiclength}
        \begin{array}{c}  \hspace{-1.7mm}
                \raisebox{-8pt}{\begingroup\makeatletter\ifx\SetFigFont\undefined%
\gdef\SetFigFont#1#2#3#4#5{%
  \reset@font\fontsize{#1}{#2pt}%
  \fontfamily{#3}\fontseries{#4}\fontshape{#5}%
  \selectfont}%
\fi\endgroup%
\begin{picture}(12925,1840)(0,-10)
\thicklines
\path(11975.500,147.000)(12013.000,12.000)(12050.500,147.000)
\put(12313.000,12.000){\arc{600.000}{3.1416}{6.2832}}
\path(11750.500,1677.000)(11713.000,1812.000)(11675.500,1677.000)
\put(11413.000,1812.000){\arc{600.000}{6.2832}{9.4248}}
\path(4962,913)(5262,913)
\path(5127.000,875.500)(5262.000,913.000)(5127.000,950.500)
\path(2262,913)(2562,913)
\path(2427.000,875.500)(2562.000,913.000)(2427.000,950.500)
\path(10363,912)(10663,912)
\path(10528.000,874.500)(10663.000,912.000)(10528.000,949.500)
\path(7663,912)(7963,912)
\path(7828.000,874.500)(7963.000,912.000)(7828.000,949.500)
\dottedline{144}(12,1813)(2112,1813)(2112,13)
	(12,13)(12,1813)
\dottedline{144}(2712,1813)(4812,1813)(4812,13)
	(2712,13)(2712,1813)
\dottedline{144}(5412,1813)(7512,1813)(7512,13)
	(5412,13)(5412,1813)
\dottedline{144}(8113,1812)(10213,1812)(10213,12)
	(8113,12)(8113,1812)
\dottedline{144}(10813,1812)(12913,1812)(12913,12)
	(10813,12)(10813,1812)
\path(7212,13)	(7212.000,75.725)
	(7212.000,129.484)
	(7212.000,175.158)
	(7212.000,213.625)
	(7212.000,272.453)
	(7212.000,313.000)

\path(7212,313)	(7200.281,383.078)
	(7193.543,425.061)
	(7185.750,468.625)
	(7176.551,511.486)
	(7165.594,551.359)
	(7137.000,613.000)

\path(7137,613)	(7093.406,655.422)
	(7035.750,704.875)
	(6972.469,745.891)
	(6912.000,763.000)

\path(6912,763)	(6851.531,739.562)
	(6788.250,699.250)
	(6730.594,653.312)
	(6687.000,613.000)

\path(6687,613)	(6658.406,551.359)
	(6647.449,511.486)
	(6638.250,468.625)
	(6630.457,425.061)
	(6623.719,383.078)
	(6612.000,313.000)

\path(6612,313)	(6612.000,272.453)
	(6612.000,213.625)
	(6612.000,175.158)
	(6612.000,129.484)
	(6612.000,75.725)
	(6612.000,13.000)

\path(6574.500,148.000)(6612.000,13.000)(6649.500,148.000)
\path(5712,1813)	(5709.907,1752.503)
	(5708.411,1700.335)
	(5707.514,1655.616)
	(5707.215,1617.467)
	(5708.411,1557.366)
	(5712.000,1513.000)

\path(5712,1513)	(5720.781,1445.504)
	(5727.143,1404.441)
	(5735.096,1361.402)
	(5744.850,1318.564)
	(5756.613,1278.100)
	(5787.000,1213.000)

\path(5787,1213)	(5827.383,1164.870)
	(5882.501,1116.051)
	(5946.119,1078.207)
	(6012.000,1063.000)

\path(6012,1063)	(6077.881,1078.207)
	(6141.499,1116.051)
	(6196.617,1164.870)
	(6237.000,1213.000)

\path(6237,1213)	(6267.387,1278.100)
	(6279.150,1318.564)
	(6288.904,1361.402)
	(6296.857,1404.441)
	(6303.219,1445.504)
	(6312.000,1513.000)

\path(6312,1513)	(6315.589,1557.366)
	(6316.785,1617.467)
	(6316.486,1655.616)
	(6315.589,1700.335)
	(6314.093,1752.503)
	(6312.000,1813.000)

\path(6354.444,1679.472)(6312.000,1813.000)(6279.495,1676.711)
\path(4512,13)	(4512.000,75.725)
	(4512.000,129.484)
	(4512.000,175.158)
	(4512.000,213.625)
	(4512.000,272.453)
	(4512.000,313.000)

\path(4512,313)	(4505.672,383.078)
	(4505.408,425.061)
	(4506.375,468.625)
	(4508.045,511.486)
	(4509.891,551.359)
	(4512.000,613.000)

\path(4512,613)	(4506.375,644.875)
	(4512.000,688.000)

\path(4512,688)	(4462.078,726.203)
	(4393.875,757.375)
	(4328.484,791.359)
	(4287.000,838.000)

\path(4287,838)	(4274.344,883.703)
	(4275.750,944.875)
	(4282.781,1008.859)
	(4287.000,1063.000)

\path(4287,1063)	(4298.250,1132.375)
	(4299.656,1174.797)
	(4287.000,1213.000)

\path(4287,1213)	(4255.125,1224.250)
	(4212.000,1213.000)

\path(4212,1213)	(4199.344,1174.797)
	(4200.750,1132.375)
	(4212.000,1063.000)

\path(4212,1063)	(4218.328,1008.859)
	(4228.875,944.875)
	(4230.984,883.703)
	(4212.000,838.000)

\path(4212,838)	(4154.578,787.141)
	(4113.123,770.764)
	(4067.625,757.375)
	(4021.424,744.689)
	(3977.859,730.422)
	(3912.000,688.000)

\path(3912,688)	(3912.000,644.875)
	(3912.000,613.000)

\path(3912,613)	(3912.000,551.359)
	(3912.000,511.486)
	(3912.000,468.625)
	(3912.000,425.061)
	(3912.000,383.078)
	(3912.000,313.000)

\path(3912,313)	(3912.000,272.453)
	(3912.000,213.625)
	(3912.000,175.158)
	(3912.000,129.484)
	(3912.000,75.725)
	(3912.000,13.000)

\path(3874.500,148.000)(3912.000,13.000)(3949.500,148.000)
\path(3012,1813)	(3012.000,1752.329)
	(3012.000,1700.037)
	(3012.000,1655.243)
	(3012.000,1617.070)
	(3012.000,1557.068)
	(3012.000,1513.000)

\path(3012,1513)	(3012.000,1446.904)
	(3012.000,1406.065)
	(3012.000,1363.000)
	(3012.000,1319.935)
	(3012.000,1279.096)
	(3012.000,1213.000)

\path(3012,1213)	(3008.115,1176.355)
	(3012.000,1138.000)

\path(3012,1138)	(3057.669,1094.926)
	(3124.500,1063.000)
	(3191.331,1031.074)
	(3237.000,988.000)

\path(3237,988)	(3250.116,935.545)
	(3248.659,872.939)
	(3241.372,811.613)
	(3237.000,763.000)

\path(3237,763)	(3221.531,690.655)
	(3219.598,649.035)
	(3237.000,613.000)

\path(3237,613)	(3274.500,597.535)
	(3312.000,613.000)

\path(3312,613)	(3329.402,649.035)
	(3327.469,690.655)
	(3312.000,763.000)

\path(3312,763)	(3306.845,811.438)
	(3298.253,872.474)
	(3296.534,935.022)
	(3312.000,988.000)

\path(3312,988)	(3374.248,1036.311)
	(3416.530,1051.007)
	(3462.000,1063.000)
	(3507.470,1074.993)
	(3549.752,1089.689)
	(3612.000,1138.000)

\path(3612,1138)	(3616.586,1176.509)
	(3612.000,1213.000)

\path(3612,1213)	(3612.000,1279.096)
	(3612.000,1319.935)
	(3612.000,1363.000)
	(3612.000,1406.065)
	(3612.000,1446.904)
	(3612.000,1513.000)

\path(3612,1513)	(3612.000,1557.068)
	(3612.000,1617.070)
	(3612.000,1655.243)
	(3612.000,1700.037)
	(3612.000,1752.329)
	(3612.000,1813.000)

\path(3649.500,1678.000)(3612.000,1813.000)(3574.500,1678.000)
\path(1812,13)	(1818.327,59.696)
	(1824.245,104.781)
	(1829.755,148.298)
	(1834.857,190.287)
	(1839.551,230.789)
	(1843.837,269.846)
	(1851.184,343.789)
	(1856.898,412.446)
	(1860.980,476.145)
	(1863.429,535.216)
	(1864.245,589.990)
	(1863.429,640.795)
	(1860.980,687.961)
	(1856.898,731.818)
	(1851.184,772.696)
	(1834.857,846.830)
	(1812.000,913.000)

\path(1812,913)	(1793.304,956.238)
	(1769.854,1004.467)
	(1741.340,1054.372)
	(1707.454,1102.637)
	(1667.885,1145.948)
	(1622.326,1180.989)
	(1570.468,1204.445)
	(1512.000,1213.000)

\path(1512,1213)	(1453.532,1204.445)
	(1401.674,1180.989)
	(1356.115,1145.948)
	(1316.546,1102.637)
	(1282.660,1054.372)
	(1254.146,1004.467)
	(1230.696,956.238)
	(1212.000,913.000)

\path(1212,913)	(1189.143,846.830)
	(1172.816,772.696)
	(1167.102,731.818)
	(1163.020,687.961)
	(1160.571,640.795)
	(1159.755,589.990)
	(1160.571,535.216)
	(1163.020,476.145)
	(1167.102,412.446)
	(1172.816,343.789)
	(1180.163,269.846)
	(1184.449,230.789)
	(1189.143,190.287)
	(1194.245,148.298)
	(1199.755,104.781)
	(1205.673,59.696)
	(1212.000,13.000)

\path(1156.467,141.636)(1212.000,13.000)(1230.768,151.849)
\path(312,1813)	(305.673,1766.304)
	(299.755,1721.219)
	(294.245,1677.702)
	(289.143,1635.713)
	(284.449,1595.211)
	(280.163,1556.154)
	(272.816,1482.211)
	(267.102,1413.554)
	(263.020,1349.855)
	(260.571,1290.784)
	(259.755,1236.010)
	(260.571,1185.205)
	(263.020,1138.039)
	(267.102,1094.182)
	(272.816,1053.304)
	(289.143,979.170)
	(312.000,913.000)

\path(312,913)	(330.693,869.762)
	(354.142,821.533)
	(382.655,771.628)
	(416.543,723.362)
	(456.112,680.052)
	(501.672,645.011)
	(553.532,621.555)
	(612.000,613.000)

\path(612,613)	(670.468,621.555)
	(722.326,645.011)
	(767.885,680.052)
	(807.454,723.362)
	(841.340,771.628)
	(869.854,821.533)
	(893.304,869.762)
	(912.000,913.000)

\path(912,913)	(934.857,979.170)
	(951.184,1053.304)
	(956.898,1094.182)
	(960.980,1138.039)
	(963.429,1185.205)
	(964.245,1236.010)
	(963.429,1290.784)
	(960.980,1349.855)
	(956.898,1413.554)
	(951.184,1482.211)
	(943.837,1556.154)
	(939.551,1595.211)
	(934.857,1635.713)
	(929.755,1677.702)
	(924.245,1721.219)
	(918.327,1766.304)
	(912.000,1813.000)

\path(967.533,1684.364)(912.000,1813.000)(893.232,1674.151)
\path(9913,12)	(9864.083,33.558)
	(9822.425,53.435)
	(9758.249,89.902)
	(9688.000,162.000)

\path(9688,162)	(9679.255,196.972)
	(9680.226,238.710)
	(9688.000,312.000)

\path(9688,312)	(9688.000,345.048)
	(9688.000,387.000)
	(9688.000,428.952)
	(9688.000,462.000)

\path(9688,462)	(9688.000,495.048)
	(9688.000,537.000)
	(9688.000,578.952)
	(9688.000,612.000)

\path(9688,612)	(9703.469,684.345)
	(9705.402,725.965)
	(9688.000,762.000)

\path(9688,762)	(9650.500,777.465)
	(9613.000,762.000)

\path(9613,762)	(9595.598,725.965)
	(9597.531,684.345)
	(9613.000,612.000)

\path(9613,612)	(9613.000,578.952)
	(9613.000,537.000)
	(9613.000,495.048)
	(9613.000,462.000)

\path(9613,462)	(9613.000,428.952)
	(9613.000,387.000)
	(9613.000,345.048)
	(9613.000,312.000)

\path(9613,312)	(9622.165,239.018)
	(9623.311,197.318)
	(9613.000,162.000)

\path(9613,162)	(9574.703,117.072)
	(9516.625,79.440)
	(9477.971,62.261)
	(9431.734,45.588)
	(9377.037,28.981)
	(9313.000,12.000)

\path(9434.814,81.229)(9313.000,12.000)(9453.067,8.484)
\path(8413,1812)	(8461.917,1790.436)
	(8503.575,1770.554)
	(8567.751,1734.082)
	(8638.000,1662.000)

\path(8638,1662)	(8646.745,1627.028)
	(8645.774,1585.290)
	(8638.000,1512.000)

\path(8638,1512)	(8638.000,1478.952)
	(8638.000,1437.000)
	(8638.000,1395.048)
	(8638.000,1362.000)

\path(8638,1362)	(8638.000,1328.952)
	(8638.000,1287.000)
	(8638.000,1245.048)
	(8638.000,1212.000)

\path(8638,1212)	(8622.531,1139.655)
	(8620.598,1098.035)
	(8638.000,1062.000)

\path(8638,1062)	(8675.500,1046.535)
	(8713.000,1062.000)

\path(8713,1062)	(8730.402,1098.035)
	(8728.469,1139.655)
	(8713.000,1212.000)

\path(8713,1212)	(8713.000,1245.048)
	(8713.000,1287.000)
	(8713.000,1328.952)
	(8713.000,1362.000)

\path(8713,1362)	(8713.000,1395.048)
	(8713.000,1437.000)
	(8713.000,1478.952)
	(8713.000,1512.000)

\path(8713,1512)	(8703.831,1584.983)
	(8702.685,1626.682)
	(8713.000,1662.000)

\path(8713,1662)	(8751.297,1706.928)
	(8809.375,1744.560)
	(8848.029,1761.739)
	(8894.266,1778.412)
	(8948.963,1795.019)
	(9013.000,1812.000)

\path(8891.186,1742.771)(9013.000,1812.000)(8872.933,1815.516)
\end{picture} }
                \hspace{-1.9mm}
        \end{array}
.
\end{equation}

\begin{lemma} \label{needles}
If the two knots $K_{1,2}$ both contain a sharp needle of
width $\epsilon$, and are the identical except possibly for the length and the
directions of their respective needles, then
\[ ||Z_m(K_1)-Z_m(K_2)||\sim\epsilon \]
where $Z_m$ is the degree $m$ piece of $Z$ and $||\cdot||$ is some fixed
norm on ${\cal A}^r_m$.
\end{lemma}

\proof{}
Clearly, the difference between $Z_m(K_1)$ and $Z_m(K_2)$ will come only
from terms in \eqref{Kints} in which one of the $z_i$'s (or $z'_i$'s)
is on the needle. So
let us show that if a knot $K$ contains a needle $N$ of width $\epsilon$,
then such terms in $Z_m(K)$ are at most proportional to $\epsilon$. Without
loss of generality we can assume that the needle $N$ points upward. If the
highest pair $(z_i,z'_i)$ that touches $N$ connects the two sides of $N$,
the corresponding diagram is $0$ in ${\bar{\cal A}}^r$ and there is
nothing to worry about. If there is no pair $(z_j,z'_j)$ that connects
the two sides of $N$ then again life is simple: in that case there are
no singularities in \eqref{Kints} so nothing big prevents
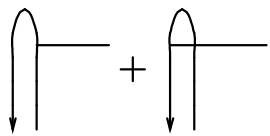
\[ 
        \setlength{\unitlength}{0.5\eepiclength}
        \begin{array}{c}  \hspace{-1.7mm}
                \raisebox{-8pt}{\begingroup\makeatletter\ifx\SetFigFont\undefined%
\gdef\SetFigFont#1#2#3#4#5{%
  \reset@font\fontsize{#1}{#2pt}%
  \fontfamily{#3}\fontseries{#4}\fontshape{#5}%
  \selectfont}%
\fi\endgroup%
\begin{picture}(3193,1531)(0,-10)
\thicklines
\path(1516,908)(1516,608)
\path(1366,758)(1666,758)
\path(316,1058)(1216,1058)
\path(1981,1058)(3181,1058)
\path(316,8)	(316.000,68.671)
	(316.000,120.963)
	(316.000,165.757)
	(316.000,203.930)
	(316.000,263.932)
	(316.000,308.000)

\path(316,308)	(316.000,374.096)
	(316.000,414.935)
	(316.000,458.000)
	(316.000,501.065)
	(316.000,541.904)
	(316.000,608.000)

\path(316,608)	(316.000,674.096)
	(316.000,714.935)
	(316.000,758.000)
	(316.000,801.065)
	(316.000,841.904)
	(316.000,908.000)

\path(316,908)	(318.350,973.747)
	(320.406,1014.281)
	(322.266,1057.070)
	(323.343,1099.975)
	(323.050,1140.858)
	(316.000,1208.000)

\path(316,1208)	(307.572,1243.621)
	(296.279,1288.275)
	(282.079,1337.713)
	(264.929,1387.689)
	(244.788,1433.953)
	(221.614,1472.258)
	(166.000,1508.000)

\path(166,1508)	(110.386,1472.258)
	(87.212,1433.953)
	(67.071,1387.689)
	(49.921,1337.713)
	(35.721,1288.275)
	(24.428,1243.621)
	(16.000,1208.000)

\path(16,1208)	(8.950,1140.858)
	(8.657,1099.975)
	(9.734,1057.070)
	(11.594,1014.281)
	(13.650,973.747)
	(16.000,908.000)

\path(16,908)	(16.000,841.904)
	(16.000,801.065)
	(16.000,758.000)
	(16.000,714.935)
	(16.000,674.096)
	(16.000,608.000)

\path(16,608)	(16.000,541.904)
	(16.000,501.065)
	(16.000,458.000)
	(16.000,414.935)
	(16.000,374.096)
	(16.000,308.000)

\path(16,308)	(16.000,263.932)
	(16.000,203.930)
	(16.000,165.757)
	(16.000,120.963)
	(16.000,68.671)
	(16.000,8.000)

\path(-21.500,143.000)(16.000,8.000)(53.500,143.000)
\path(2281,8)	(2281.000,68.671)
	(2281.000,120.963)
	(2281.000,165.757)
	(2281.000,203.930)
	(2281.000,263.932)
	(2281.000,308.000)

\path(2281,308)	(2281.000,374.096)
	(2281.000,414.935)
	(2281.000,458.000)
	(2281.000,501.065)
	(2281.000,541.904)
	(2281.000,608.000)

\path(2281,608)	(2281.000,674.096)
	(2281.000,714.935)
	(2281.000,758.000)
	(2281.000,801.065)
	(2281.000,841.904)
	(2281.000,908.000)

\path(2281,908)	(2283.350,973.747)
	(2285.406,1014.281)
	(2287.266,1057.070)
	(2288.343,1099.975)
	(2288.050,1140.858)
	(2281.000,1208.000)

\path(2281,1208)	(2272.572,1243.621)
	(2261.279,1288.275)
	(2247.079,1337.713)
	(2229.929,1387.689)
	(2209.788,1433.953)
	(2186.614,1472.258)
	(2131.000,1508.000)

\path(2131,1508)	(2075.386,1472.258)
	(2052.212,1433.953)
	(2032.071,1387.689)
	(2014.921,1337.713)
	(2000.721,1288.275)
	(1989.428,1243.621)
	(1981.000,1208.000)

\path(1981,1208)	(1973.950,1140.858)
	(1973.657,1099.975)
	(1974.734,1057.070)
	(1976.594,1014.281)
	(1978.650,973.747)
	(1981.000,908.000)

\path(1981,908)	(1981.000,841.904)
	(1981.000,801.065)
	(1981.000,758.000)
	(1981.000,714.935)
	(1981.000,674.096)
	(1981.000,608.000)

\path(1981,608)	(1981.000,541.904)
	(1981.000,501.065)
	(1981.000,458.000)
	(1981.000,414.935)
	(1981.000,374.096)
	(1981.000,308.000)

\path(1981,308)	(1981.000,263.932)
	(1981.000,203.930)
	(1981.000,165.757)
	(1981.000,120.963)
	(1981.000,68.671)
	(1981.000,8.000)

\path(1943.500,143.000)(1981.000,8.000)(2018.500,143.000)
\end{picture} }
                \hspace{-1.9mm}
        \end{array}
 \]
from being small. (Notice that these two terms appear in $Z(K)$ with opposite
signs due to the factor $(-1)^{\#P_{\downarrow}}$ but otherwise they differ
only by something proportional to $\epsilon$). If $(z_j,z'_j)$ is a pair
that does connect the two sides of $N$, it has to do so in the top (round)
part of $N$ --- otherwise $dz_j-dz'_j=0$.

So the only terms that cause some worry are those that have some $k>1$
pairs $(z_{j_1},z'_{j_1}),\ldots,(z_{j_k},z'_{j_k})$ on the top part of
$N$, with $(z_{j_k},z'_{j_k})$ being the highest of these pairs and
$(z_{j_1},z'_{j_1})$ the lowest. We might as well assume that there are no
pairs other than $(z_i,z'_i)$ that touch $N$ only once --- such pairs just
shorten the domain of integration in \eqref{Kints} without adding any
singularity in the denominator. So what we have looks like:
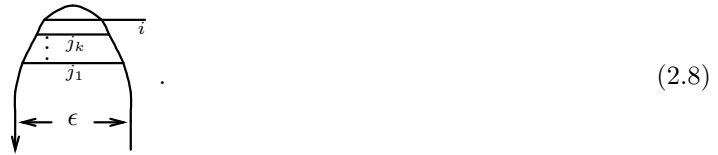
\begin{equation} \label{worry} 
        \setlength{\unitlength}{0.6\eepiclength}
        \begin{array}{c}  \hspace{-1.7mm}
                \raisebox{-8pt}{\begingroup\makeatletter\ifx\SetFigFont\undefined%
\gdef\SetFigFont#1#2#3#4#5{%
  \reset@font\fontsize{#1}{#2pt}%
  \fontfamily{#3}\fontseries{#4}\fontshape{#5}%
  \selectfont}%
\fi\endgroup%
\begin{picture}(1425,1531)(0,-10)
\thicklines
\path(837,293)(1137,293)
\path(1002.000,255.500)(1137.000,293.000)(1002.000,330.500)
\path(387,293)(87,293)
\path(222.000,330.500)(87.000,293.000)(222.000,255.500)
\path(312,1358)(1362,1358)
\path(237,1208)(987,1208)
\path(1137,908)(87,908)
\path(1212,8)	(1212.000,68.671)
	(1212.000,120.963)
	(1212.000,165.757)
	(1212.000,203.930)
	(1212.000,263.932)
	(1212.000,308.000)

\path(1212,308)	(1213.197,373.994)
	(1214.244,414.745)
	(1215.191,457.730)
	(1215.740,500.749)
	(1215.590,541.601)
	(1212.000,608.000)

\path(1212,608)	(1200.072,675.484)
	(1190.956,716.535)
	(1180.511,759.564)
	(1169.315,802.397)
	(1157.945,842.862)
	(1137.000,908.000)

\path(1137,908)	(1107.019,975.931)
	(1086.727,1016.846)
	(1064.726,1059.631)
	(1042.384,1102.212)
	(1021.070,1142.516)
	(987.000,1208.000)

\path(987,1208)	(972.341,1241.584)
	(954.476,1284.429)
	(912.000,1358.000)

\path(912,1358)	(852.992,1407.054)
	(815.675,1432.145)
	(775.185,1455.770)
	(733.046,1476.549)
	(690.786,1493.101)
	(649.928,1504.045)
	(612.000,1508.000)

\path(612,1508)	(574.072,1504.045)
	(533.214,1493.101)
	(490.954,1476.549)
	(448.815,1455.770)
	(408.325,1432.145)
	(371.008,1407.054)
	(312.000,1358.000)

\path(312,1358)	(269.520,1284.429)
	(251.655,1241.584)
	(237.000,1208.000)

\path(237,1208)	(202.930,1142.516)
	(181.616,1102.212)
	(159.274,1059.631)
	(137.273,1016.846)
	(116.981,975.931)
	(87.000,908.000)

\path(87,908)	(66.051,842.862)
	(54.681,802.397)
	(43.485,759.564)
	(33.041,716.535)
	(23.927,675.484)
	(12.000,608.000)

\path(12,608)	(8.410,541.601)
	(8.260,500.749)
	(8.809,457.730)
	(9.756,414.745)
	(10.803,373.994)
	(12.000,308.000)

\path(12,308)	(12.000,263.932)
	(12.000,203.930)
	(12.000,165.757)
	(12.000,120.963)
	(12.000,68.671)
	(12.000,8.000)

\path(-25.500,143.000)(12.000,8.000)(49.500,143.000)
\put(552,263){\makebox(0,0)[lb]{\smash{{\mathmode{\epsilon}}}}}
\put(537,758){\makebox(0,0)[lb]{\smash{{\mathmode{\scriptstyle j_1}}}}}
\put(537,1058){\makebox(0,0)[lb]{\smash{{\mathmode{\scriptstyle j_k}}}}}
\put(1287,1208){\makebox(0,0)[lb]{\smash{{\mathmode{\scriptstyle i}}}}}
\put(304,938){\makebox(0,0)[lb]{\smash{{\mathmode{\scriptstyle \vdots}}}}}
\end{picture} }
                \hspace{-1.9mm}
        \end{array}
. \end{equation}
Writing $\delta_\alpha=|z_{j_\alpha}-z'_{j_\alpha}|$,
we see that the integral corresponding to
\eqref{worry} is bounded by a constant times
\begin{equation}
  \int_0^\epsilon \frac{d\delta_1}{\delta_1}
  \int_0^{\delta_1}\frac{d\delta_2}{\delta_2}
  \cdots\int_0^{\delta_{k-1}}\frac{d\delta_k}{\delta_k}
  \int_{z_{j_k}}^{z'_{j_k}}\frac{dz_i-dz'_i}{z_i-z'_i}
  \sim \epsilon.
% \tag*{$\Box$}
\prtag
\end{equation}

Unfortunately, there is one type of deformation that \eqref{mov_crit} and
lemma \ref{needles} cannot handle --- the total number of critical points
in $K$ cannot be changed:
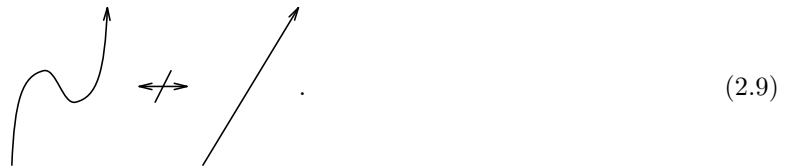
\begin{equation} \label{cubic} \psdraw{cubic}{1.5in}. \end{equation}
Even if the hump on the left figure is deformed into a needle and then this
needle is removed, a (smaller) hump still remains.

\subsubsection{The correction} Let the symbol $\infty$ stand for the embedding
$
        \setlength{\unitlength}{0.5\eepiclength}
        \begin{array}{c}  \hspace{-1.7mm}
                \raisebox{-8pt}{\begingroup\makeatletter\ifx\SetFigFont\undefined%
\gdef\SetFigFont#1#2#3#4#5{%
  \reset@font\fontsize{#1}{#2pt}%
  \fontfamily{#3}\fontseries{#4}\fontshape{#5}%
  \selectfont}%
\fi\endgroup%
\begin{picture}(1174,704)(0,-10)
\thicklines
\path(532,437)	(492.231,485.288)
	(462.360,520.264)
	(422.000,562.000)

\path(422,562)	(354.139,634.775)
	(308.795,666.974)
	(256.000,675.000)

\path(256,675)	(200.806,659.090)
	(151.799,632.032)
	(109.341,595.762)
	(73.796,552.216)
	(45.527,503.331)
	(24.895,451.042)
	(12.266,397.286)
	(8.000,344.000)

\path(8,344)	(12.247,290.594)
	(24.718,236.405)
	(45.124,183.562)
	(73.174,134.191)
	(108.577,90.422)
	(151.042,54.381)
	(200.280,28.198)
	(256.000,14.000)

\path(256,14)	(324.788,28.168)
	(381.190,75.148)
	(426.747,132.803)
	(463.000,179.000)

\path(463,179)	(520.566,250.108)
	(553.229,295.849)
	(587.000,344.500)
	(620.771,393.151)
	(653.434,438.892)
	(711.000,510.000)

\path(711,510)	(747.264,556.206)
	(792.840,613.875)
	(849.246,660.857)
	(918.000,675.000)

\path(918,675)	(973.777,660.730)
	(1023.057,634.441)
	(1065.551,598.271)
	(1100.965,554.357)
	(1129.009,504.836)
	(1149.390,451.845)
	(1161.818,397.521)
	(1166.000,344.000)

\path(1166,344)	(1161.753,290.594)
	(1149.282,236.405)
	(1128.876,183.562)
	(1100.826,134.191)
	(1065.423,90.422)
	(1022.958,54.381)
	(973.720,28.198)
	(918.000,14.000)

\path(918,14)	(849.068,27.900)
	(792.427,74.431)
	(746.823,131.997)
	(711.000,179.000)

\path(711,179)	(690.390,206.839)
	(657.000,262.000)

\path(757.397,164.267)(657.000,262.000)(692.607,126.488)
\end{picture} }
                \hspace{-1.9mm}
        \end{array}
$. Notice that
\begin{equation} \label{Zinfty}
	Z(\infty)=\psdraw{zero}{0.2in}+(\text{higher order terms})
\end{equation}
and so using power series $Z(\infty)$ can be inverted and the following
definition makes sense:

\begin{defi} \label{correction}
Let $K$ be a knot embedded in ${\bold C}\times{\bold R}$
with $c$ critical points. Notice that $c$ is always even and
set\footnote{The non-invariance of $Z(K)$ under the move \eqref{cubic} was
first noticed by R.~Bott and the first author. The correction $\tilde{Z}(K)$ is
due to Kontsevich \cite{Kontsevich:Vassiliev}.}
\[ \tilde{Z}(K)=\frac{Z(K)}{(Z(\infty))^{\frac{c}{2}}}. \]
\end{defi}

\begin{theo} \label{tildeZinv}
$\tilde{Z}(K)$ is invariant under arbitrary deformations of the knot $K$.
\end{theo}

\proof{} Clearly, $\tilde{Z}(K)$ is invariant under deformations that do
not change the number of critical points of $K$, and the only thing that
remains to be checked is its invariance under the move \eqref{cubic}. So
let $K_c$ and $K_s$ be two knots that are identical other then that in some
place $K_c$ has the figure in the left side of \eqref{cubic} while in the
same place $K_s$ has the figure on the right side of \eqref{cubic}. We need
to show that in ${\bar{\cal A}}^r$,
\[ Z(K_c)=Z(\infty)Z(K_s). \]
Using deformations as in section \ref{mov_crit_sec} we can move the `humps'
of $K_c$ to be very far from the rest of the knot, and shrink them to be
very small. This done, we can ignore contributions to $Z(K_c)$ coming from
pairings in which any of the pairs connect the humps to the rest of the
knot. Hence $Z(K_c)$ factors to a part which is the same as in $Z(K_s)$
times contributions that come from pairings that pair the `humpy' part of
$K_c$ to itself. But as the following figure shows, for the same reasons
as in section \ref{mov_crit_sec}, these contributions are
precisely $Z(\infty)$:
%{\def\@eqnnum{\reset@font \rm\qed}\let\stepcounter\@gobble
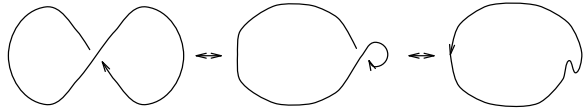
\begin{equation}
  \psdraw{infty_cubic}{3in}
% \tag*{$\Box$}
\prtag
\end{equation}
%}

\begin{exercise} Show that $\tilde{Z}(K)$ is in fact real, even though
complex numbers do appear in~\eqref{ZDef}.
\end{exercise}

\begin{hint} Use the fact that the transformation $t\rightarrow-t$,
$z\rightarrow\bar{z}$ maps a knot to an equivalent knot, while mapping
$\Omega_{n,n}$ to minus its conjugate.
\end{hint}

\begin{remark} Le and Murakami~\cite{LeMurakami:Universal}, building on
work of Drinfel'd (\cite{Drinfeld:QuasiHopf} and \cite{Drinfeld:GalQQ}),
proved that $\tilde{Z}(K)$ has rational coefficients.
\end{remark}

\subsection{Universallity of the Kontsevich integral.}

It is enough to show that if $D\in{\cal D}^0_m$ is a chord diagram of
degree $m$ underlying some $m$-singular knot $K$, then (for the natural
extension of $\tilde{Z}$ to knots with double points):
\[ \tilde{Z}(K)=D+(\text{terms of degree $>m$}). \]
In view
of \eqref{Zinfty}, it is enough to prove the same for $Z$ rather than for
$\tilde{Z}$. If two knots $K^o$ and $K^u$ are identical except that two of
their strands form an overcrossing in $K^o$ and an undercrossing in $K^u$,
it is clear that the only contributions to $Z(K^o)-Z(K^u)$ come from
pairings in which these two strands are paired. $Z(K_D)$ is a signed
sum of $Z$ evaluated on $2^m$ knots, and this sum can be partitioned in
pairs like the above $K^{o,u}$ around $m$ different crossings --- and thus
contributions to $Z(K_D)$ come only from pairings that pair the strands
near any of the $m$ double points of $K_D$. This implies that the lowest
degree contribution to $Z(K_D)$ is at least of degree $m$. In degree $m$
the pairing $P$ is determined by the above restriction. It is easy to see
that in that case $D_P=D$, and therefore the piece of degree precisely
$m$ in $Z(K_D)$ is proportional to $D$. It remains to determine
the constant of proportionality. This is a simple computation --- in
degree $1$, the difference between $Z(K^o)$ and $Z(K^u)$ comes from the
difference between integrating
\[ \frac{dz-dz'}{z-z'} \]
along a contour in which $z$ passes near but above $z'$ and along a
contour in which $z$ passes near but under $z'$. By Cauchy's theorem this
is $2\pi i$. Repeating this $m$ times for each of the $m$ double points
of $K_D$, we get $(2\pi i)^m$ and this exactly cancels the $(2\pi i)^{-m}$
in \eqref{Kints}.  \qed

\subsection{Why are we not happy?}

\begin{enumerate}
\item Why did we have to choose time axis?
\item Why did analysis (estimates for needles \dots)
come in?
\item Why did the real numbers come all together? The theorem can
be formulated over an arbitrary Abelian group. Is it
true in that generality?
\end{enumerate}

\lectitle{Physics (sketch)}

\begin{remark} This is the oldest approach, about 5-6 years old. Here we
follow the presentation in~\cite{Bar-Natan:PerturbativeCS}
and~\cite{Bar-Natan:Thesis}.
\end{remark}

\subsection{Invariants from path integrals.}

\noindent\underline{\bf Reminder:}\quad We are looking for a knot invariant
\[
\tZ\,:\,\bigl\{\text{knots}\bigr\}\,\rightarrowto{1.2cm}\,
\bar{\cA}^r\,=\,
\Span\left\{
\diag{5mm}{2}{2}{
\piccirclevecarc{1 1}{1}{30 390}
\pictranslate{1 1}{
\picveclength{0.1}
\picvecwidth{0.05}
\piccirclevecarc{-0.2 0}{0.2}{0 250}
%\piclinedash{0.2 0.1}{0.25}
\picline{1 120 polar}{-0.2 0}
\picline{1 240 polar}{-0.2 0}
\picline{1 0 polar}{-0.2 0}
\picline{1 70 polar}{1 -70 polar}
}}\right\}\,\Bigg/\,\mbox{\footnotesize$\displaystyle
\begin{array}{l}FI\\AS\\IHX\\STU\end{array}$}
\]
such that if $K$ is singular, $\tZ(K)=D_K\,+\,(\text{higher degrees})$,
where $D_K$ is the chord diagram underlying $K$.

\myitemm{\underline{Idea}\quad} Geometric invariants are cheaper that
topological ones. So introduce a geometrical structure
$A$, get an invariant and average out over all possible
choices of $A$.

\def\cZ{{\cal Z}}
\def\cD{{\cal D}}
\def\bR{{\bold R}}

\begin{example} Define
\[
\cZ_k(\bR^3,K)\,=\,\int\limits_{\Omega^1(\bR^3,\fg)}\,\cD A\,\tr_R\, hol_K(A)\,
\cdot\,I_k\left(\frac{1}{4\pi}\int_{\bR^3}\,\tr(A\land dA+\frac 34
A\land A\land A)\right)
\]
where $\fg$ is a Lie algebra, $A$ is a $\fg$-connection on $\bR^3$,
$\Omega^1(\bR^3,\fg)$ is the space of all such connections, $R$ is a
representation of $\fg$, $\pi$ is the ratio of the circumference and the
diameter of a circle, and $I_k(z)$ is the $k$'th modified Bessel function
of the first kind.
\end{example}

This is of course silly. Most of us don't even remember what a Bessel
function is, and certainly not how to integrate Bessel functions on spaces
of high dimension. None of us knows how to integrate things like that on
infinite dimensional spaces. When we toss the question to our physicist
friends, we find that they never really meant to say that integration
on infinite dimensional spaces is possible. Only that
\begin{itemize}
\item[*] for {\em very} special types of integrands there is a very
  complicated {\em formal} integration technique,
\item[*] which is very delicate and plagued with several layers of
  unexpected difficulties.
\end{itemize}

\def\myenumlabel#1#2{%
\ea\def\csname label#1\endcsname{\xdef\@currentlabel
{\csname the#1\endcsname}#2}}

The integration technique is
\begin{enumerate}
\myenumlabel{enumi}{\underline{Step \theenumi}}
\item \label{stepone}
Find something you can do in $\bR^n$ for all $n$
with a closed-form answer which depends lightly on $n$.
\item \label{steptwo}
Roughly, ``substitute $n=\infty$'' and
{\em hope} that everything still makes sense.
\end{enumerate}

Step \ref{stepone} basically restricts us to deal
with integrals of the form
\begin{equation} \label{Doable}
\int\,(\text{polynomial})\,e^{\kappa\,\left(\text{quadratic}\,+
\,\vctext{@{}c@{}}{higher order\\perturbations}{\scriptsize}\,\right)}
\end{equation}
(step \ref{steptwo} will put even further restrictions). So we're left with
\begin{equation} \label{CS34}
\cZ_k(\bR^3,K)\,=\,\int\,\cD A\,\tr_R\, hol_K(A)\,
\cdot\,\exp\left(\frac{ik}{4\pi}\int_{\bR^3}\,\tr(A\land dA+\frac 34
A\land A\land A)\right)\,,
\end{equation}
which is of the required form because
\[
hol_k(A)\,=\,\sum_{m=0}^\infty\,\int\limits_{0\le t_1\le
\dots\le t_m\le 1}\kern-1em (K^*A)(t_m)\cdot\,\dots\,\cdot
(K^*A)(t_1)
\]
is a polynomial (oh well, power series) in $A$.

\subsection{A finite dimensional analogue.} Let us start by showing how
integrals like~\eqref{Doable} are computed in $\bR^n$. By rescaling $\vec
x$ and Taylor expanding,

\begin{eqnarray*}
  \int_{\bR^n}d\vec xe^{it\left(\frac 12\lambda_{ij}x^ix^j
  +\lambda_{ijk}x^ix^jx^k\right)}
& \propto &
  \int_{\bR^n}d\vec xe^{\frac i2\lambda_{ij}x^ix^j+\frac i{\sqrt t}
  \lambda_{ijk}x^ix^jx^k}\\
&=&
  \int_{\bR^n}d\vec xe^{\frac i2\lambda_{ij}x^ix^j}\,\sum_{m=0}^\infty
  \frac {(-1)^m}{(2m)!\,t^m}\bigl(\lambda_{ijk}x^ix^jx^k\bigr)^{2m}.
\end{eqnarray*}
Picking up just one term, for simplicity:
\begin{eqnarray}
\int\limits_{\bR^n}d\vec xe^{\frac i2\lambda_{ij}x^ix^j}\,
\frac {(-1)^m}{(2m)!\,t^m}\bigl(\lambda_{ijk}x^ix^jx^k\bigr)^{2m}\quad&=&
\left(\lambda_{ijk}\,\frac{-i\partial}{\partial J_i}
\frac{-i\partial}{\partial J_j}
\frac{-i\partial}{\partial J_k}
\right)^{2m}\,\left.\int
d\vec xe^{\frac i2\lambda_{ij}x^ix^j+iJ_ix^i}\right|_{J=0} \nonumber \\
&\propto&\left.
\left(\lambda_{ijk}\,\frac{-i\partial}{\partial J_i}
\frac{-i\partial}{\partial J_j}
\frac{-i\partial}{\partial J_k}
\right)^{2m}e^{-\frac i2\lambda^{\alpha\beta}J_\alpha J_\beta}\right|_{J=0}\,,
\label{NeedToUnderstand} \end{eqnarray}
where $\lambda^{\alpha\beta}$ is the inverse of $\lambda_{ij}$.

Now comes a {\em combinatorial} challenge (notice that there are no
integrals left). We need to understand the multiple differentiations
in~\eqref{NeedToUnderstand}. When all the dust settles, this becomes:

\def\cE{{\cal E}}

\def\sum@feyn#1{\sum_{\mbox{\scriptsize$\begin{array}{c}
\scriptstyle \text{#1}\\
\scriptstyle \text{diagrams} D
\end{array}$}}}
\[
\sum@feyn{Feynman}
\left(\sum_{\text{labels}}\,\cE(D)\right),
\]
where $\cE(D)$ is defined as in the following example

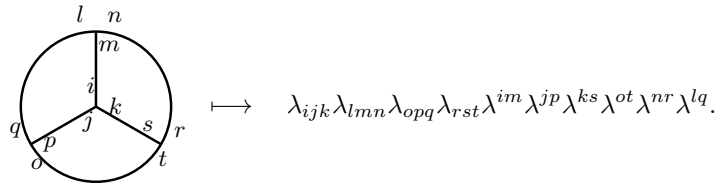
\def\@temp#1#2#3{%
\pictext{\small$#1$}{0.2 #2 polar}{#3}}

\[
\diag{10mm}{2}{2}
{
\pictranslate{1 1}{
\piccircle{0 0}{1}{}
\picline{0 0}{1 90 polar}
\picline{0 0}{1 210 polar}
\picline{0 0}{1 330 polar}
\pictranslate{1 90 polar}{
\@temp n{40}{0 0}
\@temp l{140}{-1 0}
\@temp m{280}{0 -0.5}
}
\pictranslate{1 210 polar}{
\@temp p{30}{0 -1}
\@temp q{130}{-1 0}
\@temp o{270}{0 -1}
}
\pictranslate{1 330 polar}{
\@temp r{30}{0 0}
\@temp s{120}{-1 0}
\@temp t{260}{0 -0.5}
}
\@temp i{110}{-0.5 0}
\@temp j{240}{-0.5 -0.5}
\@temp k{330}{0 0}
}
}
\quad\longmapsto\quad
\lambda_{ijk}\lambda_{lmn}\lambda_{opq}\lambda_{rst}
\lambda^{im}\lambda^{jp}\lambda^{ks}\lambda^{ot}\lambda^{nr}
\lambda^{lq}.
\]

\noindent\underline{The ellipticity problem:} We see that in
computing~\eqref{Doable}, we need to invert the quadratic piece. So to
compute~\eqref{CS34} we need to invert $A\land dA$.
But it is not invertible because $\{dC\}$ is in the radical of
this quadratic form.

\noindent\underline{Back to finite dimensions:}
If $L=\frac 12\lambda_{ij}x^ix^j\,+\,\lambda_{ijk}x^ix^jx^k$
is invariant under some $l$-dimensional group action you can
integrate over a section:

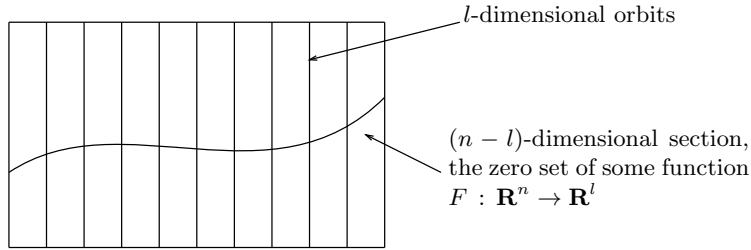
\[
\diag{5mm}{20}{6}{
\picmultigraphics{11}{1 0}{\picline{0 0}{0 6}}
\picline{0 6}{10 6}
\picline{0 0}{10 0}
\piccurve{0 2}{3 4}{7 1}{10 4}
\picvecline{12 6}{8 5}
\picvecline{11.5 2}{9.5 3}
\pictext{\small$l$-dimensional orbits}{12.1 6}{0 0}
\pictext{\small\parbox{4cm}{$(n-l)$-dimensional section, the zero set of some function
$F\,:\,\bR^n\to\bR^l$}}{11.7 2}{0 0}
}
\]

On the section, the quadratic form is non-degenerate. But the section
may not be a linear space!

\noindent{\bf Solution:} (The ``Faddeev-Popov procedure''.) Integrate against a
$\delta$-function concentrated on the section, and include a Jacobian
which measures both the volume of the orbit and the ``angle'' with which
the orbit meets the section. That is, compute

\[
  \int\,d\vec xe^{it\left(\frac 12\lambda_{ij}x^ix^j\,+\,
  \lambda_{ijk}x^ix^jx^k\right)}\delta^l(F(\vec x))\det
  \left(\frac{\partial F^a}{\partial \fg_b}\right)(x)
\]

By Fourier analysis
\[
\delta^l(F(\vec x))\,=\int d\vec\Phi e^{iF^a(\vec x)\Phi_a}.
\]
By cheating (or by introducing ``anti-commuting variables'')
\[
\det\left(\frac{\partial F^a}{\partial \fg_b}\right)\,=\,
\int\,d\vec c\,d\vec c\,e^{i\bar c_a
\frac{\partial F^a}{\partial \fg_b}c^b
}
\]
We end up back again with a Gaussian, this time non-degenerate:
\[
\int \,d\vec x\,d\vec\Phi \,d\vec c\,d\vec c
\,\exp\left(t\left(\frac 12\lambda_{ij}x^ix^j\,+\,
\lambda_{ijk}x^ix^jx^k\right)+F^a(\vec x)\Phi_a+\bar c_a
\frac{\partial F^a}{\partial \fg_b}c^b\right)\,.
\]

%\begin{eqnarray*}

\subsection{Chern-Simons perturbation theory.} Setting $\frac 34=\frac 23$,
our Lagrangian becomes the Chern-Simons-functional
\[
\int\,\tr\bigl(A\land dA+\frac 23A\land A\land A\bigr)\,,
\]
which is gauge-invariant. Applying Faddeev-Popov (with some $F$)
and then crunching Feynman diagrams, we get a messier Lagrangian,
whose perturbation theory has the following general form:

\[
\cZ_k(K)\sim\sum_{m=0}^\infty\frac 1{k^m}\hspace{-2mm}\sum@feyn{degree $m$}
\hspace{-2mm}W(D)\hspace{-2mm}
\rx{-1em}\sumint_{\mbox{\scriptsize \begin{tabular}{@{}c@{}}%
labels in $\{1,2,3\},$\\$\{1,2,\dots,\dim\fg\},$\\$\bR^3$ and $S^1$
\end{tabular}}}
\rx{-1em}\cE(D) \qquad
\raisebox{-3mm}{$\left(\mbox{\scriptsize
\parbox{4.2cm}{%
\noindent{\tiny $\sumint$}
  is a symbol that should have long been introduced into
  mathematics. It means ``sum over discrete variables and integrate over
  continuous ones''.
}}\right)$}
\]
Here $W(D)$ is the Lie-algebraic weight of $D$ as in the previous
lecture, and $\cE(D)$ is a horrifying expression which is a big product
of ``vertex terms'' corresponding to $\frac 23\,A\land A\land A$, ``edge
terms'' that look like
$ \epsilon^{ijk}\,\frac{x^k-y^k}{||x-y||^3}\,, $
corresponding to the singular integral-kernel of an inverse of $A\land
dA$,
and additional terms coming from the holonomy of $A$ along $K$.
\vspace{3mm}

\begin{example} When $m=2$, we get (roughly):
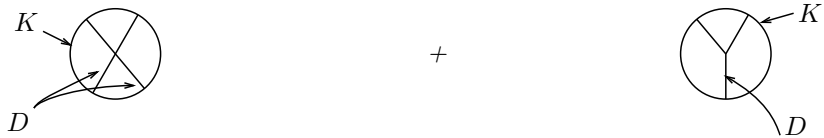
\begin{center}
  \begin{tabular}{*3c}
  \parbox[b]{3cm}{\mbox{\diag{6mm}44
  {\pictranslate{2 2}{
  \picline{1 60 polar}{1 240 polar}
  \picline{1 130 polar}{1 310 polar}
  \piccircle{0 0}{1}{}
  }
  \picvecline{0.4 2.5}{1.0 2.2}
  \picveccurve{0.2 0.8}{0.3 1.0}{1.3 1.3}{2.4 1.3}
  \picveccurve{0.2 0.8}{0.3 1.0}{0.8 1.2}{1.6 1.6}
  \pictext{$K$}{0.3 2.5}{-1 0}
  \pictext{$D$}{0.1 0.7}{-1 -1}
  }}}
  &+&
  \parbox[b]{3cm}{\mbox{\diag{6mm}44
  {\pictranslate{2 2}{
  \piccircle{0 0}{1}{}
  \picline{1 60 polar}{0 240 polar}
  \picline{1 130 polar}{0 310 polar}
  \picline{1 270 polar}{0 310 polar}
  }
  \picveccurve{3.2 0.2}{3.1 0.5}{2.9 1.0}{2.0 1.5}
  \picvecline{3.5 2.9}{2.8 2.7}
  \pictext{$K$}{3.6 2.9}{0 -0.5}
  \pictext{$D$}{3.3 0.2}{0 0}
  }}}
  \\[7mm]
  \\
  \parbox[t]{7cm}{4-fold integration
  along $K$ of 4 copies of
  \[
  \epsilon^{ijk}\,\frac{x^k-y^k}{||x-y||^3}\,
  \]
  }&&
  \parbox[t]{7cm}{3-fold integration  along $K$,
  1-fold integration  on $\bR^3$ of three copies of
  $\displaystyle 
  \epsilon^{ijk}\,\frac{x^k-y^k}{||x-y||^3}\,
  $ and summation over many indices.
  }
  \end{tabular}
\end{center}
\end{example}

\myitemm{\underline{Idea}\quad} There ought to be a direct ``differentiating
under the integral'' proof of the invariance of $\cZ_k(K)$. That
proof will use properties of the map $D\,\to\,W(D)$, and it seems
that all that those can be is $AS$, $IHX$, and $STU$.

Therefore,

\[
\cZ(K)\,=\,\sum_{m=0}^\infty\,\sum_D\,D\cdot\,\sumint_\text{labels}\,
\cE(D)\,\in\,\bar{\cA}\,=
\,\faktor{\left\{\vctext{@{}c@{}}{Feynman\\diagrams}{\small}%
\right\}}{\vcbox{\hbox{\small$\begin{array}{c}
AS\\IHX\\STU\end{array}$}}}
\]
ought to be an invariant (and a local computation near the double points
shows it to be a universal Vassiliev invariant).
%\newpage%\vspace{4mm plus 14mm}

{\parskip0pt
\noindent\underline{\bf Problems}\nopagebreak
\begin{enumerate}
\item $\displaystyle\int\cE(D)$ is naively {\bf divergent}.
\item When differentiating under the integral $\displaystyle\int\cE(D)$
gets \underline{worse}, and in fact, the result is non-zero. I.~e.,
we have an {\em ikke-invariant}.\footnote{The Danish work ``ikke'' is better
suited for our purposes than the English ``not'', as it is not homophonous
to ``knot''.}
\end{enumerate}

\noindent\underline{\bf Solutions}
\begin{enumerate}
\item Work harder to show convergence.
\item Add a local correction factor, in the same spirit as
of $Z(\infty)$ of the previous lecture (but very different).
\end{enumerate}

\noindent\underline{\bf Costs}
\begin{enumerate}
%\tracingmacros1
\item%\tracingmacros0
Lose some on elegance.
\item Reintroduce the Framing-Independence ($FI$) relation.
\end{enumerate}
}

\myitemm{\underline{\bf History}\quad} The first knot invariants
of this type were written (with no invariance proof) by Guadagnini,
Martellini and Mintchev~\cite{GuadagniniMartelliniMintchev:Aspects,
GuadagniniMartelliniMintchev:LinkInvariants}, following Witten's
discovery~\cite{Witten:JonesPolynomial} that the Jones polynomial can
be written in terms of the Chern-Simons quantum field theory. The same
invariants were independently written (together with an invariance proof)
somewhat later by the first author~\cite{Bar-Natan:PerturbativeCS},
who was later~\cite{Bar-Natan:Thesis} able to write a general
invariance proof in all orders of perturbation theory using only
the $STU$, $AS$, and $IHX$ relations, but assuming without proof the
convergence of all the integrals appearing.  Rather complete results on
perturbative invariants of 3-manifolds were obtained later by Axelrod and
Singer~\cite{AS1, AS2} and by Kontsevich (mostly unpublished). Recently
Bott and Taubes~\cite{BottTaubes:SelfLinking} reformulated the results
of~\cite{Bar-Natan:PerturbativeCS} in a much cleaner and prettier
topological language and suggested how this can be continued in higher
orders, and Thurston~\cite{Thurston:SeniorThesis} was able to complete
their work and write a proof of the Fundamental Theorem in these terms.

\vspace{4mm}

\def\bZ{{\bold Z}}

\subsection{Why are we not happy?}

\begin{enumerate}
\item In the Bott-Taubes formulation, much of the
mess is gone, and the integrals become evaluations of the volume
form of $S^2$ on various reasonably natural chains constructed
out of configuration space of points on the knot and elsewhere in $\bR^3$.
But still, this approach is very complicated and not quite the first
thing you would come up with.
\item The relationship with the Kontsevich-KZ approach
and with the Reshetikhin-Turaev invariants is still unclear.
\item The usual problem~--
what if you wanted to work over $\bZ/3\bZ$ or over $\bZ$?
\end{enumerate}

\lectitle{Algebra (sketch)}

\begin{remark} In this lecture we follow~\cite{Bar-Natan:NAT},
\cite{Cartier:Construction, Kassel:QuantumGroups, LeMurakami:Universal,
Piunikhin:Combinatorial, Stoimenow:HarrisonKohomologie}, and of course
Drinfel'd~\cite{Drinfeld:QuasiHopf, Drinfeld:GalQQ}.
\end{remark}

\subsection{Motivation from the Kontsevich-KZ integrals}
\label{motivation}

Recall the Kontsevich-KZ integrals, roughly given as
\[
\sumint_{\vctext{c}{pairings\\$P$}{\scriptsize}}\,D_P\,\cdot\,\prod\,
\frac{dz_i-dz'_i}{z_i-z'_i}\,\in\,\bar\cA^r
\]
and the Chern-Simons integrals
\[
\sumint_{\vctext{c}{Feynman\\diagrams\\$D$}{\scriptsize}}D\,\cdot\,
\cE(D)\,\in\,\bar\cA
\]
The KKZ integrals are ``multiplicative'' in the sense of
\eqref{multiplicative}.
The CS integrals probably also become multiplicative once
the knot is sufficiently ``stretched''.

Maybe these integrals can be evaluated by first deforming
the knot into some favorable position, and then by cutting along the time
slices and computing each piece separately?

For example, here is a better presentation for the trefoil:

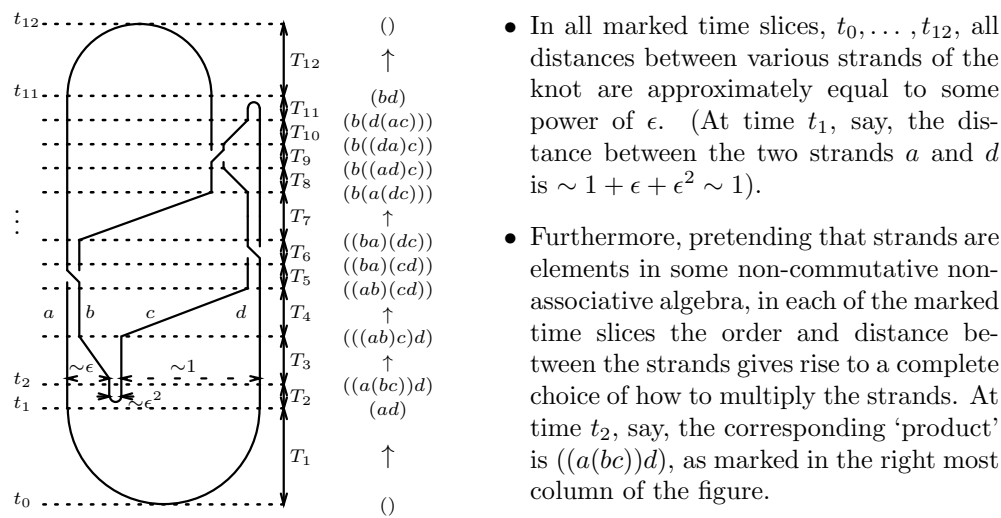
\[ 
\def\c#1{\parbox{0.65in}{\[ \scriptstyle (#1) \]}}
\def\ua{\parbox{0.65in}{\[ \uparrow \]}}
\def\sua{\parbox{0.65in}{\[ \scriptstyle\uparrow \]}}

        \setlength{\unitlength}{0.5\eepiclength}
        \begin{array}{c}  \hspace{-1.7mm}
                \raisebox{-8pt}{\begingroup\makeatletter\ifx\SetFigFont\undefined%
\gdef\SetFigFont#1#2#3#4#5{%
  \reset@font\fontsize{#1}{#2pt}%
  \fontfamily{#3}\fontseries{#4}\fontshape{#5}%
  \selectfont}%
\fi\endgroup%
\begin{picture}(4869,6195)(0,-10)
\thicklines
\put(3000.000,4986.000){\arc{150.000}{3.1416}{6.2832}}
\put(1575.000,5136.000){\arc{1800.000}{3.1416}{6.2832}}
\put(1275.000,1386.000){\arc{150.000}{6.2832}{9.4248}}
\put(1875.000,1236.000){\arc{2400.000}{6.2832}{9.4248}}
\path(1200,1611)(825,2136)
\path(675,1536)(675,2811)
\path(825,2136)(825,2811)(675,2961)
\path(675,2961)(675,4236)
\path(825,2961)(825,3336)
\path(2925,3636)(2925,3936)(2625,4236)(2625,4311)
\path(825,3336)(2475,3936)(2475,4311)(2625,4461)
\path(2625,4461)(2625,4536)(2925,4836)(2925,4986)
\path(3075,3636)(3075,4986)
\path(2475,4461)(2475,5136)
\path(675,5136)(675,4236)
\path(675,1536)(675,1236)
\path(3075,1536)(3075,1236)
\path(1200,1611)(1200,1386)
\path(1350,1611)(1350,1386)
\dottedline{150}(375,1236)(3375,1236)
\dottedline{150}(375,2136)(3375,2136)
\dottedline{150}(375,2736)(3375,2736)
\dottedline{150}(375,3036)(3375,3036)
\dottedline{150}(375,3336)(3375,3336)
\dottedline{150}(375,3936)(3375,3936)
\dottedline{150}(375,4236)(3375,4236)
\dottedline{150}(375,4536)(3375,4536)
\dottedline{150}(375,4836)(3375,4836)
\dottedline{150}(375,5136)(3375,5136)
\dottedline{150}(375,36)(3375,36)
\dottedline{150}(375,6036)(3375,6036)
\path(3412.500,1101.000)(3375.000,1236.000)(3337.500,1101.000)
\path(3375,1236)(3375,36)
\path(3337.500,171.000)(3375.000,36.000)(3412.500,171.000)
\path(3337.500,1371.000)(3375.000,1236.000)(3412.500,1371.000)
\path(3375,1236)(3375,1536)
\path(3412.500,1401.000)(3375.000,1536.000)(3337.500,1401.000)
\path(3337.500,1671.000)(3375.000,1536.000)(3412.500,1671.000)
\path(3375,1536)(3375,2136)
\path(3412.500,2001.000)(3375.000,2136.000)(3337.500,2001.000)
\path(3412.500,5901.000)(3375.000,6036.000)(3337.500,5901.000)
\path(3375,6036)(3375,5136)
\path(3337.500,5271.000)(3375.000,5136.000)(3412.500,5271.000)
\path(3337.500,2871.000)(3375.000,2736.000)(3412.500,2871.000)
\path(3375,2736)(3375,3036)
\path(3412.500,2901.000)(3375.000,3036.000)(3337.500,2901.000)
\path(3337.500,3171.000)(3375.000,3036.000)(3412.500,3171.000)
\path(3375,3036)(3375,3336)
\path(3412.500,3201.000)(3375.000,3336.000)(3337.500,3201.000)
\path(3337.500,4071.000)(3375.000,3936.000)(3412.500,4071.000)
\path(3375,3936)(3375,4236)
\path(3412.500,4101.000)(3375.000,4236.000)(3337.500,4101.000)
\path(3337.500,4371.000)(3375.000,4236.000)(3412.500,4371.000)
\path(3375,4236)(3375,4536)
\path(3412.500,4401.000)(3375.000,4536.000)(3337.500,4401.000)
\path(3337.500,4671.000)(3375.000,4536.000)(3412.500,4671.000)
\path(3375,4536)(3375,4836)
\path(3412.500,4701.000)(3375.000,4836.000)(3337.500,4701.000)
\path(3337.500,4971.000)(3375.000,4836.000)(3412.500,4971.000)
\path(3375,4836)(3375,5136)
\path(3412.500,5001.000)(3375.000,5136.000)(3337.500,5001.000)
\path(3337.500,2271.000)(3375.000,2136.000)(3412.500,2271.000)
\path(3375,2136)(3375,2736)
\path(3412.500,2601.000)(3375.000,2736.000)(3337.500,2601.000)
\path(3337.500,3471.000)(3375.000,3336.000)(3412.500,3471.000)
\path(3375,3336)(3375,3936)
\path(3412.500,3801.000)(3375.000,3936.000)(3337.500,3801.000)
\dottedline{150}(375,1536)(3375,1536)
\path(3075,3636)(3075,3261)
\path(1350,1611)(1350,2136)(2925,2736)
\path(2925,3636)(2925,3261)(3075,3111)(3075,1536)
\path(2925,3111)(2925,2736)
\path(1050,1386)(1200,1386)
\path(1065.000,1348.500)(1200.000,1386.000)(1065.000,1423.500)
\path(1500,1386)(1350,1386)
\path(1485.000,1423.500)(1350.000,1386.000)(1485.000,1348.500)
\path(795.000,1641.000)(675.000,1611.000)(795.000,1581.000)
\dashline{90.000}(675,1611)(1200,1611)
\path(1080.000,1581.000)(1200.000,1611.000)(1080.000,1641.000)
\path(1470.000,1641.000)(1350.000,1611.000)(1470.000,1581.000)
\dashline{90.000}(1350,1611)(3075,1611)
\path(2955.000,1581.000)(3075.000,1611.000)(2955.000,1641.000)
\put(3450,3561){\makebox(0,0)[lb]{\smash{{\mathmode{\scriptstyle T_7}}}}}
\put(3450,561){\makebox(0,0)[lb]{\smash{{\mathmode{\scriptstyle T_1}}}}}
\put(3450,1311){\makebox(0,0)[lb]{\smash{{\mathmode{\scriptstyle T_2}}}}}
\put(3450,1761){\makebox(0,0)[lb]{\smash{{\mathmode{\scriptstyle T_3}}}}}
\put(3450,2361){\makebox(0,0)[lb]{\smash{{\mathmode{\scriptstyle T_4}}}}}
\put(3450,2811){\makebox(0,0)[lb]{\smash{{\mathmode{\scriptstyle T_5}}}}}
\put(3450,3111){\makebox(0,0)[lb]{\smash{{\mathmode{\scriptstyle T_6}}}}}
\put(3450,4011){\makebox(0,0)[lb]{\smash{{\mathmode{\scriptstyle T_8}}}}}
\put(3450,4311){\makebox(0,0)[lb]{\smash{{\mathmode{\scriptstyle T_9}}}}}
\put(3450,4611){\makebox(0,0)[lb]{\smash{{\mathmode{\scriptstyle T_{10}}}}}}
\put(3450,4911){\makebox(0,0)[lb]{\smash{{\mathmode{\scriptstyle T_{11}}}}}}
\put(3450,5511){\makebox(0,0)[lb]{\smash{{\mathmode{\scriptstyle T_{12}}}}}}
\put(1950,1686){\makebox(0,0)[lb]{\smash{{\mathmode{\scriptstyle \sim 1}}}}}
\put(1425,1236){\makebox(0,0)[lb]{\smash{{\mathmode{\scriptstyle \sim \epsilon^2}}}}}
\put(690,1686){\makebox(0,0)[lb]{\smash{{\mathmode{\scriptstyle \sim \epsilon}}}}}
\put(0,6036){\makebox(0,0)[lb]{\smash{{\mathmode{\scriptstyle t_{12}}}}}}
\put(0,5136){\makebox(0,0)[lb]{\smash{{\mathmode{\scriptstyle t_{11}}}}}}
\put(0,36){\makebox(0,0)[lb]{\smash{{\mathmode{\scriptstyle t_0}}}}}
\put(0,1236){\makebox(0,0)[lb]{\smash{{\mathmode{\scriptstyle t_1}}}}}
\put(0,1536){\makebox(0,0)[lb]{\smash{{\mathmode{\scriptstyle t_2}}}}}
\put(0,3411){\makebox(0,0)[lb]{\smash{{\mathmode{\vdots}}}}}
\put(375,2361){\makebox(0,0)[lb]{\smash{{\mathmode{\scriptstyle a}}}}}
\put(900,2361){\makebox(0,0)[lb]{\smash{{\mathmode{\scriptstyle b}}}}}
\put(2775,2361){\makebox(0,0)[lb]{\smash{{\mathmode{\scriptstyle d}}}}}
\put(1650,2361){\makebox(0,0)[lb]{\smash{{\mathmode{\scriptstyle c}}}}}
\put(3900,36){\makebox(0,0)[lb]{\smash{{\mathmode{\c{}}}}}}
\put(3900,1236){\makebox(0,0)[lb]{\smash{{\mathmode{\c{ad}}}}}}
\put(3900,1536){\makebox(0,0)[lb]{\smash{{\mathmode{\c{(a(bc))d}}}}}}
\put(3900,2136){\makebox(0,0)[lb]{\smash{{\mathmode{\c{((ab)c)d}}}}}}
\put(3900,2736){\makebox(0,0)[lb]{\smash{{\mathmode{\c{(ab)(cd)}}}}}}
\put(3900,3336){\makebox(0,0)[lb]{\smash{{\mathmode{\c{(ba)(dc)}}}}}}
\put(3900,3036){\makebox(0,0)[lb]{\smash{{\mathmode{\c{(ba)(cd)}}}}}}
\put(3900,3936){\makebox(0,0)[lb]{\smash{{\mathmode{\c{b(a(dc))}}}}}}
\put(3900,4236){\makebox(0,0)[lb]{\smash{{\mathmode{\c{b((ad)c)}}}}}}
\put(3900,4536){\makebox(0,0)[lb]{\smash{{\mathmode{\c{b((da)c)}}}}}}
\put(3900,4836){\makebox(0,0)[lb]{\smash{{\mathmode{\c{b(d(ac))}}}}}}
\put(3900,5136){\makebox(0,0)[lb]{\smash{{\mathmode{\c{bd}}}}}}
\put(3900,6036){\makebox(0,0)[lb]{\smash{{\mathmode{\c{}}}}}}
\put(3900,636){\makebox(0,0)[lb]{\smash{{\mathmode{\ua}}}}}
\put(3900,5586){\makebox(0,0)[lb]{\smash{{\mathmode{\ua}}}}}
\put(3900,1836){\makebox(0,0)[lb]{\smash{{\mathmode{\sua}}}}}
\put(3900,2436){\makebox(0,0)[lb]{\smash{{\mathmode{\sua}}}}}
\put(3900,3636){\makebox(0,0)[lb]{\smash{{\mathmode{\sua}}}}}
\end{picture} }
                \hspace{-1.9mm}
        \end{array}
\qquad
        \parbox{2.8in}{\begin{itemize}
        \item In all marked time slices, $t_0,\ldots,t_{12}$, all distances
               between various strands of the knot are approximately
                equal to some power of $\epsilon$. (At time $t_1$, say, the
               distance between the two strands $a$ and $d$ is
                $\sim 1+\epsilon+\epsilon^2\sim 1$).
        \item Furthermore, pretending that strands are elements in some
                non-commutative non-associative algebra, in each of the
                marked time slices the order and distance between the
                strands gives rise to a complete choice of how to
                multiply the strands.  At time $t_2$, say, the
                corresponding `product' is $((a(bc))d)$, as marked in
                the right most column of the figure.
        \end{itemize}}
\]

\par\noindent
\begin{myitemize}
\item In each of the time intervals $T_1,\ldots,T_{12}$ only one change
        occurs to the `product' corresponding to the strands, and only
        three types of changes occur:
        \begin{enumerate}
        \def\theenumi{(\roman{enumi})} \setcounter{enumi}{1}
        \def\theenum{\theenumi\refstepcounter{enumi}}
        \item[\theenum]
                {\em Pair creation (annihilation)}, in which a pair of
                {\em neighboring strands} is created (or annihilated).
                {\em Neighboring strands} are strands for which
                the distance between them is smaller than the distance
                between them and any other strand.  (intervals $T_1$,
                $T_2$, $T_{11}$, and $T_{12}$).
        \item[\theenum]
                {\em Braiding morphism}, in which two neighboring strands
                are braided. (intervals $T_5$, $T_6$, and $T_9$).
        \item[\theenum]
                {\em Associativity morphism}, in which the associative
                law is applied once. (intervals $T_3$, $T_4$, $T_7$,
                $T_8$, and $T_{10}$).
        \end{enumerate}
\end{myitemize}

In this presentation, the computation of the KKZ integral in each time
interval is relatively simple. For example,
\def\T#1{\myitem[1cm]{$T_{#1}$:}}

\T{12} $Id$ by the $FI$ relation.

\T{11} The left two strands are too far to matter (the contribution of
diagrams with chords ending on them is too small), and the rest is as
in $T_{12}$.

\T{10} The left strand is irrelevant, so we are left with understanding
what is the Kontsevich integral on
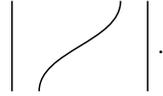
\[
\diag{6mm}32{
\picline{3 0}{3 2}
\picline{0 0}{0 2}
\piccurve{0.6 0}{0.6 0.9}{2.4 1.1}{2.4 2}
}\,.
\]
The result will be some $\Phi$ of the following shape
\[
\Phi\,=:\,Z\left(\,
\diag{4mm}32{
\picline{3 0}{3 2}
\picline{0 0}{0 2}
\piccurve{0.6 0}{0.6 0.9}{2.4 1.1}{2.4 2}
}\,
\right)\,=\,\sum\,\langle\text{ coefficients }\rangle\enspace
\langle\text{ diagrams like \enspace\diag{5mm}{3}{2}{
\picmultigraphics{3}{1.5 0}{\picvecline{0 0}{0 2}}
\picline{1.5 0.4}{3 0.4}
\picline{1.5 1.6}{3 1.6}
\picline{1.5 1.3}{0 1.3}
\picline{1.5 0.7}{0 0.7}
}\enspace}\rangle
\]

\myitem[5cm]{\hbox to 1cm{$T_9$:\hss}\vbox to 1em{\hbox{\kern5mm\diag{4mm}{6}{6}{
\picellipse{3 1}{3 1}{}
\picellipse{3 5}{3 1}{}
\picline{0 1}{0 5}
\picline{6 1}{6 5}
\pictranslate{3 3}{
  \picrotate{2 3 atan}{
    \picscale{13 q 1}{
      \picellipsearc{0 0}{1 0.70}{180 0}
    }
  }
}
\piclinedash{0.2 0.1}{0.1}
\pictranslate{3 3}{
  \picscale{1 -1}{
    \picrotate{2 3 atan}{
      \picscale{13 q 1}{
	\picellipsearc{0 0}{1 0.70}{180 0}
      }
    }
  }
}
}}\vss}}
After ignoring the left most and the right most strands, we
have a braiding. We may assume that the two braiding strands are
parametrized uniformly around a cylinder as in the figure, and then
$(dz-dz')/(z-z')=d\theta=dt$, so we can easily compute the value of
the Kontsevich integral by directly integrating over the simplex
$t_{\text{\scriptsize min}}< t_1<\ldots<t_m<t_{\text{\scriptsize
max}}$. We get

\myitem[1cm]{}
\[
Z\,\left(\,\diag{6mm}{1}{1}{
\picmultivecline{-5 1 -1 0}{1 0}{0 1}
\picmultivecline{-5 1 -1 0}{0 0}{1 1}
}\,\right)\enspace = \enspace
% \left(\,\exp\frac{1}{2}\diag{6mm}{1}{1}
% {\picvecline{0 0}{0 1}
% \picvecline{1 0}{1 1}
% \picline{0 0.5}{1 0.5}
% }\,\right)\,\cdot\,
% \diag{6mm}{1}{1}{
% \picvecline{0 0}{1 1}
% \picvecline{1 0}{0 1}
% \picfillgraycol{0}
% }\\
% &=&\enspace
\diag{6mm}{1}{1}{
\picvecline{0 0}{1 1}
\picvecline{1 0}{0 1}
}\,+\,\frac1{2\cdot 1!}\,
\diag{6mm}{1}{1}{
\picvecline{0 0}{1 1}
\picvecline{1 0}{0 1}
\picline{0.25 0.25}{0.75 0.25}
}\,+\,\frac1{2^2\cdot 2!}\,
\diag{6mm}{1}{1}{
\picvecline{0 0}{1 1}
\picvecline{1 0}{0 1}
\picline{0.15 0.15}{0.85 0.15}
\picline{0.28 0.28}{0.72 0.28}
}\,+\,\frac1{2^3\cdot 3!}\,
\diag{6mm}{1}{1}{
\picvecline{0 0}{1 1}
\picvecline{1 0}{0 1}
\picline{0.12 0.12}{0.88 0.12}
\picline{0.24 0.24}{0.76 0.24}
\picline{0.36 0.36}{0.64 0.36}
}\,+\,\dots\enspace
 = \enspace R\,\cdot\,
\diag{6mm}{1}{1}{
\picvecline{0 0}{1 1}
\picvecline{1 0}{0 1}
\picfillgraycol{0}
}\,,
\]
% \\[2mm]
where $R\,=\,\exp\left(\frac{1}{2}\diag{6mm}{1}{1}
{\picvecline{0 0}{0 1}
\picvecline{1 0}{1 1}
\picline{0 0.5}{1 0.5}
}\,\right)$\,.
% \end{eqnarray*}

\T8 Here we have a variation of the inverse of $\Phi$, obtained by
``placing $\Phi^{-1}$ on strands 2, 3, and 4''. Symbolically, we write it
as $(\Phi^{234})^{-1}$.

\T7 In this time interval, the rightmost two strands are too close to each
other for the other strands to tell them apart. In the Kontsevich integral,
chords whose right end is on the $3^\text{rd}$ strand appear with the
same weight as chords whose right end is on the $4^\text{th}$ strand, and
with the same weight as the corresponding chords that appear in the
computation of $\Phi$. This means that
\[
Z\,\left(\,
\diag{4mm}{3.4}2{
\picline{3 0}{3 2}
\picline{3.3 0}{3.3 2}
\picline{0 0}{0 2}
\piccurve{0.6 0}{0.6 0.9}{2.4 1.1}{2.4 2}
%\piclinedash{0.2 0.2}{0.1}
%\picline{
}\right)\,=\,
(1\ot1\ot\Delta)(\Phi)
\,,
\]
where the operation $\Delta$ doubles a strand and sums over all possible ways
of `lifting' the chords that were connected to it to the two offspring
chords.

\subsection{Relations between $R$ and $\Phi$.}

Now try to reconstruct the Kontsevich integral algebraically. As is clear
from section~\ref{motivation}, to know the Kontsevich integral on all knots
(and, in fact, on all ``parenthesized tangles'', the kind of objects that
appear between any two time slices in knot presentations such as in
section~\ref{motivation}), it is enough to compute only two quantities, $R$
and $\Phi$, once and for all. As the computation of $\Phi$ appears hard,
let's just assign an arbitrary value to it (and to $R$ as well), and check
what axioms these $R$ and $\Phi$ have to satisfy so that the computation
algorithm implicitly defined in section~\ref{motivation} really does yield
an invariant. One can check that the axioms are as follows:

\subsubsection{Axioms for $R$ and $\Phi$.}

\begin{figure}[htpb]
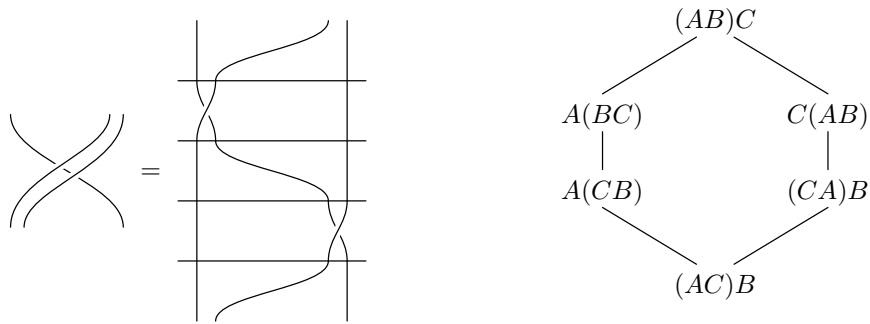

\[
  \diag{5mm}{3}{3}{\picscale{0.6 0.6}{
    \piccurve{5 0}{5 2}{0 3}{0 5}
    \picmulticurve{-7 1 -1 0}{0.6 0}{0.6 1.7}{5 2.7}{5 5}
    \picmulticurve{-7 1 -1 0}{0 0}{0 2.3}{4.4 3.3}{4.4 5}
  }}\,=\, \diag{5mm}{5}{8}{
    \pic@pss{/default@brdrm [-7 1 -1 0] D} \picscale{1 0.8}{
       \lbraid{0.75 7}{0.5 2} \lbraid{4.25 3}{0.5 2}
       \picline{0.5 0}{0.5 6} \picline{0.5 8}{0.5 10}
       \picline{4.5 0}{4.5 2} \picline{4.5 4}{4.5 10}
      \picmultigraphics{2}{0 8}{\piccurve{1.0 0}{1.0 0.9}{4.0 1.1}{4.0 2.0}}
      \piccurve{4.0 4}{4.0 4.9}{1.0 5.1}{1.0 6.0}
      \picmultigraphics{4}{0 2}{\picline{0 2}{5 2}}
    }
  }
  \hspace{20mm}
  \diag{5mm}{10}{10}{
    \picputtext{5 9}{$(AB)C$} \picputtext{2.0 6.5}{$A(BC)$}
    \picputtext{2.0 4.5}{$A(CB)$} \picputtext{8 4.5}{$(CA)B$}
    \picputtext{8 6.5}{$C(AB)$} \picputtext{5 2}{$(AC)B$}
    \picline{4.5 8.5}{2 7} \picline{5.5 8.5}{8 7} \picline{8 5}{8 6}
    \picline{2 5}{2 6} \picline{2 4}{4.5 2.5} \picline{8 4}{5.5 2.5}
  }
\]
\caption{The braid equality leading to the ``$+$'' hexagon axiom, and the
  associativities hexagon that gave it its name. To get the ``$-$'' hexagon
  axiom, simply flip all crossings.
}
\label{HexagonAxioms} \end{figure}

\begin{itemize}

\item The pentagon axiom:
\begin{equation}
  \Phi^{123}\cdot(1\ot\Delta\ot 1)(\Phi)\cdot\Phi^{234}\cdot
  (1\ot 1\ot \Delta)(\Phi^{-1})\cdot(\Delta\ot1\ot1)(\Phi^{-1})\,=\,1
\tag{$\pentagon$} \end{equation}
Needed because of the funny presentation of the trivial braid shown in
figure~\ref{PentagonAxiom}.

\begin{figure}[htpb]
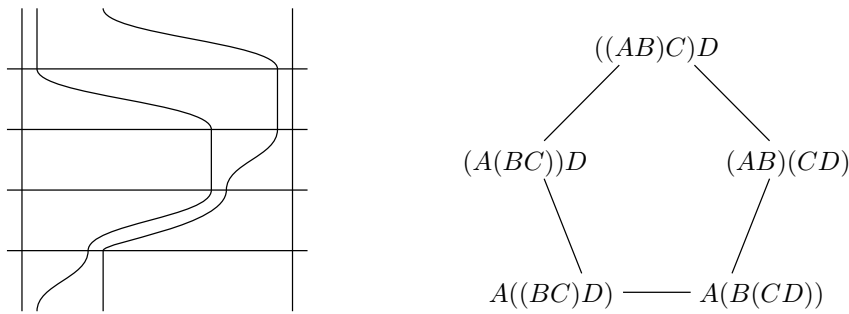

\[
  \diag{5mm}{8}{8}{
   %\pictranslate{0 1}{
   \picscale{0.8 0.8}{
    \picmultigraphics{4}{0 2}{\picline{0 2}{10 2}}
    \picline{0.5 0}{0.5 10} \picline{9.5 0}{9.5 10} \picline{3.2 0}{3.2 2}
    \picmultigraphics{2}{6.3 4}{ \piccurve{1.0 0}{1.0 0.9}{2.7 1.1}{2.7 2.0} }
    \piccurve{2.7 2.0}{2.7 2.9}{6.8 3.1}{6.8 4.0}
    \piccurve{3.2 2.0}{3.2 2.3}{7.3 2.8}{7.3 4.0}
    \picline{6.8 4.0}{6.8 6.0} \piccurve{6.8 6.0}{6.8 6.9}{1.0 7.1}{1.0 8.0}
    \piccurve{9.0 8.0}{9.0 8.9}{3.2 9.1}{3.2 10.0} \picline{1.0 8.0}{1.0 10.0}
    \picline{9.0 8.0}{9.0 6.0}
   }
   %}
  }
  \hspace{20mm}
  \diag{5mm}{10}{8}{
   \pictranslate{0 -2}{
    \picputtext[m]{5 9}{$((AB)C)D$} \picputtext[m]{1.5 6.0}{$(A(BC))D$}
    \picputtext[m]{8.5 6.0}{$(AB)(CD)$} \picputtext[m]{2.2 2.5}{$A((BC)D)$}
    \picputtext[m]{7.8 2.5}{$A(B(CD))$}
    \picline{6 8.5}{8 6.5} \picline{4 8.5}{2 6.5} \picline{3.0 3.0}{2 5.5}
    \picline{7.0 3.0}{8 5.5} \picline{4.1 2.5}{5.9 2.5}
   }
  }
\]
\caption{
  Left: The presentation of the trivial braid that leads to the
  pentagon axiom. Right: The reason for the name ``pentagon''. When the
  arrangements at the indicated time slices of the strands on the left
  are written as `products', we get a pentagon of associativities.
}
\label{PentagonAxiom} \end{figure}

\item The hexagon axioms: (See figure~\ref{HexagonAxioms})
\begin{equation}
  (\Delta\ot1)R^{\pm1}\,=\,\Phi\cdot(R^{23})^{\pm1}\cdot
  (\Phi^{-1})^{132}\cdot(R^{13})^{\pm1}\cdot\Phi^{312}
\tag{$\hexagons$} \end{equation}

\end{itemize}

\subsubsection{Automatic relations between $R$ and $\Phi$.}
\label{FreeRelations}
There are two other types of relations, that $R$ and $\Phi$ satisfy
automatically (and hence do not impose constrains on their possible
values).

\begin{itemize}
\item Locality in space relations: Events that happen far away from each
other commute. For example:
\[ 
        \setlength{\unitlength}{0.5\eepiclength}
        \begin{array}{c}  \hspace{-1.7mm}
                \raisebox{-8pt}{\begingroup\makeatletter\ifx\SetFigFont\undefined%
\gdef\SetFigFont#1#2#3#4#5{%
  \reset@font\fontsize{#1}{#2pt}%
  \fontfamily{#3}\fontseries{#4}\fontshape{#5}%
  \selectfont}%
\fi\endgroup%
\begin{picture}(5723,1240)(0,-10)
\thicklines
\path(3611,13)(3611,613)
\path(3311,13)(3311,613)
\path(3611,1213)	(3611.000,1164.014)
	(3611.000,1138.000)

\path(3611,1138)	(3613.903,1101.100)
	(3611.000,1063.000)

\path(3611,1063)	(3583.133,1022.072)
	(3540.605,980.995)
	(3495.775,943.421)
	(3461.000,913.000)

\path(3461,913)	(3426.225,882.579)
	(3381.395,845.005)
	(3338.867,803.928)
	(3311.000,763.000)

\path(3311,763)	(3308.097,724.900)
	(3311.000,688.000)

\path(3311,688)	(3311.000,661.986)
	(3311.000,613.000)

\path(3311,1213)	(3311.000,1164.014)
	(3311.000,1138.000)

\path(3311,1138)	(3308.097,1101.100)
	(3311.000,1063.000)

\path(3311,1063)	(3333.564,1031.729)
	(3386.000,988.000)

\path(3611,613)	(3611.000,661.986)
	(3611.000,688.000)

\path(3611,688)	(3613.903,724.900)
	(3611.000,763.000)

\path(3611,763)	(3588.436,794.271)
	(3536.000,838.000)

\path(5711,13)(5711,1213)
\path(5411,613)(5411,1213)
\drawline(4811,613)(4811,613)
\path(4811,13)(4811,1213)
\path(5411,613)	(5414.075,564.561)
	(5411.000,538.000)

\path(5411,538)	(5388.413,486.637)
	(5353.044,431.086)
	(5309.152,372.741)
	(5261.000,313.000)
	(5212.848,253.259)
	(5168.956,194.914)
	(5133.587,139.363)
	(5111.000,88.000)

\path(5111,88)	(5107.925,61.446)
	(5111.000,13.000)

\path(1512,1212)(1512,12)
\path(1812,612)(1812,12)
\drawline(2412,612)(2412,612)
\path(2412,1212)(2412,12)
\path(12,1212)(12,612)
\path(312,1212)(312,612)
\path(12,12)	(12.000,60.986)
	(12.000,87.000)

\path(12,87)	(9.097,123.900)
	(12.000,162.000)

\path(12,162)	(39.867,202.928)
	(82.395,244.005)
	(127.225,281.579)
	(162.000,312.000)

\path(162,312)	(196.775,342.421)
	(241.605,379.995)
	(284.133,421.072)
	(312.000,462.000)

\path(312,462)	(314.903,500.100)
	(312.000,537.000)

\path(312,537)	(312.000,563.014)
	(312.000,612.000)

\path(312,12)	(312.000,60.986)
	(312.000,87.000)

\path(312,87)	(314.903,123.900)
	(312.000,162.000)

\path(312,162)	(289.436,193.271)
	(237.000,237.000)

\path(12,612)	(12.000,563.014)
	(12.000,537.000)

\path(12,537)	(9.097,500.100)
	(12.000,462.000)

\path(12,462)	(34.564,430.729)
	(87.000,387.000)

\path(1812,612)	(1808.925,660.439)
	(1812.000,687.000)

\path(1812,687)	(1834.587,738.363)
	(1869.956,793.914)
	(1913.848,852.259)
	(1962.000,912.000)
	(2010.152,971.741)
	(2054.044,1030.086)
	(2089.413,1085.637)
	(2112.000,1137.000)

\path(2112,1137)	(2115.075,1163.554)
	(2112.000,1212.000)

\put(2712,612){\makebox(0,0)[lb]{\smash{{\mathmode{=}}}}}
\end{picture} }
                \hspace{-1.9mm}
        \end{array}

  \qquad\mbox{or}\qquad
  R^{12}\Phi^{345}=\Phi^{345}R^{12}.
\]
These relations are a consequence of the fact that in chord diagrams the
chords are not time-ordered; only their ends are ordered. And thus chords
whose ends are on different strands always commute.
\item Locality in scale relations: Events that happen at different scales
commute. The third Reidemeister move is an example for such a relation:
\[ 
        \setlength{\unitlength}{0.5\eepiclength}
        \begin{array}{c}  \hspace{-1.7mm}
                \raisebox{-8pt}{\begingroup\makeatletter\ifx\SetFigFont\undefined%
\gdef\SetFigFont#1#2#3#4#5{%
  \reset@font\fontsize{#1}{#2pt}%
  \fontfamily{#3}\fontseries{#4}\fontshape{#5}%
  \selectfont}%
\fi\endgroup%
\begin{picture}(2424,1239)(0,-10)
\thicklines
\path(12,912)(12,1212)
\path(2412,12)(2412,312)
\path(912,912)(912,987)
\path(762,1137)(762,1212)
\path(1662,12)(1662,87)
\path(1512,237)(1512,312)
\path(162,12)	(168.126,72.982)
	(174.741,125.498)
	(182.065,170.425)
	(190.316,208.643)
	(210.483,268.466)
	(237.000,312.000)

\path(237,312)	(298.885,361.332)
	(339.549,382.430)
	(383.258,401.475)
	(427.434,418.710)
	(469.501,434.380)
	(537.000,462.000)

\path(537,462)	(604.499,489.620)
	(646.566,505.290)
	(690.742,522.525)
	(734.451,541.570)
	(775.115,562.668)
	(837.000,612.000)

\path(837,612)	(863.511,655.534)
	(883.676,715.357)
	(891.928,753.575)
	(899.253,798.502)
	(905.871,851.018)
	(912.000,912.000)

\path(12,12)	(18.126,72.982)
	(24.741,125.498)
	(32.065,170.425)
	(40.316,208.643)
	(60.483,268.466)
	(87.000,312.000)

\path(87,312)	(148.885,361.332)
	(189.549,382.430)
	(233.257,401.475)
	(277.434,418.710)
	(319.501,434.380)
	(387.000,462.000)

\path(387,462)	(454.499,489.620)
	(496.566,505.290)
	(540.742,522.525)
	(584.451,541.570)
	(625.115,562.668)
	(687.000,612.000)

\path(687,612)	(713.511,655.534)
	(733.676,715.357)
	(741.928,753.575)
	(749.253,798.502)
	(755.871,851.018)
	(762.000,912.000)

\path(912,12)	(907.915,72.982)
	(902.758,125.498)
	(896.309,170.425)
	(888.349,208.643)
	(867.015,268.466)
	(837.000,312.000)

\path(837,312)	(805.476,337.515)
	(760.984,357.349)
	(698.250,373.258)
	(658.394,380.290)
	(612.000,387.000)

\path(312,537)	(264.287,562.121)
	(223.489,584.543)
	(160.001,623.048)
	(87.000,687.000)

\path(87,687)	(67.548,718.906)
	(49.736,763.526)
	(31.806,826.133)
	(22.247,865.829)
	(12.000,912.000)

\path(762,912)	(757.643,960.086)
	(762.000,987.000)

\path(762,987)	(791.594,1025.306)
	(837.000,1062.004)
	(882.406,1098.700)
	(912.000,1137.000)

\path(912,1137)	(916.357,1163.914)
	(912.000,1212.000)

\path(1662,312)	(1668.129,372.982)
	(1674.747,425.498)
	(1682.072,470.425)
	(1690.324,508.643)
	(1710.489,568.466)
	(1737.000,612.000)

\path(1737,612)	(1798.885,661.332)
	(1839.549,682.430)
	(1883.258,701.475)
	(1927.434,718.710)
	(1969.501,734.380)
	(2037.000,762.000)

\path(2037,762)	(2104.499,789.620)
	(2146.566,805.290)
	(2190.742,822.525)
	(2234.451,841.570)
	(2275.115,862.668)
	(2337.000,912.000)

\path(2337,912)	(2363.511,955.534)
	(2383.676,1015.358)
	(2391.928,1053.575)
	(2399.253,1098.503)
	(2405.871,1151.018)
	(2412.000,1212.000)

\path(1512,312)	(1518.129,372.982)
	(1524.747,425.498)
	(1532.072,470.425)
	(1540.324,508.643)
	(1560.489,568.466)
	(1587.000,612.000)

\path(1587,612)	(1648.885,661.332)
	(1689.549,682.430)
	(1733.258,701.475)
	(1777.434,718.710)
	(1819.501,734.380)
	(1887.000,762.000)

\path(1887,762)	(1954.499,789.620)
	(1996.566,805.290)
	(2040.742,822.525)
	(2084.451,841.570)
	(2125.115,862.668)
	(2187.000,912.000)

\path(2187,912)	(2213.511,955.534)
	(2233.676,1015.358)
	(2241.928,1053.575)
	(2249.253,1098.503)
	(2255.871,1151.018)
	(2262.000,1212.000)

\path(2412,312)	(2407.915,372.982)
	(2402.758,425.498)
	(2396.309,470.425)
	(2388.349,508.643)
	(2367.015,568.466)
	(2337.000,612.000)

\path(2337,612)	(2305.476,637.515)
	(2260.984,657.349)
	(2198.250,673.258)
	(2158.394,680.290)
	(2112.000,687.000)

\path(1812,837)	(1764.287,862.121)
	(1723.489,884.543)
	(1660.001,923.048)
	(1587.000,987.000)

\path(1587,987)	(1567.548,1018.906)
	(1549.736,1063.526)
	(1531.806,1126.133)
	(1522.247,1165.829)
	(1512.000,1212.000)

\path(1512,12)	(1507.650,60.086)
	(1512.000,87.000)

\path(1512,87)	(1541.594,125.307)
	(1587.000,162.007)
	(1632.406,198.704)
	(1662.000,237.000)

\path(1662,237)	(1666.350,263.914)
	(1662.000,312.000)

\put(1062,537){\makebox(0,0)[lb]{\smash{{\mathmode{=}}}}}
\end{picture} }
                \hspace{-1.9mm}
        \end{array}
,
  \qquad\mbox{or}\qquad
  (\Delta\ot 1)R\cdot R^{12} = R^{12}\cdot(\Delta\ot 1)R.
\]
These relations are consequences of the $4T$ relation, written in the form
\[ 
        \setlength{\unitlength}{0.5\eepiclength}
        \begin{array}{c}  \hspace{-1.7mm}
                \raisebox{-8pt}{\begingroup\makeatletter\ifx\SetFigFont\undefined%
\gdef\SetFigFont#1#2#3#4#5{%
  \reset@font\fontsize{#1}{#2pt}%
  \fontfamily{#3}\fontseries{#4}\fontshape{#5}%
  \selectfont}%
\fi\endgroup%
\begin{picture}(5425,939)(0,-10)
\thicklines
\path(12,12)(12,912)
\path(42.000,792.000)(12.000,912.000)(-18.000,792.000)
\path(312,12)(312,912)
\path(342.000,792.000)(312.000,912.000)(282.000,792.000)
\path(912,12)(912,912)
\path(942.000,792.000)(912.000,912.000)(882.000,792.000)
\path(1513,12)(1513,912)
\path(1543.000,792.000)(1513.000,912.000)(1483.000,792.000)
\path(1813,12)(1813,912)
\path(1843.000,792.000)(1813.000,912.000)(1783.000,792.000)
\path(2413,12)(2413,912)
\path(2443.000,792.000)(2413.000,912.000)(2383.000,792.000)
\path(3013,12)(3013,912)
\path(3043.000,792.000)(3013.000,912.000)(2983.000,792.000)
\path(3313,12)(3313,912)
\path(3343.000,792.000)(3313.000,912.000)(3283.000,792.000)
\path(3913,12)(3913,912)
\path(3943.000,792.000)(3913.000,912.000)(3883.000,792.000)
\path(4513,12)(4513,912)
\path(4543.000,792.000)(4513.000,912.000)(4483.000,792.000)
\path(4813,12)(4813,912)
\path(4843.000,792.000)(4813.000,912.000)(4783.000,792.000)
\path(5413,12)(5413,912)
\path(5443.000,792.000)(5413.000,912.000)(5383.000,792.000)
\path(12,612)(312,612)
\path(1512,612)(1812,612)
\path(3012,312)(3312,312)
\path(4512,312)(4812,312)
\path(12,312)(912,312)
\path(1812,312)(2412,312)
\path(3012,612)(3912,612)
\path(4812,612)(5412,612)
\put(1062,387){\makebox(0,0)[lb]{\smash{{\mathmode{+}}}}}
\put(2562,387){\makebox(0,0)[lb]{\smash{{\mathmode{=}}}}}
\put(4062,387){\makebox(0,0)[lb]{\smash{{\mathmode{+}}}}}
\end{picture} }
                \hspace{-1.9mm}
        \end{array}
. \]
\end{itemize}

\subsection{An aside on quasi-Hopf algebras.}
Just for the sake of completeness, let us spend just around one page on
recalling where Drinfel'd first found the pentagon and the hexagon
equations~\cite{Drinfeld:QuasiHopf}. The context is superficially very
different, but the equations turn out to be exactly the same (though they
are about different kinds of objects). The technique we use for solving
these equations in section~\ref{solve} is nothing but Drinfel'd's technique
of~\cite{Drinfeld:GalQQ}, adopted to our situation.

If $R_1$ and $R_2$ are representations of some algebra $A$,
$R_1\ot R_2$ is a representation of $A\ot A$, but
(in general) not of $A$. To be able to take the tensor product
of representations, we need to have a morphism $\Delta:A\enspace \to A\ot A$
called ``the co-product''.

If we want the tensor product of representations to
be associative, $\Delta$ must be ``co-associative'':
\[
  (\Delta\ot 1)\Delta a\,=\, (1\ot\Delta)\Delta a\qquad\forall\,a\in A
\]
If this happens, $(A,\Delta)$ is a ``Hopf-Algebra''.

In~\cite{Drinfeld:QuasiHopf}, Drinfel'd suggested relax the condition
of associativity.  Instead of
\begin{itemize}
\item[*] ``$R_1\ot(R_2\ot R_3)$ and $(R_1\ot R_2)\ot R_3$
are the same''
\end{itemize}
we only require (roughly)
\begin{itemize}
\item[*] ``$R_1\ot(R_2\ot R_3)$ and $(R_1\ot R_2)\ot R_3$
are equivalent in a functorial way''.
\end{itemize}
This leads to a relaxation of the co-associativity condition on $A$:
\[
  \exists\,\Phi\in A^{\ot 3}\quad\text{s.t.}\quad\forall\,a\in A
  \quad (\Delta\ot 1)\Delta a
  =\Phi^{-1}\,\bigl((1\ot\Delta)\Delta a\bigr)\,\Phi
\]
But then $\Phi$ needs to have some properties, if diagrams
like
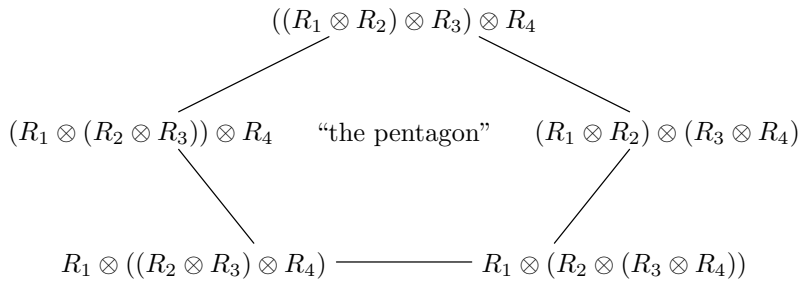
\[
\diag{5mm}{20}{10}{
\pictranslate{0 -0.6}{
\picputtext[m]{10 9}{$((R_1\ot R_2)\ot R_3)\ot R_4$}
\picputtext[m]{3 6.0}{$(R_1\ot(R_2\ot R_3))\ot R_4$}
\picputtext[m]{17 6.0}{$(R_1\ot R_2)\ot(R_3\ot R_4)$}
\picputtext[m]{4.4 2.5}{$R_1\ot((R_2\ot R_3)\ot R_4)$}
\picputtext[m]{15.6 2.5}{$R_1\ot(R_2\ot(R_3\ot R_4))$}
\picline{12 8.5}{16 6.5}
\picline{8 8.5}{4 6.5}
\picline{6.0 3.0}{4 5.5}
\picline{14.0 3.0}{16 5.5}
\picline{8.2 2.5}{11.8 2.5}
\picputtext[m]{10 6.0}{\begin{tabular}{c}``the pentagon''\end{tabular}}
}
}\hspace{1cm}
\]
are to be commutative. A triple $(A,\Delta,\Phi)$ for which these
conditions are satisfied is called ``a quasi-Hopf algebra''.

Similarly, one may assume a relaxed form of commutativity for $\ot$,
introduce an $R\in A^{\ot 2}$, and see what $R$ has to satisfy for the
hexagons to hold. The resulting gadget $(A,\Delta,\Phi,R)$ is called
``a quasitriangular quasi-Hopf algebra''.  When the conditions on $R\in
A^{\ot 2}$ and $\Phi\in A^{\ot 3}$ are written explicitly, they are
{\em formally} identical to the pentagon and the hexagon equations that
we wrote.

Now, back to our construction.  We need to find $R$ and $\Phi$ satisfying
$\pentagon$ and $\hexagons$.

\subsection{Constructing a pair $(R,\Phi)$.} \label{solve}

\myitemm{\underline{Idea:}\quad} Use the grading of chord
diagrams. Take $R_1=1+\frac{1}{2}\diag{1em}{1}{1}
{\picvecline{0 0}{0 1}
\picvecline{1 0}{1 1}
\picline{0 0.5}{1 0.5}
}
$ and $\Phi_1=\diag{1em}{1.5}1
{\picmultigraphics{3}{0.75 0}{\picvecline{0 0}{0 1}}
}
$, and work inductively, degree by degree, to find $R$ and $\Phi$.

Assume, $(R_m,\Phi_m)$ satisfy $\pentagon$ and $\hexagons$ up to degree
$m$, and let $\mu$ and $\psi_\pm$ be the corresponding error in degree
$m+1$ of putting $(R_m,\Phi_m)$ into $\pentagon$ and $\hexagons$ modulo
degree $m+1$.  Set $\Phi_{m+1}\,=\,\Phi_m+\phi$ and $R_{m+1}=R_m+r$.

We need to solve the two equations
\begin{equation} \label{NeedToSolve1}
  \mu = \phi^{234}-(\Delta\ot1\ot1)\phi+(1\ot\Delta\ot1)\phi-
    (1\ot1\ot\Delta)\phi+\phi^{123}
\end{equation}
and
\begin{equation} \label{NeedToSolve2}
  \psi_\pm =
    \phi^{123}-\phi^{132}+\phi^{312}\pm(r^{23}-(\Delta\ot1)r+r^{13}),
\end{equation}
which are the linearizations of $\pentagon$ and of $\hexagons$.

Notice that the first equation is $\mu=d\phi$ for the differential
\[
  d\,=\,\sum_{i=0}^{n+1}\,(-1)^{i}\,d^n_i\,:\,\cA(\nup)\entspr
  \left\{
    \vcbox{\kern0.6ex\hbox{\,%
     $\displaystyle 
     \underbrace{
      \diag{4mm}{3}{2}{
        \picmultigraphics{4}{1 0}{ \picvecline{0 0}{0 2} }
        \picline{0 1.3}{1 1.3}\picline{2 1.2}{3 1.2}\picline{2 0.7}{1 0.7}
      }%
    }_n%
    $\,}\kern0.6ex%
   }
  \right\}\bigg/\,4T\,\longrightarrow\, \cA((n+1)\!\uparrow)\entspr
  \left\{
   \vcbox{\kern0.6ex\hbox{\,%
    $\displaystyle
    \underbrace{
      \diag{4mm}{4}{2}{
        \picmultigraphics{5}{1 0}{ \picvecline{0 0}{0 2} }
        \picline{0 1.3 }{1 1.3} \picline{2 1.2 }{3 1.2}
        \picline{2 0.7 }{1 0.7} \picline{3 0.5 }{4 0.5}
      }%
    }_{n+1}%
    $\,}\kern0.6ex%
   }
  \right\}\bigg/\,4T
\]

defined by

\[
  d^n_{0}(\diag{1em}{1}{1}{\picbox{0.5 0.5}{1.0 1.0}{D}})
  = \diag{1em}{1.5}{1}{\picbox{1.0 0.5}{1.0 1.0}{D}
    \picvecline{0 0}{0 1}},
  \qquad
  d^n_{n+1}(\diag{1em}{1}{1}{\picbox{0.5 0.5}{1.0 1.0}{D}})
  = \diag{1em}{1.5}{1}{\picbox{0.5 0.5}{1.0 1.0}{D}
    \picvecline{1.5 0}{1.5 1}}
\]
\[ d^n_{i}(D) = 
  (1\ot\dots\ot1\ot
  \vbox{\hbox{$\;\scriptsize i$}\hbox{$\Delta$}}
  \ot1\ot\dots\ot1)(D) \quad 1 \le i \le n\,\qquad
  \left(\vctext{@{\,}c@{\,}}{double\\$i$'th strand}{\small}\right)\,, 
\]

If we are to have any hope of solving~\eqref{NeedToSolve1}
and~\eqref{NeedToSolve2}, we must find
relations between $\mu$ and $\psi_\pm$! In particular, we'd better be
able to prove that $d\mu=0$.

Idea: One and the same morphism, say
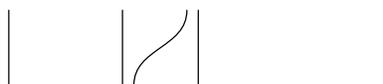
\[
\diag{5mm}{10}{2}{
\picscale{1 0.4}{
  \picline{0 0}{0 5}
  \picline{3 0}{3 5}
  \picline{5 0}{5 5}
  \picline{10 0}{10 5}
  \piccurve{3.3 0}{3.3 2.5}{4.7 2.5}{4.7 5}
 }
}
\]
appears in more than one variant of the pentagon, in which strands have
been doubled or added on the left or on the right (and also im a few
variants of the hexagons). So start from a schematic form of (say) the
hexagon, in which we put tildes on top of letters instead of bothering to
put all superscripts and $\Delta$-symbols in place:
\[ \tilde{\Phi}\tilde{R}\tilde{\Phi}\tilde{R}\tilde{\Phi}\tilde{R}
  = I+\tilde{\psi}.
\]
expand one of the $\tilde{\Phi}$'s on the left hand side
using (say) a $\pentagon$ and add an error term on the right:
\[
  \tilde{\Phi}\tilde{R}\tilde{\Phi}\tilde{\Phi}\tilde{\Phi}\tilde{\Phi}
  \tilde{R}\tilde{\Phi}\tilde{R} = I+\tilde{\psi}+\tilde{\mu}
\]
Keep going this way while simplifying whenever you can, using some more
variants of the pentagons and the hexagons, at the cost of some more
error terms on the right, or using locality relations at no cost at all.
If you're lucky, you can cancel all factors on the left (in a different
way than you have expanded them), and you get to something like this:
\[
I\enspace=\enspace I+\phi^{234}\pm (1\ot\Delta)\psi\pm\,\dots
\]
or
\[
0\enspace=\enspace \phi^{234}\pm (1\ot\Delta)\psi\pm\,\dots
\]
which is a relation of the kind we wanted.

There better be a systematic way of doing that! Here it is:

Let $CA_n$ (the $n$'th Commuto-Associahedron) be the two
dimensional CW complex made of the following cells:
\begin{list}{}{\setlength{\leftmargin}{16pt}
        \setlength{\labelwidth}{12pt}
        \setlength{\labelsep}{4pt}}

\item[\em 0-cells:] All possible products of $n$ elements
$a_1,\ldots,a_n$ in a non associative non commutative algebra.

\item[\em 1-cells:] All basic associativities and commutativities between
such products.

\item[\em 2-cells:]
\begin{list}{--}{\setlength{\leftmargin}{28pt}%
        \setlength{\labelwidth}{4pt}%
        \setlength{\labelsep}{4pt}}%
\item[] \leavevmode\kern-9.2mm --\kern0.95ex
% please don't make the following line longer than 1 LaTeX line
% else I would have more problems with this hack
% there is surely a better solution but I didn't want to
% bother with the macros themselves
        A pentagon sealing every pentagon of the kind appearing in the
	pentagon relation.
\item A hexagon sealing every hexagon of the kind appearing in
	the hexagon relation.
\item A square sealing every locality relation of the type
        considered in section~\ref{FreeRelations}. Notice that every
        locality relation can be written as a product of four $R$'s and
        $\Phi$'s, and so it corresponds to a square in the 1-skeleton
        of $C\!A_n$.
\end{list}
\end{list}
\def\@temp@{
\rbraid{0.3 3.25}{0.6 0.75}
\picmultigraphics{2}{0 1.125}{
  \picmultigraphics{2}{0.6 0}{
    \picline{0 2.5}{0 2.875}
  }
}}
\myitem[6.4cm]{\kern2cm\vbox to 1em{\kern6mm\hbox{\diag{4mm}{2}{8}{
\picscale{2 3 : 1}{
\pic@pss{/default@brdrm [-7 1 -1 0] D}
\picscale{1 1}{\@temp@
}
\picline{4 2.5}{4 4}
\pictranslate{4 4}{
  \picscale{-1 1}{
    \@temp@
    %\rbraid{0.15 3.25}{0.6 0.75}
}}
\piccurve{3.4 0}{3.4 1.1}{0 1.1}{0 2.5}
\piccurve{4.0 0}{4.0 1.4}{0.6 1.4}{0.6 2.5}
\picmulticurve{-7 1 -1 0}{0 0}{0 1.3}{4 1.3}{4 2.5}
\pictranslate{4 4}{
  \picscale{-1 1}{
    \piccurve{3.4 0}{3.4 1.1}{0 1.1}{0 2.5}
    \piccurve{4.0 0}{4.0 1.4}{0.6 1.4}{0.6 2.5}
    \picmulticurve{-7 1 -1 0}{0 0}{0 1.3}{4 1.3}{4 2.5}
  }
}
    \picline{4 2.5}{4 4}
    \picline{0 6.5}{0 8}
\picline{-0.3 2.5}{4.3 2.5}
\picline{-0.3 4.0}{4.3 4.0}
\picline{-0.3 6.5}{4.3 6.5}
}}
\raisebox{-1.8mm}{$\hskip1cm\rightarrowto{8mm}$}}\vss}}
\null
\[
(1\ot\Delta)R\cdot R^{23}\cdot (1\ot\Delta)(R^{-1})\cdot
(R^{23})^{-1}\,=\,1
\]
or
\[
\diag{5mm}{8}{5}{
%\picfillgraycol{0.5 0.3 0.7}
\picshade{\picbox{4 2.5}{5 4}{}}{-30}{-16}
\picbox{4 2.5}{5 4}{}
%\picgraycol{1}
\picfillgraycol{1}
\picnolines
\picfilledbox{1.5 0.5}{2.5 1}{$(BC)A$}
\picfilledbox{6.5 0.5}{2.5 1}{$(CB)A$}
\picfilledbox{1.5 4.5}{2.5 1}{$A(BC)$}
\picfilledbox{6.5 4.5}{2.5 1}{$A(CB)$}
}
\]\vspace{3mm}

\noindent For an example, see figure~\ref{CA3}.

\begin{figure}[htpb]
\[ 
        \setlength{\unitlength}{0.525\eepiclength}
        \begin{array}{c}  \hspace{-1.7mm}
                \raisebox{-8pt}{\begingroup\makeatletter\ifx\SetFigFont\undefined%
\gdef\SetFigFont#1#2#3#4#5{%
  \reset@font\fontsize{#1}{#2pt}%
  \fontfamily{#3}\fontseries{#4}\fontshape{#5}%
  \selectfont}%
\fi\endgroup%
\begin{picture}(5473,4770)(0,-10)
\thicklines
\put(2175,4461){\blacken\ellipse{90}{90}}
\put(2175,4461){\ellipse{90}{90}}
\put(1275,3861){\blacken\ellipse{90}{90}}
\put(1275,3861){\ellipse{90}{90}}
\put(675,2961){\blacken\ellipse{90}{90}}
\put(675,2961){\ellipse{90}{90}}
\put(3375,4461){\blacken\ellipse{90}{90}}
\put(3375,4461){\ellipse{90}{90}}
\put(675,1761){\blacken\ellipse{90}{90}}
\put(675,1761){\ellipse{90}{90}}
\put(1275,861){\blacken\ellipse{90}{90}}
\put(1275,861){\ellipse{90}{90}}
\put(2175,261){\blacken\ellipse{90}{90}}
\put(2175,261){\ellipse{90}{90}}
\put(3375,261){\blacken\ellipse{90}{90}}
\put(3375,261){\ellipse{90}{90}}
\put(4275,861){\blacken\ellipse{90}{90}}
\put(4275,861){\ellipse{90}{90}}
\put(4875,1761){\blacken\ellipse{90}{90}}
\put(4875,1761){\ellipse{90}{90}}
\put(4875,2961){\blacken\ellipse{90}{90}}
\put(4875,2961){\ellipse{90}{90}}
\put(4275,3861){\blacken\ellipse{90}{90}}
\put(4275,3861){\ellipse{90}{90}}
\texture{cccccccc 0 0 0 cccccccc 0 0 0 
	cccccccc 0 0 0 cccccccc 0 0 0 
	cccccccc 0 0 0 cccccccc 0 0 0 
	cccccccc 0 0 0 cccccccc 0 0 0 }
\shade\path(3375,261)(2175,261)(2175,4461)
	(3375,4461)(3375,261)
\path(3375,261)(2175,261)(2175,4461)
	(3375,4461)(3375,261)
\texture{c0c0c0c0 0 0 0 0 0 0 0 
	c0c0c0c0 0 0 0 0 0 0 0 
	c0c0c0c0 0 0 0 0 0 0 0 
	c0c0c0c0 0 0 0 0 0 0 0 }
\shade\path(2175,4461)(2175,261)(1275,861)
	(675,1761)(675,2961)(1275,3861)(2175,4461)
\path(2175,4461)(2175,261)(1275,861)
	(675,1761)(675,2961)(1275,3861)(2175,4461)
\shade\path(3375,4461)(3375,261)(4275,861)
	(4875,1761)(4875,2961)(4275,3861)(3375,4461)
\path(3375,4461)(3375,261)(4275,861)
	(4875,1761)(4875,2961)(4275,3861)(3375,4461)
\path(1275,3861)(4875,1761)
\path(675,2961)(4275,861)
\path(4275,3861)(675,1761)
\path(4875,2961)(1275,861)
\put(4350,3861){\makebox(0,0)[lb]{\smash{{\mathmode{\scriptstyle c(ab)}}}}}
\put(4950,1611){\makebox(0,0)[lb]{\smash{{\mathmode{\scriptstyle (cb)a}}}}}
\put(1875,36){\makebox(0,0)[lb]{\smash{{\mathmode{\scriptstyle b(ac)}}}}}
\put(3225,36){\makebox(0,0)[lb]{\smash{{\mathmode{\scriptstyle b(ca)}}}}}
\put(4950,2886){\makebox(0,0)[lb]{\smash{{\mathmode{\scriptstyle c(ba)}}}}}
\put(3150,4611){\makebox(0,0)[lb]{\smash{{\mathmode{\scriptstyle (ca)b}}}}}
\put(1875,4611){\makebox(0,0)[lb]{\smash{{\mathmode{\scriptstyle (ac)b}}}}}
\put(600,3861){\makebox(0,0)[lb]{\smash{{\mathmode{\scriptstyle a(cb)}}}}}
\put(4350,636){\makebox(0,0)[lb]{\smash{{\mathmode{\scriptstyle (bc)a}}}}}
\put(0,2886){\makebox(0,0)[lb]{\smash{{\mathmode{\scriptstyle a(bc)}}}}}
\put(0,1611){\makebox(0,0)[lb]{\smash{{\mathmode{\scriptstyle (ab)c}}}}}
\put(600,711){\makebox(0,0)[lb]{\smash{{\mathmode{\scriptstyle (ba)c}}}}}
\end{picture} }
                \hspace{-1.9mm}
        \end{array}

\parbox{3.6in}{
	\caption{The third commuto-associahedron $C\!A_3$. It is made
	by gluing three sets of a square and two hexagons each into
	a 12-gon.  Only one of these sets is shaded in the figure; the
	other two are obtained from it by rotations by $60^\circ$ and
	$120^\circ$ respectively. Topologically, the result is a circle
	with three disks glued in, and has the homotopy type of a wedge
	of two spheres.}
\label{CA3} }
\]
\end{figure}

The faces corresponding to rules applied for obtaining a relation
(as described above) form a subcomplex of $CA_n$ homeomorphic to a
closed surface.  If we note each of the error terms in $\pentagon$
and $\hexagons$ on the corresponding faces in $CA_n$ (for locality
relation this error is 0) this relation says that the sum of the terms
on the faces of this surface vanishes.  So, to find out (maximally)
how many (independent) relations we may get, we need to know $b_2=\dim
H_2(CA_n)$. MacLane's coherence theorem says that $CA_n$ is simply
connected, so we are left with a simple counting for determining
$\chi(CA_n)$. We find that
\[
b_2\quad=\quad\#\{\text{ vertices }\}\,-\,
\#\{\text{ edges }\}\,+
\#\{\text{ faces }\}\,-\,1
\]

\begin{example}
Let $n=3$. $CA_3$ is a circle with 3 disks attached (each of
them separated into a rectangle and two hexagons), as shown in
figure~\ref{CA3}.  So
\[
b_2(CA_3)\quad=\quad 12\,-\,18\,+\,9\,-\,1\quad=\quad2\,,
\]
and we can hope to find 2 relations. They turn out to be
\begin{equation} \label{P12eq1}
  \psi^{123}-\psi^{132}+\psi^{213}-\psi^{231}=0
\end{equation}
and
\begin{equation} \label{P12eq2}
  \psi^{213}-\psi^{231}+\psi^{312}-\psi^{321}=0.
\end{equation}
\end{example}

Similarly, the polyhedra in figures \ref{CA4}, \ref{CA44}, and~\ref{CA5}
prove equation~\ref{P20eq}, \ref{P32eq}, and~\ref{L5eq} respectively.

\begin{equation} \label{P20eq}
        \mu^{1234}-\mu^{1243}+\mu^{1423}-\mu^{4123}
        =\psi_{+}^{234}-(\Delta\otimes 1\otimes 1)\psi_{+}
        +(1\otimes\Delta\otimes 1)\psi_{+}-\psi_{+}^{124}
\end{equation}

\def\narrowot{\!\otimes\!}
\begin{equation} \label{P32eq}
       \mu^{1234}\!-\!\mu^{1324}\!+\!\mu^{1342}
               \!+\!\mu^{3124}\!-\!\mu^{3142}\!+\!\mu^{3412}
       \!=\!(1\narrowot 1\narrowot\Delta)\psi_{+}
               \!-\!\psi_{+}^{124} \!-\! \psi_{+}^{123}
       \!-\!(\Delta\narrowot 1\narrowot 1)\psi_{-}
               \!+\!\psi_{-}^{234} \!+\! \psi_{-}^{134}
\end{equation}

\begin{equation} \label{L5eq}
       \mu^{2345}\!-\!(\Delta\narrowot 1\narrowot 1\narrowot 1)(\mu)
       \!+\!(1\narrowot\Delta\narrowot 1\narrowot 1)(\mu)
       \!-\!(1\narrowot 1\narrowot\Delta\narrowot 1)(\mu)
       \!+\!(1\narrowot 1\narrowot 1\narrowot\Delta)(\mu)
       \!-\!\mu^{1234} \!=\! 0
\end{equation}

{
  \def\Tdii{(\Delta 11)}
  \def\Tidi{(1\Delta 1)}
  \def\Tiid{(11\Delta)}
  \def\Tddi{(\Delta^{\uparrow\uparrow\uparrow} 1)}
  \def\Tdi{(\Delta 1)}
  \def\Tdiii{(\Delta 111)}
  \def\Tidii{(1\Delta 11)}
  \def\Tiidi{(11\Delta 1)}
  \def\Tiiid{(111\Delta)}

  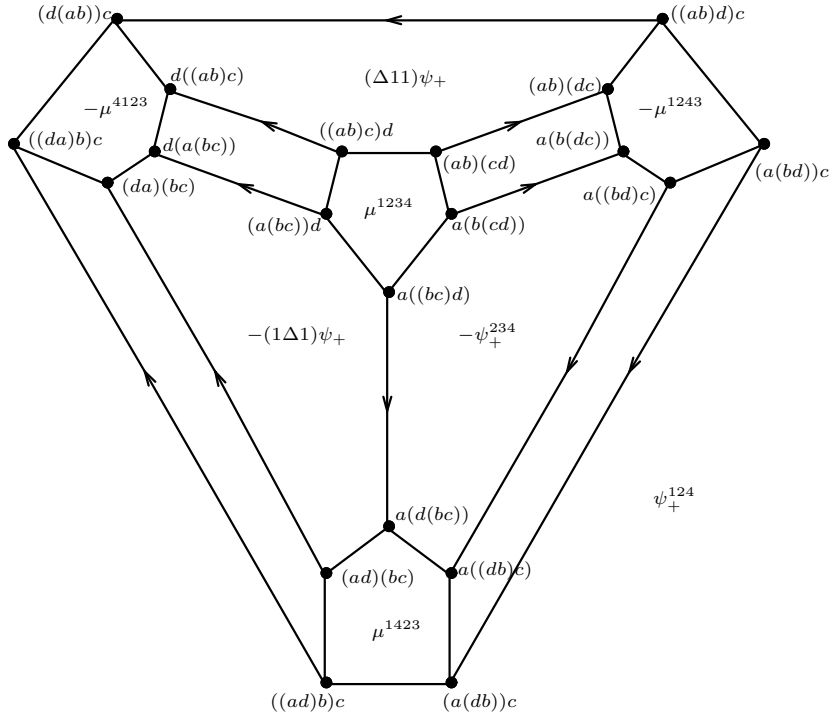
\begin{figure}[htpb]
  \[ \qquad
        \setlength{\unitlength}{0.65\eepiclength}
        \begin{array}{c}  \hspace{-1.7mm}
                \raisebox{-8pt}{\begingroup\makeatletter\ifx\SetFigFont\undefined%
\gdef\SetFigFont#1#2#3#4#5{%
  \reset@font\fontsize{#1}{#2pt}%
  \fontfamily{#3}\fontseries{#4}\fontshape{#5}%
  \selectfont}%
\fi\endgroup%
\begin{picture}(7916,6810)(0,-10)
\thicklines
\put(53,5451){\blacken\ellipse{90}{90}}
\put(53,5451){\ellipse{90}{90}}
\put(1043,6651){\blacken\ellipse{90}{90}}
\put(1043,6651){\ellipse{90}{90}}
\put(1553,5976){\blacken\ellipse{90}{90}}
\put(1553,5976){\ellipse{90}{90}}
\put(1403,5376){\blacken\ellipse{90}{90}}
\put(1403,5376){\ellipse{90}{90}}
\put(953,5076){\blacken\ellipse{90}{90}}
\put(953,5076){\ellipse{90}{90}}
\put(3203,5376){\blacken\ellipse{90}{90}}
\put(3203,5376){\ellipse{90}{90}}
\put(3053,4776){\blacken\ellipse{90}{90}}
\put(3053,4776){\ellipse{90}{90}}
\put(3653,4026){\blacken\ellipse{90}{90}}
\put(3653,4026){\ellipse{90}{90}}
\put(4253,4776){\blacken\ellipse{90}{90}}
\put(4253,4776){\ellipse{90}{90}}
\put(4103,5376){\blacken\ellipse{90}{90}}
\put(4103,5376){\ellipse{90}{90}}
\put(5903,5376){\blacken\ellipse{90}{90}}
\put(5903,5376){\ellipse{90}{90}}
\put(5753,5976){\blacken\ellipse{90}{90}}
\put(5753,5976){\ellipse{90}{90}}
\put(6278,6651){\blacken\ellipse{90}{90}}
\put(6278,6651){\ellipse{90}{90}}
\put(7253,5451){\blacken\ellipse{90}{90}}
\put(7253,5451){\ellipse{90}{90}}
\put(6353,5076){\blacken\ellipse{90}{90}}
\put(6353,5076){\ellipse{90}{90}}
\put(3653,1776){\blacken\ellipse{90}{90}}
\put(3653,1776){\ellipse{90}{90}}
\put(3053,1326){\blacken\ellipse{90}{90}}
\put(3053,1326){\ellipse{90}{90}}
\put(4253,276){\blacken\ellipse{90}{90}}
\put(4253,276){\ellipse{90}{90}}
\put(4253,1326){\blacken\ellipse{90}{90}}
\put(4253,1326){\ellipse{90}{90}}
\put(3053,276){\blacken\ellipse{90}{90}}
\put(3053,276){\ellipse{90}{90}}
\path(6263,6636)(3638,6636)
\path(3773.000,6673.500)(3638.000,6636.000)(3773.000,6598.500)
\path(7238,5436)(5963,3261)
\path(5998.921,3396.429)(5963.000,3261.000)(6063.623,3358.500)
\path(6338,5061)(5363,3261)
\path(5394.325,3397.565)(5363.000,3261.000)(5460.272,3361.844)
\path(3038,1311)(1988,3261)
\path(2085.021,3159.915)(1988.000,3261.000)(2018.986,3124.358)
\path(5363,3261)(4238,1311)
\path(4088,5361)(4913,5661)
\path(4798.943,5579.622)(4913.000,5661.000)(4773.313,5650.107)
\path(4238,4761)(5063,5061)
\path(4948.943,4979.622)(5063.000,5061.000)(4923.313,5050.107)
\path(3038,4761)(2213,5061)
\path(2352.687,5050.107)(2213.000,5061.000)(2327.057,4979.622)
\path(3188,5361)(2438,5661)
\path(2577.272,5645.680)(2438.000,5661.000)(2549.417,5576.044)
\path(5738,5961)(4913,5661)
\path(5063,5061)(5888,5361)
\path(2438,5661)(1538,5961)
\path(1388,5361)(2213,5061)
\path(1988,3261)(938,5061)
\path(3638,2886)(3638,1761)
\path(3638,4011)(3638,2886)
\path(3600.500,3021.000)(3638.000,2886.000)(3675.500,3021.000)
\path(1313,3261)(38,5436)
\path(5963,3261)(4238,261)
\path(3038,261)(1313,3261)
\path(1412.803,3162.660)(1313.000,3261.000)(1347.785,3125.275)
\path(1013,6636)(3638,6636)
\path(1013,6636)(38,5436)(938,5061)
	(1388,5361)(1538,5961)(1013,6636)
\path(3188,5361)(4088,5361)(4238,4761)
	(3638,4011)(3038,4761)(3188,5361)
\path(6263,6636)(5738,5961)(5888,5361)
	(6338,5061)(7238,5436)(6263,6636)
\path(3638,1761)(3038,1311)(3038,261)
	(4238,261)(4238,1311)(3638,1761)(3638,1761)
\put(3728,3951){\makebox(0,0)[lb]{\smash{{\mathmode{\scriptstyle a((bc)d)}}}}}
\put(2978,5526){\makebox(0,0)[lb]{\smash{{\mathmode{\scriptstyle ((ab)c)d}}}}}
\put(6353,6651){\makebox(0,0)[lb]{\smash{{\mathmode{\scriptstyle ((ab)d)c}}}}}
\put(1478,5376){\makebox(0,0)[lb]{\smash{{\mathmode{\scriptstyle d(a(bc))}}}}}
\put(4253,4626){\makebox(0,0)[lb]{\smash{{\mathmode{\scriptstyle a(b(cd))}}}}}
\put(2513,36){\makebox(0,0)[lb]{\smash{{\mathmode{\scriptstyle ((ad)b)c}}}}}
\put(4163,36){\makebox(0,0)[lb]{\smash{{\mathmode{\scriptstyle (a(db))c}}}}}
\put(278,6651){\makebox(0,0)[lb]{\smash{{\mathmode{\scriptstyle (d(ab))c}}}}}
\put(1553,6051){\makebox(0,0)[lb]{\smash{{\mathmode{\scriptstyle d((ab)c)}}}}}
\put(713,5736){\makebox(0,0)[lb]{\smash{{\mathmode{\scriptstyle -\mu^{4123}}}}}}
\put(3488,711){\makebox(0,0)[lb]{\smash{{\mathmode{\scriptstyle \mu^{1423}}}}}}
\put(3413,4761){\makebox(0,0)[lb]{\smash{{\mathmode{\scriptstyle \mu^{1234}}}}}}
\put(4313,3561){\makebox(0,0)[lb]{\smash{{\mathmode{\scriptstyle -\psi_+^{234}}}}}}
\put(3413,6036){\makebox(0,0)[lb]{\smash{{\mathmode{\scriptstyle \Tdii\psi_+}}}}}
\put(6188,1986){\makebox(0,0)[lb]{\smash{{\mathmode{\scriptstyle \psi_+^{124}}}}}}
\put(188,5436){\makebox(0,0)[lb]{\smash{{\mathmode{\scriptstyle ((da)b)c}}}}}
\put(1088,4986){\makebox(0,0)[lb]{\smash{{\mathmode{\scriptstyle (da)(bc)}}}}}
\put(2288,4611){\makebox(0,0)[lb]{\smash{{\mathmode{\scriptstyle (a(bc))d}}}}}
\put(4163,5211){\makebox(0,0)[lb]{\smash{{\mathmode{\scriptstyle (ab)(cd)}}}}}
\put(4988,5961){\makebox(0,0)[lb]{\smash{{\mathmode{\scriptstyle (ab)(dc)}}}}}
\put(5513,4911){\makebox(0,0)[lb]{\smash{{\mathmode{\scriptstyle a((bd)c)}}}}}
\put(7163,5136){\makebox(0,0)[lb]{\smash{{\mathmode{\scriptstyle (a(bd))c}}}}}
\put(3188,1236){\makebox(0,0)[lb]{\smash{{\mathmode{\scriptstyle (ad)(bc)}}}}}
\put(3713,1836){\makebox(0,0)[lb]{\smash{{\mathmode{\scriptstyle a(d(bc))}}}}}
\put(5063,5436){\makebox(0,0)[lb]{\smash{{\mathmode{\scriptstyle a(b(dc))}}}}}
\put(4313,1311){\makebox(0,0)[lb]{\smash{{\mathmode{\scriptstyle a((db)c)}}}}}
\put(2288,3561){\makebox(0,0)[lb]{\smash{{\mathmode{\scriptstyle -\Tidi\psi_+}}}}}
\put(6038,5736){\makebox(0,0)[lb]{\smash{{\mathmode{\scriptstyle -\mu^{1243}}}}}}
\end{picture} }
                \hspace{-1.9mm}
        \end{array}
 \]
  \caption{\label{CA4}A relation from $CA_4$}
  \end{figure}

  \begin{figure}[htpb]
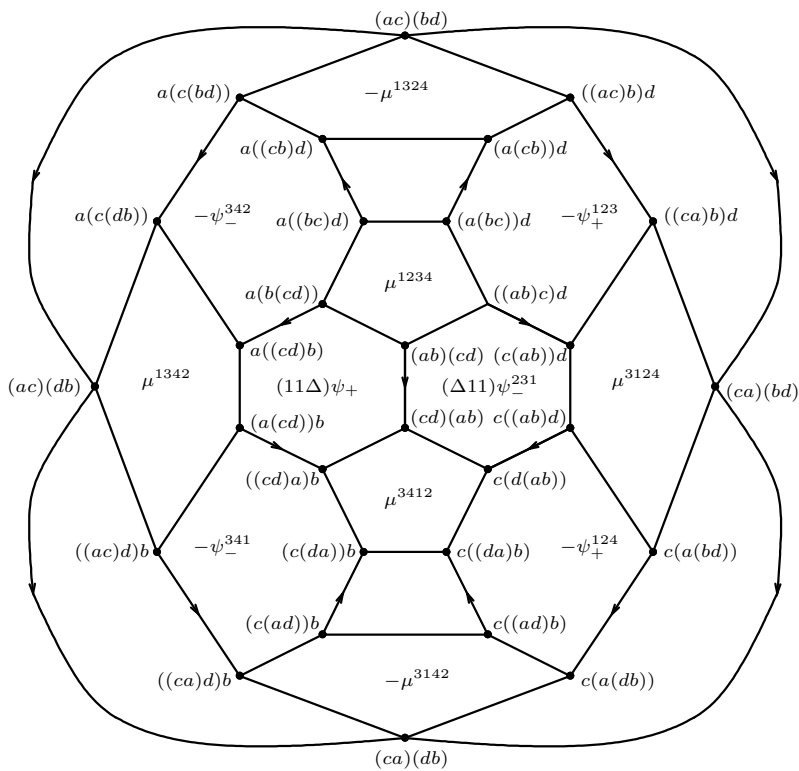

  \[ 
        \setlength{\unitlength}{0.43\eepiclength}
        \begin{array}{c}  \hspace{-1.7mm}
                \raisebox{-8pt}{\input draws/P32.tex }
                \hspace{-1.9mm}
        \end{array}
 \]
  \caption{\label{CA44}Another relation from $CA_4$}
  \end{figure}

  \begin{figure}[htpb]
  \[ 
        \setlength{\unitlength}{0.6\eepiclength}
        \begin{array}{c}  \hspace{-1.7mm}
                \raisebox{-8pt}{\begingroup\makeatletter\ifx\SetFigFont\undefined%
\gdef\SetFigFont#1#2#3#4#5{%
  \reset@font\fontsize{#1}{#2pt}%
  \fontfamily{#3}\fontseries{#4}\fontshape{#5}%
  \selectfont}%
\fi\endgroup%
\begin{picture}(6849,4095)(0,-10)
\thicklines
\put(450,3861){\blacken\ellipse{90}{90}}
\put(450,3861){\ellipse{90}{90}}
\put(450,261){\blacken\ellipse{90}{90}}
\put(450,261){\ellipse{90}{90}}
\put(6450,261){\blacken\ellipse{90}{90}}
\put(6450,261){\ellipse{90}{90}}
\put(6450,3861){\blacken\ellipse{90}{90}}
\put(6450,3861){\ellipse{90}{90}}
\put(5670,3111){\blacken\ellipse{90}{90}}
\put(5670,3111){\ellipse{90}{90}}
\put(5670,1011){\blacken\ellipse{90}{90}}
\put(5670,1011){\ellipse{90}{90}}
\put(1200,1011){\blacken\ellipse{90}{90}}
\put(1200,1011){\ellipse{90}{90}}
\put(1200,3111){\blacken\ellipse{90}{90}}
\put(1200,3111){\ellipse{90}{90}}
\put(1500,2061){\blacken\ellipse{90}{90}}
\put(1500,2061){\ellipse{90}{90}}
\put(5400,2061){\blacken\ellipse{90}{90}}
\put(5400,2061){\ellipse{90}{90}}
\put(3450,1311){\blacken\ellipse{90}{90}}
\put(3450,1311){\ellipse{90}{90}}
\put(3450,2811){\blacken\ellipse{90}{90}}
\put(3450,2811){\ellipse{90}{90}}
\put(2850,1911){\blacken\ellipse{90}{90}}
\put(2850,1911){\ellipse{90}{90}}
\put(4050,2211){\blacken\ellipse{90}{90}}
\put(4050,2211){\ellipse{90}{90}}
\path(6450,261)(6450,3861)(450,3861)
	(450,261)(6450,261)
\path(450,3861)(1200,3111)
\path(1200,1011)(450,261)
\path(5700,1011)(6450,261)
\path(5700,3111)(6450,3861)
\path(3450,2811)(4050,2211)(5400,2061)
\path(1500,2061)(2850,1911)(3450,1311)
\path(4050,2211)(2850,1911)
\path(1200,3111)(3450,2811)(5700,3111)
	(5400,2061)(5700,1011)(3450,1311)
	(1200,1011)(1500,2061)(1200,3111)
\put(1200,3186){\makebox(0,0)[lb]{\smash{{\mathmode{\scriptstyle a((bc)(de))}}}}}
\put(0,3936){\makebox(0,0)[lb]{\smash{{\mathmode{\scriptstyle a(b(c(de))}}}}}
\put(5850,3936){\makebox(0,0)[lb]{\smash{{\mathmode{\scriptstyle a(b((cd)e))}}}}}
\put(1575,2136){\makebox(0,0)[lb]{\smash{{\mathmode{\scriptstyle (a(bc))(de)}}}}}
\put(4125,2286){\makebox(0,0)[lb]{\smash{{\mathmode{\scriptstyle (a((bc)d))e}}}}}
\put(5850,36){\makebox(0,0)[lb]{\smash{{\mathmode{\scriptstyle (ab)((cd)e)}}}}}
\put(0,36){\makebox(0,0)[lb]{\smash{{\mathmode{\scriptstyle (ab)(c(de))}}}}}
\put(4725,3261){\makebox(0,0)[lb]{\smash{{\mathmode{\scriptstyle a((b(cd))e)}}}}}
\put(1200,786){\makebox(0,0)[lb]{\smash{{\mathmode{\scriptstyle ((ab)c)(de)}}}}}
\put(3675,2661){\makebox(0,0)[lb]{\smash{{\mathmode{\scriptstyle a(((bc)d)e)}}}}}
\put(1950,2511){\makebox(0,0)[lb]{\smash{{\mathmode{\scriptstyle \Tidii\mu}}}}}
\put(3900,1461){\makebox(0,0)[lb]{\smash{{\mathmode{\scriptstyle -\mu^{1234}}}}}}
\put(525,1686){\makebox(0,0)[lb]{\smash{{\mathmode{\scriptstyle \Tiiid\mu}}}}}
\put(3000,636){\makebox(0,0)[lb]{\smash{{\mathmode{\scriptstyle -\Tdiii\mu}}}}}
\put(4800,786){\makebox(0,0)[lb]{\smash{{\mathmode{\scriptstyle ((ab)(cd))e}}}}}
\put(3225,3336){\makebox(0,0)[lb]{\smash{{\mathmode{\scriptstyle \mu^{2345}}}}}}
\put(5475,1986){\makebox(0,0)[lb]{\smash{{\mathmode{\scriptstyle -\Tiidi\mu}}}}}
\put(4350,1836){\makebox(0,0)[lb]{\smash{{\mathmode{\scriptstyle (a(b(cd)))e}}}}}
\put(2175,1386){\makebox(0,0)[lb]{\smash{{\mathmode{\scriptstyle (((ab)c)d)e}}}}}
\put(1800,1761){\makebox(0,0)[lb]{\smash{{\mathmode{\scriptstyle ((a(bc))d)e}}}}}
\end{picture} }
                \hspace{-1.9mm}
        \end{array}
 \]
  \caption{The Stasheff polyhedron --- a relation from $CA_5$ \label{CA5}}
  \end{figure}
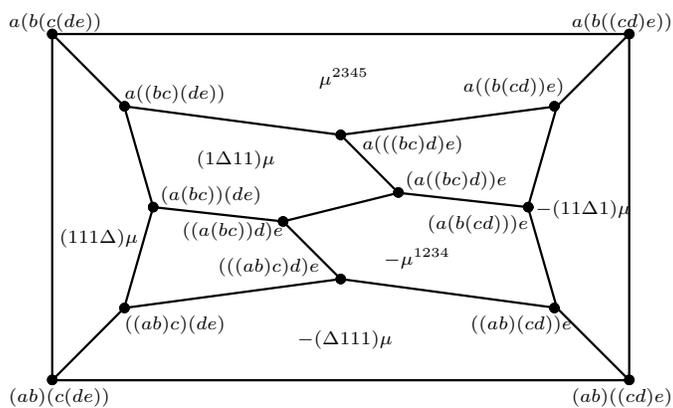
}

Notice that equation~\eqref{L5eq} simply says that $d\mu=0$. For
simplicity, let's pretend now that only the pentagon has to be solved, and
that only equation~\eqref{L5eq} is given. In reality we also need to
solve~\eqref{NeedToSolve2}, and we're also given equations~\eqref{P12eq1},
\eqref{P12eq2}, \eqref{P20eq}, and~\eqref{P32eq}. This additional
requirement and that additional information makes matters more complicated,
but the principles remain the same. Anyway, with our simplifying
assumptions, equations~\eqref{NeedToSolve1} and~\eqref{L5eq} together mean that
we're left with showing that $H^4(\cA(\nup))=0$.

Without our simplifying assumption we end up needing to show that some
easily defined subcomplex of $\cA(\nup)$ has vanishing cohomology,
$H^4_{\text{\it sub}}(\cA(\nup))=0$.

There are two possible interpretations for $\cA(\nup)$~---
\begin{itemize}
\item[*] allowing non-horizontal chords:\hfill\mbox{\rule[-2mm]{0mm}{2mm}}\\
In this case it is known that $H^4_{\text{\it sub}}(\cA(\nup))\,=\,0$, but explicit
computations are almost impossible.
\item[*] allowing only horizontal chords:\hfill\mbox{\rule[-2mm]{0mm}{2mm}}\\
In this case explicit computations are easy, but we don't know
how to compute $H^4_{\text{\it sub}}(\cA(\nup))$.
See~\cite{Stoimenow:HarrisonKohomologie} for a partial result on this
problem related to some combinatorial properties of a free resolution of
$\cA(\nup)$.
\end{itemize}

\subsection{Why are we not happy?}

\begin{itemize}
\item Why is it that we can compute $H^4_{\text{\it sub}}(\cA(\nup))$ only in
  the less natural case in which non-horizontal chords are allowed? We know
  that a horizontal-chord-only $\Phi$ does exist; Drinfel'd constructed
  one using the KZ connection in~\cite{Drinfeld:QuasiHopf}. But we still
  don't have a proof of this fact that does not use analysis.
\item The algorithm we sketched here finds a pair $(R,\Phi)$. From the
  considerations in section~\ref{motivation} (and
  from~\cite{Drinfeld:QuasiHopf}) we know that we should be able to
  take $R=\,\exp\left(\frac{1}{2}\diag{1.5em}{1}{1} {\picvecline{0 0}{0
  1} \picvecline{1 0}{1 1} \picline{0 0.5}{1 0.5} }\,\right)$. But we
  don't know how to reproduce this fact algebraically.
  % isn't "algebraically" better than "algebraicly" ?
\item And we still don't know anything about $\bZ/3\bZ$ and many other
  rings.
\end{itemize}

\ifx\undefined\bysame
        \newcommand{\bysame}{\leavevmode\hbox to3em{\hrulefill}\,}
\fi

\end{document}